\documentclass[12pt,psamsfonts,leqno,oneside,letterpaper]{amsart}
\usepackage[dvips,text={6.5truein,9truein},left=1truein,top=1truein]{geometry}
\usepackage{amssymb,amsmath,amscd,enumerate,
}
\usepackage[pdftex]{graphicx}
\usepackage{url}

\usepackage[colorlinks,linkcolor=blue,citecolor=blue,pdfstartview=FitH]{hyperref}
\input xy
\xyoption{all}
\SelectTips{cm}{12}

\usepackage{color}
\usepackage{mathrsfs}
\newcommand\redden[1]{{\color{red}#1}}

\newcommand\invisiblecomment[1]{\empty}

\newcommand\tinymissingref[1]{\empty}
\newcommand\abstractcomment[1]{\empty}

\theoremstyle{definition}
\newtheorem{para}{}[section]

\newtheorem{remark}[para]{Remark}
\newtheorem{reformulation}[para]{Reformulation}
\newtheorem{remarks}[para]{Remarks}
\newtheorem{notation}[para]{Notation}

\newtheorem{convention}[para]{Convention}
\newtheorem{definition}[para]{Definition}
\newtheorem{squeezingtheorem}[para]{Squeezing Theorem}
\newtheorem{definitions}[para]{Definitions}
\newtheorem{definitionnotation}[para]{Definition and Notation}
\newtheorem{notationdefinitionremark}[para]{Notation, Definition and Remark}

\newtheorem{remarksnotation}[para]{Remarks and Notation}

\newtheorem{remarknotation}[para]{Remark and Notation}
\newtheorem{notationremark}[para]{Notation and Remark}
\newtheorem{notationreviewremarks}[para]{Notation, Review and Remarks}
\newtheorem{definitionremark}[para]{Definition and Remark}
\newtheorem{definitionsremarks}[para]{Definitions and Remarks}
\newtheorem{notationremarks}[para]{Notation and Remarks}
\newtheorem{definitionsnotation}[para]{Definitions and Notation}

\newtheorem{reviewdefinition}[para]{Review and Definition}
\newtheorem{definitionremarks}[para]{Definition and Remarks}
\newtheorem{definitionnotationremarks}[para]{Definition, Notation and Remarks}
\newtheorem{definitionsnotationremarks}[para]{Definitions, Notation and Remarks}

\newcommand\Alternatives{\begin{enumerate}[(i)]}
\newcommand\EndAlternatives{\end{enumerate}}
\newcommand\Conditions{\begin{enumerate}[(1)]}
\newcommand\EndConditions{\end{enumerate}}

\theoremstyle{plain}
\newtheorem{theorem}[para]{Theorem}

\newtheorem{lemma}[para]{Lemma}
\newtheorem{remarkdefinition}[para]{Remark and Definition}
\newtheorem{remarkdefinitions}[para]{Remark and Definitions}
\newtheorem{proposition}[para]{Proposition}

\newtheorem{corollary}[para]{Corollary}
\newtheorem{conjecture}[para]{Conjecture}

\newtheorem*{unnumberedremark}{Remark}

\newtheorem*{whoknowsclaim}{3.12.1}

\newtheorem*{PropositionA}{Proposition A}
\newtheorem*{PropositionB}{Proposition B}
\newtheorem*{TheoremC}{Theorem C}

\newtheorem{claim}[equation]{}
\numberwithin{equation}{para}
\numberwithin{figure}{section}
\numberwithin{specialremark}{para}
\numberwithin{specialnumber}{para}

\newcommand\Number{\begin{para}}
\newcommand\EndNumber{\end{para}}
\newcommand\Definition{\begin{definition}}
\newcommand\EndDefinition{\end{definition}}
\newcommand\Definitions{\begin{definitions}}
\newcommand\DefinitionsNotation{\begin{definitionsnotation}}
\newcommand\NotationDefinitionRemark{\begin{notationdefinitionremark}}
\newcommand\EndNotationDefinitionRemark{\end{notationdefinitionremark}}

\newcommand\DefinitionNotation{\begin{definitionnotation}}
\newcommand\RemarksNotation{\begin{remarksnotation}}
\newcommand\Reformulation{\begin{reformulation}}
\newcommand\EndRemarksNotation{\end{remarksnotation}}
\newcommand\EndReformulation{\end{reformulation}}
\newcommand\RemarkNotation{\begin{remarknotation}}
\newcommand\EndRemarkNotation{\end{remarknotation}}
\newcommand\NotationRemark{\begin{notationremark}}
\newcommand\EndDefinitionNotationRemarks{\end{definitionnotationremarks}}
\newcommand\NotationReviewRemarks{\begin{notationreviewremarks}}
\newcommand\DefinitionRemark{\begin{definitionremark}}
\newcommand\DefinitionRemarks{\begin{definitionremarks}}    
\newcommand\DefinitionsRemarks{\begin{definitionsremarks}}
\newcommand\DefinitionNotationRemarks{\begin{definitionnotationremarks}}
\newcommand\DefinitionsNotationRemarks{\begin{definitionsnotationremarks}}
\newcommand\EndDefinitionsNotationRemarks{\end{definitionsnotationremarks}}

\newcommand\NotationRemarks{\begin{notationremarks}}
\newcommand\EndNotationRemark{\end{notationremark}}
\newcommand\EndNotationReviewRemarks{\end{notationreviewremarks}}
\newcommand\EndDefinitionRemark{\end{definitionremark}}
\newcommand\EndDefinitionRemarks{\end{definitionremarks}}
\newcommand\EndDefinitionsRemarks{\end{definitionsremarks}}
\newcommand\EndNotationRemarks{\end{notationremarks}}
\newcommand\EndRemarkDefinition{\end{remarkdefinition}}
\newcommand\EndRemarkDefinitions{\end{remarkdefinitions}}
\newcommand\RemarkDefinition{\begin{remarkdefinition}}
\newcommand\RemarkDefinitions{\begin{remarkdefinitions}}
\newcommand\EndDefinitionsNotation{\end{definitionsnotation}}
\newcommand\EndDefinitionNotation{\end{definitionnotation}}
\newcommand\ReviewDefinition{\begin{reviewdefinition}}
\newcommand\EndReviewDefinition{\end{reviewdefinition}}
\newcommand\EndDefinitions{\end{definitions}}
\newcommand\Theorem{\begin{theorem}}
\newcommand\EndTheorem{\end{theorem}}
\newcommand\Conjecture{\begin{conjecture}}
\newcommand\EndConjecture{\end{conjecture}}
\newcommand\Remark{\begin{remark}}
\newcommand\EndRemark{\end{remark}}
\newcommand\Remarks{\begin{remarks}}
\newcommand\EndRemarks{\end{remarks}}
\newcommand\Convention{\begin{convention}}
\newcommand\EndConvention{\end{convention}}
\newcommand\Notation{\begin{notation}}
\newcommand\EndNotation{\end{notation}}
\newcommand\Lemma{\begin{lemma}}
\newcommand\EndLemma{\end{lemma}}
\newcommand\Proposition{\begin{proposition}}
\newcommand\EndProposition{\end{proposition}}
\newcommand\Corollary{\begin{corollary}}
\newcommand\EndCorollary{\end{corollary}}
\newcommand\Claim{\begin{claim}}
\newcommand\EndClaim{\end{claim}}
\newcommand\Proof{\begin{proof}}
\newcommand\EndProof{\end{proof}}
\newcommand\Equation{\begin{equation}}
\newcommand\EndEquation{\end{equation}}

\newcommand\Bullets{\begin{itemize}}
\newcommand\EndBullets{\end{itemize}}

\newcommand\calT{{\mathcal T}}
\newcommand\caln{{\mathcal N}}

\newcommand\frakc{\mathfrak c}

\renewcommand\epsilon{\varepsilon}

\newcommand\ADS{\rm ADS}
\newcommand\DS{\rm DS}

\newcommand\mni{m_{\rm ni}}
\newcommand\minv{m_{\rm inv}}

\newcommand\rk{\mathop{\rm rk}\nolimits}
\newcommand\discup{\mathbin{\rotatebox[origin=c]{90}{$\vDash$}}}

\newcommand\redeta{\eta}
\newcommand\redLambdaplus{\Lambda^+}
\newcommand\redi{i}

\newcommand\redj{j}

\newcommand\redTheoremC{Theorem C}
\newcommand\reallyredj{j}

\newcommand\redA{{A}}
\newcommand\trulyredA{A}
\newcommand\redPL{PL\ }
\newcommand\reddJ{J}
\newcommand\redalpha{\alpha}
\newcommand\redbeta{\beta}
\newcommand\redWbang{W^!}
\newcommand\redV{V}
\newcommand\wasA{B}

\newcommand\wasB{A}
\newcommand\wasa{b}
\newcommand\wasb{a}
\newcommand\redUpsilon{\Upsilon}
\newcommand\redscrSB{\scrS^B}
\newcommand\redstar{^*}
\newcommand\redlambdatwostar{\lambda_2^*}
\newcommand\redPi{\Pi}
\newcommand\reddeltaB{\delta^B}
\newcommand\Uneg{U\redsubtwo}

\newcommand\reallyredR{R}
\newcommand\redX{X}
\newcommand\redFprime{F'}

\newcommand\reddestX{X}
\newcommand\reddestA{I}
\newcommand\reddestB{II}
\newcommand\reddestK{K}
\newcommand\redOmega{ \Omega}

\newcommand\redRdagger{\Omega_1}
\newcommand\redWdagger{W_1}
\newcommand\redcalwdagger{\calw_1}
\newcommand\redcalA{\cala}
\newcommand\reallyredcala{{\cala}}
\newcommand\redcale{\cale}
\newcommand\redlambdadagger{\lambda_2^*}
\newcommand\redcaladagger{{\cale_1}}
\newcommand\redCdagger{C_1}

\newcommand\redWdaggerstar{W_1^{\diamondsuit}}
\newcommand\redcalwdaggerstar{{\mathcal W}^\diamondsuit_1}
 \newcommand\redcaladaggerstar{{\cale_1^\diamondsuit}}
\newcommand\redWdaggerplus{W_1^+}

\newcommand\reddaggerprime{^\prime_{1}}

\newcommand\redp{p}
\newcommand\redRstar{R^*}
\newcommand\redUsigma{U,\sigma}

\newcommand\hatB{\widehat B}

\newcommand\maybeR{R}

\newcommand\chibar{\overline\chi}

\newcommand\inter{\mathop{\rm int}}

\newcommand\std{{\rm std}}
\newcommand\maybeast{\ast}
\newcommand\maybeprime{'}

\newcommand\reddenedone{1}

\newcommand\tH{\widetilde H}
\newcommand\tM{\widetilde M}
\newcommand\tf{\widetilde f}
\newcommand\teta{\widetilde \eta}

\newcommand\tC{\widetilde C}

\newcommand\tN{\widetilde N}

\newcommand\tY{\widetilde Y}

\newcommand\redPsi{\Psi}
\newcommand\tpsi{\widetilde\psi}
\newcommand\redtpsi{\tpsi}
\newcommand\redcalU{\calU}
\newcommand\redtphi{\tphi}

\newcommand\tQ{\widetilde Q}

\newcommand\tphi{{\widetilde \phi}}

\newcommand\ta{thickened annulus}

\newcommand\hatf{{\widehat f}}
\newcommand\hatV{{\widehat V}}
\newcommand\hatcalB{{\widehat \calB}}
\newcommand\hatpsi{{\widehat \psi}}

\newcommand\talpha{{\widetilde\alpha}}

\newcommand\tbeta{{\widetilde\beta}}

\newcommand\calB{{\mathcal B}}

\newcommand\cala{{\mathcal A}}

\newcommand\calp{{\mathcal P}}

\newcommand\calf{{\mathcal F}}

\newcommand\calL{{\mathcal L}}

\newcommand\tg{\mathop{\rm tg}}

\newcommand\calt{{\mathcal T}}
\newcommand\calb{{\mathcal B}}

\newcommand\calg{{\mathcal G}}
\newcommand\cale{{\mathcal E}}

\newcommand\redgamma{\gamma}

\newcommand\redone{1}
\newcommand\redtwo{2}
\newcommand\redsubone{_{1}}
\newcommand\redsubtwo{_{2}}

\newcommand\calu{{\mathcal U}}
\newcommand\calU{{\mathcal U}}

\newcommand\calw{{\mathcal W}}
\newcommand\calv{{\mathcal V}}

\newcommand\Csmooth{C^{\rm smooth}}
\newcommand\Csmoothprime{C'_{\rm smooth}}

\newcommand\ZZ{{\bf Z}}

\newcommand\QQ{{\mathbb Q }}

\newcommand\cals{{\mathscr S}}

\newcommand\scrA{{\mathscr A}}

\newcommand\scrS{{\mathscr S}}
\newcommand\scrT{{\mathscr T}}
\newcommand\scrV{{\mathscr V}}
\newcommand\scrL{{\mathscr L}}

\newcommand\Ltriv{{\mathscr L}^{\|}}
\newcommand\Lacyl{C}

\newcommand\vol{\mathop{\rm vol}}

\newcommand\image{\mathop{{\rm Im}}}

\newcommand\Fr{\mathop{\rm Fr}}

\newcommand\realize{realize}
\newcommand\realiz{realiz}
\newcommand\abs{a-b splitting}
\newcommand\abt{a-b triad}
\newcommand\anabs{an \abs}
\newcommand\anabt{an \abt}

\newcommand\pagelike{pagelike}
\newcommand\bindinglike{bindinglike}

\newcommand\jay{J}

\newcommand\boldE{{\Phi}}
\newcommand\boldW{{\Sigma}} 
\newcommand\boldV{{\Sigma_0}}
\newcommand\boldU{{\Sigma}}
\newcommand\FF{{\bf F}}

\newcommand\newredcalT{\mathcal T}
\newcommand\reddenedG{G}
\newcommand\redgnought{g_0}
\newcommand\redg{g}
\newcommand\redQ{Q}

\DeclareFontFamily{U}{rcjhbltx}{}
\DeclareFontShape{U}{rcjhbltx}{m}{n}{<->rcjhbltx}{}
\DeclareSymbolFont{hebrewletters}{U}{rcjhbltx}{m}{n}

\let\aleph\relax\let\beth\relax
\let\gimel\relax\let\daleth\relax

\DeclareMathSymbol{\aleph}{\mathord}{hebrewletters}{39}
\DeclareMathSymbol{\beth}{\mathord}{hebrewletters}{98}
\DeclareMathSymbol{\gimel}{\mathord}{hebrewletters}{103}
\DeclareMathSymbol{\daleth}{\mathord}{hebrewletters}{100}

\DeclareMathSymbol{\lamed}{\mathord}{hebrewletters}{108}
\DeclareMathSymbol{\mem}{\mathord}{hebrewletters}{109}
\DeclareMathSymbol{\ayin}{\mathord}{hebrewletters}{96}
\DeclareMathSymbol{\tsadi}{\mathord}{hebrewletters}{118}
\DeclareMathSymbol{\qof}{\mathord}{hebrewletters}{114}
\DeclareMathSymbol{\shin}{\mathord}{hebrewletters}{152}

\begin{document}

\title
[
Actualizing subgroups  in small submanifolds
]
{
Actualizing subgroups of 3-manifold groups in
homologically small submanifolds
}

\author{Rosemary K. Guzman}
\address{Department of Mathematics
\\
University of Illinois\\
1409 W. Green St.\\
Urbana, IL 61801}
\email{rguzma1@illinois.edu}

\author{Peter B. Shalen}
\address{Department of Mathematics, Statistics, and Computer Science
(M/C 249)\\
University of Illinois at Chicago\\
851 S. Morgan St.\\
Chicago, IL 60607-7045}
\email{petershalen@gmail.com}

\maketitle

\begin{abstract}
Let $Y$ be a simple $3$-manifold, and let $A$ be a finitely generated, freely indecomposable
subgroup of $\pi_1(Y)$. Set $\eta=\dim H_1(A;\FF_2)$.  
Suppose that {\it either} (a) $\partial
Y\ne\emptyset$ {\it or} (b) 
$\dim H_1(Y;\FF_2)\ge3\eta^2-4\eta+4$.

Under these hypotheses, we show that
$A$ is carried
by some
compact,
connected three-dimensional
submanifold
$Z$ of $\inter Y$ 
such that 
(1) $\partial Z$ is  non-empty, and each of its components is incompressible in $Y$;
(2) the Euler characteristic of $Z$ is  bounded below by $1-\eta$; and
(3) $\dim H_1(Z;\FF_2)\le
3\eta^2-4\eta+1$.

The conclusion implies that any boundary component of $Z$ is an
incompressible surface of genus at most $\eta$. In Case (b), this
should be compared with earlier results proved by Agol-Culler-Shalen
and Culler-Shalen, which provide a surface of genus at most $\eta$
under weaker hypotheses (the lower bound on $\dim H_1(Y;\FF_2)$ being
linear in $\eta$ rather than quadratic), but do not give any relationship between the given subgroup $A$ and this surface.

In a forthcoming paper we will apply the result to give a new upper bound for the ratio of the rank of the mod 2 homology of a closed, orientable hyperbolic $3$-manifold to the volume of the manifold.

\end{abstract}




\section{Introduction}

Theorem 1.1 of \cite{singular-two} asserts that if $Y$ is a closed,
orientable, hyperbolic 3-manifold, and if $\pi_1(Y)$ has a subgroup isomorphic to a
genus-$g$ (orientable) surface group for a given positive integer $g$,
and if the dimension of $H_1(Y;\FF_2)$ is at least $\max(3g - 1, 6)$, then
$Y$ contains a connected (and in particular non-empty) incompressible closed surface whose genus is
at 
most $g$. (The
piecewise linear category is the default category for manifolds
throughout the present paper. 
To say that a closed PL surface $F\subset Y$ is {\it incompressible}
means that $F$ is not a sphere and that the inclusion homomorphism
$\pi_1(F)\to\pi_1(Y)$ is injective. This term, and others explained in this
introduction, will be defined in a more formal setting in the body of
the paper.)

As is pointed out in \cite{singular-two}
(and in its predecessor \cite{acs-singular}, where a weaker result is
proved), Theorem 1.1 of \cite{singular-two}  may be regarded as a partial analogue of Dehn's lemma
for closed incompressible surfaces.

From
Theorem 1.1 of \cite{singular-two}, it is not hard to deduce a result
about more general subgroups of $\pi_1(Y)$ when $Y$ is a closed,
orientable, hyperbolic $3$-manifold. We define a group to be {\it
  freely indecomposable} if it is not trivial or infinite cyclic, and
is not the free product of two non-trivial subgroups. 
The following consequence of \cite[Theorem
1.1]{singular-two}  will be proved (using slightly different
terminology) 
as Proposition \ref{real Proposition A} of this paper.

\begin{PropositionA}
Let $Y$ be a closed, orientable, hyperbolic $3$-manifold, and
let $A$ be a 
finitely generated,
freely  indecomposable subgroup of $\pi_1(Y)$. Suppose
that 
$\dim H_1(Y;\FF_2)\ge\max(3\redeta-1,6)$, where 
$\redeta=\dim H_1(A;\FF_2)$.
Then $Y$ contains a connected closed incompressible surface of genus
at most $\redeta$.
\end {PropositionA}

The proof of 
Proposition A given in the body of the paper uses the 
following result, which is stated and proved as Proposition \ref{real
  Proposition B}. It is a strict improvement over
\cite[Prop. 1.1]{sw}. 

\begin{PropositionB}
Let $p$ be a prime, and let $M$ be a closed $3$-manifold. Assume that
either $M$ is orientable or $p=2$. Let $E$ be
a finitely generated subgroup of $\pi_1(M)$, and set
$\eta=\dim H_1(E;\FF_p)$. If 
$\dim H_1(M;\FF_p)
\ge\max(3,\eta+2)$, then $E$ has infinite index in $\pi_1(M)$. In fact, $E$ is contained in infinitely many
distinct subgroups of $\pi_1(M)$ having finite index.
\end{PropositionB}

In retrospect, Proposition B could have been stated and proved in
\cite{sw}, as its proof is only a slight modification of the proof of
Proposition 1.1 of that paper. 
Proposition B 
seems likely to find other applications.

Neither \cite[Theorem 1.1]{singular-two} nor Proposition A above
gives any direct connection between the subgroup of $\pi_1(Y)$
appearing in the hypothesis and the incompressible surface appearing
in 
the conclusion. The main result of this paper, which we state below as
\redTheoremC\ 
 and prove in the body of the paper (in a slightly stronger
form) as Theorem \ref{main
  desingular}, remedies this situation, but at the cost of requiring a
quadratic lower bound on $\dim H_1(Y;\FF_2)$ in terms of 
$\dim H_1(A;\FF_2)$, 
rather than the linear lower bound required for
\cite[Theorem 1.1]{singular-two} or Proposition A.

The natural setting for \redTheoremC\  is purely topological and is  broader than that of closed hyperbolic
manifolds. In the body of this paper, we define a $3$-manifold to be
{\it simple} if it is  compact, orientable, irreducible and
boundary-irreducible,  
and its fundamental group
is infinite and has no rank-$2$ free
abelian subgroups. For a closed, connected, orientable $3$-manifold,
the existence of a hyperbolic metric implies simplicity (and the
converse is included in Perelman's geometrization theorem), but there
are many simple $3$-manifolds with non-empty boundary.

A subgroup $A$ of the fundamental group of a path-connected
space $Y$ will be said to be {\it
  carried by}
a path-connected subset  $Z$  of $Y$ if $A$ is conjugate to a subgroup of the image of
the inclusion homomorphism $\pi_1(Z)\to\pi_1(Y)$. 

\begin{TheoremC}
Let $Y$ be a simple 
$3$-manifold, 
and let $A$ be a finitely generated, freely indecomposable
subgroup of $\pi_1(Y)$. Set $\eta=\dim H_1(A;\FF_2)$.  
Suppose that {\it either} (a) $\partial
Y\ne\emptyset$ {\it or} (b) 
$\dim H_1(Y;\FF_2)\ge3\eta^2-4\eta+4$.
Then 
$A$ is carried
by some
compact,
connected,
three-dimensional
submanifoldtion
$Z$ of $\inter Y$ 
such that 
\begin{itemize}
\item
$\partial Z$ is  non-empty, and each of its components is incompressible in $Y$;
\item
The Euler characteristic of $Z$ is  bounded below by $1-\eta$; and
\item
$\dim H_1(Z;\FF_2)\le
3\eta^2-4\eta+1$.
\end{itemize}
\end{TheoremC}

Note that since the submanifold $Z$ given by Theorem C has Euler
characteristic at least $1-\eta$, its boundary has Euler
characteristic at least $2-2\eta$; and that since $Z$ has no sphere boundary
components, an arbitrary boundary component has Euler
characteristic at least $2-2\eta$, and therefore has genus at most
$\eta$. Thus $Z$, which carries $A$, provides the promised
relationship between the subgroup $A$ of the hypothesis and an
incompressible surface of genus at most $\eta$. Note also that the
conclusion provides a quadratic upper bound for $\dim H_1(Z;\FF_2)$ in
terms of $\eta$. Thus $Z$ is the ``homologically small submanifold''
referred to in the title.


If $Y$ and $A$ satisfy the hypotheses of \redTheoremC, and if $Y$ is
closed (for example), then the incompressibility of the boundary components of the
submanifold $Z$ given by the theorem, together with the irreducibility
of $Y$, implies that each component of $\partial Z$ has non-positive
Euler characteristic; since the Euler characteristic of $\partial Z$ is twice that
of $Z$, an arbitrary component $F$ of $\partial Z$ has Euler
characteristic at least $2-2\eta$, and is therefore an incompressible
surface of genus at most $\eta$. In particular the genus of $F$ is
bounded above by  
$\dim H_1(A;\FF_2)$. 
Thus the hypotheses of \redTheoremC,
which are stronger than those of Proposition A, yield a stronger
conclusion. 
In particular, the fact that $A$ is carried by $Z$, of which
$F$ is an arbitrary boundary component, provide the promised connection
between $A$ and $F$.

In \cite{acs-singular}, \cite{singular-two}, \cite{ds}, \cite{CDS},
\cite{fourfree}, \cite{kfreevol}, \cite{ratioII}, and \cite{moreDS},
the topological results of \cite{acs-singular} and  \cite{singular-two} were applied to the problem of relating the
volume of a 
hyperbolic $3$-manifold $M$ to its
topological invariants, specifically the dimension of
$H_1(M;\FF_2)$. (The arguments given in these papers implicitly
establish certain cases of Proposition A.) \redTheoremC\  above is in large part motivated by the same
problem. In fact, in a subsequent paper we expect to apply \redTheoremC\ 
to obtain an upper bound for 
$\vol(M)/\dim(H_1(M;\FF_2))$ for an
arbitrary 
closed, orientable hyperbolic $3$-manifold
$M$ which is very nearly as strong as the bound obtained in \cite{ratioII}
under certain topological restrictions on $M$.

 The proof of \redTheoremC\  will be given in Section \ref{tower
  section}, after we lay the groundwork in the earlier sections. Like
the proofs of the topological results in \cite{acs-singular} and
\cite{singular-two},  the proof of \redTheoremC\  uses a tower of
two-sheeted coverings analogous to the one used by Shapiro and
Whitehead in their proof of Dehn's lemma
\cite{shapiro-whitehead}. 
In sketching the proof here, we shall use some
conventions that are introduced formally in the body of the paper. For a
compact manifold $X$ we write $\chibar(X)=-\chi(X)$, where $\chi(X)$
denotes the Euler characteristic, and $\rk_p(X)=\dim H_1(X;\FF_p)$ for any
prime $p$. 

The sketch that will be given here is organized
somewhat differently from the proof of Lemma \ref{key lemma} that will
be given in the body of the paper. The point of view taken in this
sketch is intended to appeal directly to the intuition in
explaining why the upper bound for $\rk_2(Z)$ is  quadratic in $\eta$,
given a suitable quadratic lower bound on $\rk_2(Y)$ (or the
assumption $\partial Y\ne\emptyset$). The more technical mode of
organization used in the body of the paper is more convenient for the
purpose of giving a formal proof with good explicit quadratic bounds.

Under the hypotheses of \redTheoremC, we show
that there is a piecewise linear map $\phi$ from some compact PL space
$K$ to $Y$ such that the homomorphism
$\phi_\sharp:\pi_1(K)\to\pi_1(Y)$ is injective and has image conjugate
to $A$. It then suffices to produce a submanifold $Z$ of $\inter Y$
having the properties stated in \redTheoremC, such that $\phi$ is
homotopic to a map $g$ whose image is contained in $Z$.

The methods of \cite{acs-singular} provide a ``tower''
$\calt = (M_0,N_0,\Pi_1,M_1,N_1,\Pi_2,\ldots,\Pi_n,M_n,N_n)$, where
the ``height'' $n$
 is a non-negative integer; the $M_j$ and $N_j$ are simple
 $3$-manifolds, and the ``base'' $M_0$ is $Y$;
the manifold $N_j$ is a submanifold of $\inter M_j$ for each $j$, and $\partial
N_j$ is incompressible in $M_j$; and $\Pi_j :M_j
\to N_{j-1}$ is a two-sheeted covering map for
$j=1,\ldots,n$.
For each $j$, let $\iota_j : N_j \to M_j$ denote the
inclusion map. For the purpose of the present sketch, let us  define a
map $\psi_j$ from $N\doteq N_n$ to $N_j$ for $j=0,\ldots,n$ by setting 
$\psi_j=\redPi_{\reallyredj+1}\circ\iota_{\reallyredj+1}\circ
\redPi_{\reallyredj+2}\circ\cdots\circ\iota_{n-1}\circ
\redPi_n\circ\iota_n$ for $j<n$, and defining $\psi_n$ to be the
identity map.
The tower given by the methods of \cite{acs-singular} has the property that there is a PL map
$\tphi:K\to N$ (called a ``homotopy lift'' of $\phi$) such that
$\psi_0\circ\tphi:K\to M_0$ is homotopic to
$\phi$, and $\tphi_*:H_1(K;\FF_2)\to H_1(N;\FF_2)$ is surjective.
(Of course the presence of mod $2$ coefficients here reflects the fact
that the $\Pi_j$ are two-sheeted covering maps.)
The proof of \redTheoremC\  also involves other, more technical properties of the tower
which will be omitted from the present sketch.

The basic strategy underlying the proof of \redTheoremC\  is to construct
recursively, for $i=0,\ldots,n$, a map
$\redgamma_{n-\redi}:
K\to N_{n-\redi}$ which is homotopic to $\psi_{n-i}\circ\tphi$, and
has its image contained in a 
submanifold  $R_{n-\redi}$ of
  $\inter N_{n-\redi}$ 
having certain topological properties.
When we set $i=n$ we obtain a 
map $g\doteq\gamma_0:K\to M$, homotopic to $\phi$, whose image is
contained in $Z\doteq R_0$; the topological information about $R_{n-i}$,
specialized to the case $i=n$, will show that $Z$ has the properties
stated in
\redTheoremC, thus completing the proof of the theorem.

To begin the recursive construction of the $\gamma_{n-i}$ and $R_{n-i}$,
we simply define $R_n$ to be the complement of a small open boundary collar
in $N$, and set $\gamma_n=\tphi$.

In the recursion step, we are given an index $i$ with $0\le i<n$, a
submanifold $R_{n-i}$ of $\inter N_{n-i}$ having certain
topological properties, and a map $\redgamma_{n-\redi}:
K\to N_{n-\redi}$ which is homotopic to $\psi_{n-i}\circ\tphi$, and
has its image contained in $R_{n-\redi}$. In particular $R_{n-i}$ is a
submanifold
of  $\inter M_{n-i}$, and $M_{n-i}$ is in turn  a two-sheeted covering
of
$N_{n-(i+1)}$, with covering projection $\Pi_{n-i}$. We are required to find a
submanifold $R_{n-(i+1)}$ of $\inter N_{n-(i+1)}$ having
topological properties analogous to those of $R_{n-i}$, 
and a map $\redgamma_{n-(\redi+1)}:
K\to N_{n-(\redi+1)}$ which is homotopic to $\psi_{n-(i+1)}\circ\tphi=\Pi_{n-i} \circ(\psi_{n-i}\circ\tphi)$, and
has  image contained in $R_{n-(\redi+1)}$. 

This situation is abstracted in Lemma \ref{key
  lemma}, which
is the central technical result
of the paper and provides a large part of the recursion step in the
proof of \redTheoremC. 
In that lemma, we are given
a simple $3$-manifold  $Q$,     a
connected
two-sheeted covering $\tQ$  of $Q$ (which  itself is automatically simple),
a three-dimensional 
submanifold $R$ of $\inter\tQ$ with certain topological properties,
a compact, connected PL space $K$ with freely indecomposable
fundamental group, and a $\pi_1$-injective \redPL map $f:K\to \tQ$ 
with  $f(K)\subset R$. 
The lemma provides a
three-dimensional
 submanifold $T$ of
$\inter Q$ having topological properties partly analogous to those of $R$, 
and a \redPL map
from
$K$ to $Q$ which has image contained in $T$ and
is homotopic to 
$\redPi\circ f$, where $\redPi:\tQ\to Q$ denotes the covering
projection.

The  construction of $T$ used in the proof of Lemma \ref{key
  lemma}, is
relatively simple. One may assume that
$R$ has been chosen within its isotopy class in such a way that it has transverse intersection with its image under the non-trivial
deck transformation $\tau$ of $\tQ$, and that this intersection has no
components which are trivial or redundant in a suitable sense. The short Section
\ref{good reps} is devoted to establishing the existence of the
required isotopy, which is encoded in Lemma \ref{bestest curve
  reduction}.
Given this choice of $T$ within its isotopy class, $\Pi(R)$ is a
submanifold of $Q$. The submanifold $T$ is constructed by modifying
$\Pi(R)$ in such a way that its boundary components become incompressible; the
modifications involved are addition of $2$-handles, removal of
$1$-handles, capping off boundary $2$-spheres and discarding
components. The
free indecomposability of $\pi_1(K)\cong A$ is used to show that the
modifications can be done in such a way that $T$ is connected and 
$\redPi\circ f$ is homotopic to a map with image contained in $T$.
Some of the relevant topological properties of $T$ follow rather
easily from the construction: 
The incompressibility of the
components of $\partial T$, combined with the connectedness of $T$
and the simplicity of $Q$,
implies that $T$ is simple. 
It is also easy to show that $0\le\chibar(T)\le\eta-1$. 
It is important to observe that a priori $\partial T$  could be empty.

A deeper topological property of $T$
is expressed
in terms of a decomposition, canonical up to isotopy, of a general simple $3$-manifold
$X$.
 This decomposition, the theory of which
is developed in Section \ref{book section}, is called an a-b splitting
for reasons explained in that section. It is
defined by a system $\cale_X$ of disjoint, properly embedded, $\pi_1$-injective annuli  in $X$ which
separate $X$ into two (possibly disconnected) 
submanifolds $C_X$ and $ W_X$; 
here $(C_X,\cale_X)$ is ``acylindrical'' in a suitable sense, and $W_X$ is a ``book of $I$-bundles.''
By definition, a book of $I$-bundles $W$ (referred to more formally in
the body of the paper as the underlying manifold of a book of $I$-bundles) itself contains 
a disjoint system $\cala$ of properly embedded, $\pi_1$-injective annuli  which
separate $W$ into a disjoint union of ``bindings'' and a disjoint
union of ``pages''; each binding is a solid torus, each page has the
structure of an $I$-bundle over a connected surface of strictly
negative Euler characteristic, and $\cala$ is
the union of the vertical boundaries of the pages.  In the a-b
splitting of a simple manifold $X$, each component of $\cale_X$ is
contained in the intersection of $\partial X$ with a binding of $W_X$.

The
a-b splitting of a simple $3$-manifold is closely related to
the characteristic submanifold theory \cite{Jo}, \cite{JS}. Our
account of it involves a number of facts about the characteristic
submanifold which are widely known but are not  systematically
presented in the literature. We have attempted to remedy this by
proving these facts in Section \ref{suction}.

The deeper property of $T$ which we have mentioned is that one of two
alternatives holds: either (i)
 $\chibar(T)=\chibar(R)$,
and the pairs $(C_T,\redcale_T)$ and $(C_R,\redcale_R)$ are
homeomorphic; or
(ii) $\chibar(T)<\chibar(R)$, and 
for every prime $p$, the quantity $\rk_p(T)$ is bounded above
by an explicit linear combination (given in the statement of Lemma
\ref{key lemma}) of the quantities $\rk_p(C_R)$,
$\chibar(R)$ and $\chibar(T)$.

The proof that either Alternative (i) or Alternative (ii) holds
breaks up into two cases.  Our choice of $R$ within its isotopy class
guarantees that every component of $R\cap\tau(R)$ has non-positive
Euler characteristic and is $\pi_1$-injective in $\tQ$. We first
consider the case in which every component of $R\cap\tau(R)$ has 
Euler characteristic $0$. In this case one can show that the
components of $R\cap\tau(R)$ are solid tori, and that  for a suitable
choice of the a-b splitting within its isotopy class,
$R\cap\tau(R)$ is contained in both $ W_R$ and $\tau(W_R)$, and is
disjoint from $\cale_R$ and $\tau(\cale_R)$. 
The proof of this fact uses properties of the a-b splitting of $R$, which in turn follow from well-known properties
of the characteristic submanifold. It also uses, in an essential way, that $\partial
R\ne\emptyset$, which is a crucial hypothesis of Lemma \ref{key
  lemma}. (We shall see that this is an issue in the recursion argument.)

In the situation where the
components of $R\cap\tau(R)$ are solid tori, contained in both $ W_R$ and $\tau(W_R)$ and 
disjoint from $\cale_R$ and $\tau(\cale_R)$, one can show that  $\Pi$
maps $C_R$ homeomorphically onto a submanifold of $Q$, while
$\Pi(W_R)$  inherits the structure of a ``generalized
book of $I$-bundles'' in which the bindings are Seifert fibered
spaces. The modifications of $\Pi(R)$ which are involved in the construction
of $T$
may involve enlarging $\Pi(R)$ by adding a disjoint union of solid
tori; this union of solid tori, whose boundaries are contained in the
boundary of $\Pi(W_R)$, will be denoted here by $\scrT$. The union of $\Pi(W_R)$ with
$\scrT$ is a genuine book of $I$-bundles. If  $\Pi(R)\cup\scrT$ has
incompressible boundary in $Q$, it is equal to $T$; and in this
subcase we show that $C_T=\Pi(C_R)$ and that $W_T=W_R\cup\scrT$, which
implies that Alternative (i) holds. If $\Pi(R)\cup\scrT$ has a
compressible boundary in $Q$, then elementary homological arguments
show that Alternative (ii) holds.

We now turn to the case in which some component of $R\cap\tau(R)$ has
strictly negative
Euler characteristic. In this case, one shows that Alternative
(ii) holds. The first part of Alternative
(ii)---that 
$\chibar(T)<\chibar(R)$---is rather straightforward. The proof of the
second part involves an analysis of the manifold $\Pi(R)$.
Up to homeomorphism, one can think of $\Pi(R)$
as being obtained from $R$ by certain self-gluings involving the
involution $\tau$. The self-gluings can be divided into two stages.
First, one performs the gluings involving components of
$R\cap\tau(R)$ having
Euler characteristic $0$, giving a $3$-manifold which is denoted
$\Omega_1$ in the proof of Lemma \ref{key lemma}.
Then one obtains $\Pi(R)$ from $\Omega_1$ by self-gluings involving
homeomorphic copies of one or more components of
$R\cap\tau(R)$ having
strictly negative
Euler characteristic. The topological structure of $\Omega_1$ in this
case is very similar to the topological structure of $\Pi(R)$ in the
previous 
case, and this can be used to bound $\rk_p(\Omega_1)$  above
by an explicit linear combination of  $\rk_p(C_R)$,
$\chibar(R)$ and $\chibar(T)$. This is combined with  elementary
homological arguments to
show that the second part of Alternative (ii) holds in this case.
This completes our sketch of the proof of Lemma \ref{key
  lemma}.

We have mentioned that  Lemma \ref{key
  lemma} requires the hypothesis $\partial R\ne\emptyset$. On the
other hand, the conclusion of  Lemma \ref{key
  lemma} does not include the assertion that $\partial
T\ne\emptyset$. This means that Lemma \ref{key
  lemma}  does not by itself cover the recursion step in the proof of
\redTheoremC. We shall therefore first indicate how one obtains the
conclusion of \redTheoremC\  under the assumption that $\partial
R_{n-i}\ne\emptyset$ for each $i$, and then indicate a bootstrapping
argument to show that this
assumption always holds.

Under the assumption that $\partial
R_{n-i}\ne\emptyset$ for each $i$, it follows from Lemma \ref{key
  lemma} that for every $i$ with $0\le i<n$, one of the alternatives
(i) and (ii) stated above holds with $R_{n-i}$ and $R_{n-(i+1)}$ in
place of $R$ and $T$. For the purpose of this sketch, let us say that
an index $i$ is {\it special} if Alternative (ii) holds. Let us write
$c_i=\chibar(R_{n-i})$ and $r_i=\rk_2(R_{n-i})$ for each $i$. We have
$c_{i+1}\le c_i$ for every $i<n$, and this inequality is strict if $i$
is special; furthermore, it is not hard to show that $c_i\ge0$ for
every $i\le n$, and that $c_0<\eta$. It follows that 
$c_i<\eta$ for each
index $i$, and that
the number of
special indices is less than $\eta$.

If $i$ is a special index, then by specializing Alternative (ii) to
the case $p=2$, we deduce $r_{i+1}$ is bounded
above by a linear combination of the quantities $\rk_2(C_R)$,
$c_i$ and $c_{i+1}$, and therefore by a linear combination of
$\rk_2(C_R)$ and $\eta$. Using homological results about the a-b splitting that are established in Section \ref{book section}, one can also show that
$\rk_2(C_R)$ exceeds $r_i$ by at most an explicit multiple of $\eta$,
and therefore that: 
\newline
 {\bf (I)} {\it Whenever $i$ is special, $r_{i+1}-r_i$ is bounded 
above by  an explicit multiple of $\eta$.}

On
the other hand, if $i_1$ and $i_2$ are indices with $0\le i_1<i_2\le
n$, and if there is no special index $i$ with $i_1\le i<i_2$, then
Alternative (i) must hold for each $i$ with $i_1\le i<i_2$; this
implies that the pairs $(C_{R_{i_1}},\cale_{R_{i_1}})$ and
$(C_{R_{i_2}},\cale_{R_{i_2}})$ are homeomorphic.
Combining this observation with homological results about the a-b splitting that are established in Section \ref{book section},
along with the inequalities $c_{i_1}<\eta$ and $c_{i_2}<\eta$, one
can now deduce that:
\newline
{\bf (II)} {\it  If $0\le i_1<i_2\le
n$, and if there is no special index $i$ with $i_1\le i<i_2$,  then 
$r_{i+1}-r_i$ is bounded 
above by  an explicit multiple of $\eta$.} 

Combining the facts (I) and
(II) with the fact that the number of
special indices is less than $\eta$, 
we deduce that $\rk_2(Z)=r_n$ is bounded by a sum with fewer than
$\eta$ terms, in which each term is bounded by an explicit multiple of
$\eta$. Thus we obtain the promised quadratic upper bound for $\rk_2(Z)$.

The proof that the condition  $\partial R_{n-i}\ne\emptyset$ always
holds may be thought of as an argument by  strong
induction on $i$. Suppose that $\partial R_{n-k}$ is known to be non-empty whenever $0\le k<i$; we wish to prove
$\partial R_{n-i}\ne\emptyset$. If 
we assume $\partial R_{n-i}=\emptyset$, then $N_{n-i}$ is
closed and $R_{n-i}= N_{n-i}$.
When $N_{n-i}$ is
closed, the methods of \cite{sw}, together with the quadratic lower
bound for  $\rk_2(Y)$ in terms of $\eta$ that is given by the hypotheses,
and the
technical properties of the tower that were omitted above, give a lower
bound on $\rk_2(N_{n-i})$. 
On the other hand, by  applying the strong induction hypothesis, and
using the same
argument that gave an upper  bound on $\rk_2(Z)$ in the argument above
under the extra assumption that the $R_{n-i}$ all have non-empty
boundary, one can obtain an upper  bound on $\rk_2(R_{n-i})=
\rk_2(N_{n-i})$. The  upper bound so obtained turns out to be less
than the lower bound obtained by the  methods of \cite{sw},  so that
we obtain a contradiction. This completes our sketch of the proof of
\redTheoremC.

In the course of this sketch we have described the contents of
Sections \ref{suction} and \ref{good
  reps}. 
We have also mentioned that Section \ref{book section} includes the
needed material about a-b splittings. The earlier parts of Section
\ref{book section} cover useful properties of books of $I$-bundles and
generalized books of $I$-bundles, some of which refine results proved
in \cite{acs-singular}.

Section \ref{pushing section} is devoted to the proof of
Lemma \ref{key lemma}, which has been featured prominently in the
sketch above. 
The section includes some preliminary material:
Subsections \ref{before gist}---\ref{composition phenomenon} are 
needed for the gluing step which is described in the summary of the proof of Lemma
\ref{key lemma} given above, and Lemma \ref{gist} provides the generalized book of
$I$-bundles which is mentioned in that summary.
Section \ref{tower section} contains some background material about
towers, and the proof of \redTheoremC\  (stated there as Theorem
\ref{main desingular}).

Section \ref{prelim section} contains needed background in
$3$-manifold theory. Several of the results are variants of well-known
results in the field, or are facts that are widely believed to be true
but for which references are hard to find. Section \ref{prelim section}  also contains
the proof of Proposition A (stated there as Proposition \ref{real Proposition A}).


\section{Preliminaries}\label{prelim section}

\Notation\label{oiler} If $X$ is a  space having the homotopy type
of a finite CW complex, we shall write $\chi(X)$ for the Euler
characteristic of $X$, and set $\chibar(X)=-\chi(X)$.
If $\Gamma$ is a group such that $K\doteq K(\Gamma,1)$ has the homotopy type
of a finite CW complex, we shall set $\chi(\Gamma)=\chi(K)$  and 
$\chibar(\Gamma)=-\chi(\Gamma)=\chibar(K)$.

For any space $X$ we denote the number of components of $X$ by
$\frakc(X)$. If $A$ is a subset of a space $X$, we shall denote by
$\Fr_XA$ (or by $\Fr A$ when the space $X$ is understood) the frontier
of $A$ in $X$, defined to be $\overline{A}\cap\overline{X-A}$.

The boundary and interior of a manifold $X$ will be denoted
respectively by $\partial X$ and $\inter X$.

\EndNotation

\NotationRemarks\label{redirects here}
Let $p$ be a prime number. If $X$ is a space having the homotopy type
of a finite CW complex, we shall write $\rk_p(X)$ for the dimension of
the $\FF_p$-vector space $H_1(X;\FF_p)$. 
Similarly, if $\Gamma$ is a
finitely generated group, we shall write $\rk_p(\Gamma)$ for the dimension of
the $\FF_p$-vector space $H_1(\Gamma;\FF_p)$.

If $X$ and $Y$ are spaces,
and $f:X\to Y$ is a map, we shall write $\rk_p(f;X, Y)$ for the
rank of the linear map $f_*:H_1(X;\FF_p)\to H_1(Y;\FF_p)$. We
shall write $\rk_p(f)$ in place of $\rk_p(f;X,Y)$ in cases where there is
no danger of ambiguity. If $X$ is a subspace of a space $Y$, we shall
write $\rk_p(X,Y)$ for $\rk_p(i;X,Y)$, where $i:X\to Y$ is the inclusion
map. Thus $\rk_p(X,Y)$ is the rank of the inclusion
homomorphism $H_1(X;\FF_p)\to H_1(Y;\FF_p)$.

Note that if 
$f:X\to Y$ is any map between spaces, and if there exist a space
$Z$ and maps $g:X\to Z$ and $h:Z\to Y$ such that $h\circ g=f$, then
$\rk_p(f)\le\rk_p(Z)$ for any prime $p$.

\EndNotationRemarks

\Definition\label{pie one invective}
A map $f:X\to Y$ between pathwise connected spaces will be termed {\it
  $\pi_1$-injective} if for some (and hence for every) basepoint $x\in
X$, the homomorphism $f_\sharp:\pi_1(X,x)\to\pi_1(Y,f(x))$ is
injective. For arbitrary spaces  $X$ and $Y$,
a  map $f:X\to Y$  will be termed {\it
  $\pi_1$-injective} if for every component $K$ of $X$, the map $f|K$
is a $\pi_1$-injective map from $K$ to the component of $Y$ containing
$f(K)$. A subspace $X$ of a space $Y$ is termed {\it $\pi_1$-injective} if
the inclusion map $X\to Y$ is $\pi_1$-injective.
\EndDefinition

\Convention
In this paper, all statements and arguments involving manifolds are to be interpreted
in the piecewise linear category, except for one brief passage in the
proof of Lemma \ref{erev chag}. Manifolds and submanifolds are
understood to be PL, as are maps, homotopies and isotopies between
manifolds.

\EndConvention

\DefinitionRemarks\label{keep it simple}
A $3$-manifold $X$ is said to be {\it irreducible} if $X$ is connected
and every $2$-sphere in $X$ bounds a ball.
The $3$-manifold $X$ is said to be {\it
  boundary-irreducible} if $\partial X$ is $\pi_1$-injective in
$X$. (According to Dehn's Lemma and the Loop Theorem, $X$ is
boundary-irreducible if and only if for every properly embedded disk
$D\subset X$, the simple closed curve $\partial D$ bounds a disk in
$\partial X$.)

A $3$-manifold $Q$ will be termed {\it simple} if it is
compact, orientable, irreducible and
boundary-irreducible,  
and 
$\pi_1(Q)$ 
is infinite and has no rank-$2$ free
abelian subgroups. (This definition is almost equivalent to the
one given in \cite[Definitions 1.10]{acs-singular}; the only
difference is that a $3$-ball is simple according to the definition
used in \cite{acs-singular}, but is not simple according to the
definition used in this paper.)

It is a standard consequence of the Sphere Theorem that an irreducible
$3$-manifold with infinite fundamental group is aspherical. In
particular, every simple $3$-manifold is aspherical.

If $Q$ is simple, it is irreducible and is not a $3$-ball; hence
$\partial Q$ has no $2$-sphere component.
Furthermore, since $Q$ is boundary-irreducible and $\pi_1(M)$ has no rank-$2$ free
abelian subgroups,
$\partial Q$ has no torus component. It follows that
$\chibar(Q)=\chibar(\partial Q)/2\ge0$ for any simple $3$-manifold
$Q$, and that $\chibar(Q)>0$ for the case of a simple $3$-manifold
$Q$ with $\partial Q\ne\emptyset$.

If $Q$ is a simple $3$-manifold then every connected finite-sheeted
covering space $\tQ$ of $Q$ is simple. Indeed, according to
\cite[Theorem 3]{msy}, the irreducibility of $Q$ implies that
$\tQ$ is irreducible. It is immediate that the other defining
properties of a simple manifold are inherited by $\tQ$ from $Q$.
\EndDefinitionRemarks

\DefinitionRemarks\label{don't squeeze-a da fruit}
Let $F$ be a closed (not necessarily connected) surface
in the interior of an irreducible $3$-manifold $Q$. 
As in \cite{acs-singular}, we 
will say that $F$ is {\it
  incompressible} in $Q$ if it is $\pi_1$-injective and has no
$2$-sphere component. Note that with this definition, the empty
surface in $Q$ is incompressible.
As in \cite{acs-singular}, we will avoid the term ``incompressible''
for surfaces that are not closed.

If $F$ is  any (possibly disconnected) closed surface
in the interior of an irreducible $3$-manifold $Q$, we define a {\it
  compressing disk} for $F$ to be a disk $D\subset\inter Q$ such that
$D\cap F=\partial D$, and $\partial D$ does not bound a disk in
$F$. It is a standard consequence of Dehn's Lemma and the Loop Theorem
that $F$ is incompressible if and only if $F$ has no compressing disk
and no $2$-sphere component.

Note that if $X$ is a compact submanifold of the interior of an irreducible
$3$-manifold $ Q$, and $\partial X$ is incompressible
in $Q$, then $X$ is $\pi_1$-injective in $Q$.

\EndDefinitionRemarks

\Remark\label{more fruit}
If $X$ is a compact, connected submanifold of the interior of an irreducible
$3$-manifold $ M$, and $\partial X$ is incompressible
in $M$, then $X$ is simple. (Indeed, the boundary
irreducibility of $X$ follows from the non-existence of compressing
disks for $\partial X$. The non-existence of rank-$2$ free abelian
subgroups in $\pi_1(M)$ follows from the $\pi_1$-injectivity of $X$ in
$M$ (see \ref{don't squeeze-a da fruit}). To show that $X$ is
irreducible, we note that any  $2$-sphere in $\inter X$ bounds a ball
$B\subset M$. If $B$ were not contained in $X$, it would contain a
component of $\partial X$, contradicting incompressibility; hence
$B\subset X$, and the irreducibility of $X$ is established. Finally,
to show that $\pi_1(X)$ is infinite, we note that if $\partial
X=\emptyset$ then $X=M$ and the assertion is immediate; and if
$\partial X\ne\emptyset$, then since any component of $\partial X$ has
strictly positive genus, $M$ has positive first betti number and
therefore has an infinite fundamental group.)
\EndRemark

The following result will be used in Section \ref{book section}.

\Lemma\label{because g and c-r}
Let $Q$ be a simple $3$-manifold, and let $L$ be a finitely generated
subgroup of $\pi_1(Q)$. Suppose that $L$ has an infinite cyclic normal
subgroup. Then $L$ is itself infinite cyclic.
\EndLemma

\Proof
By applying the theorem of \cite{scott} to the covering space of $Q$
defined by $L$, we see that
 $L$
 is isomorphic to the fundamental group of a compact, connected, orientable
$3$-manifold $N$. According to \cite[Theorem 3.15]{hempel}, $N$ may be
written as a connected sum $M_1\#\cdots\# M_n$, where each $M_i$ is
prime (in the sense defined on p. 27 of \cite{hempel}). If at least
two of the $M_i$ were non-simply connected, then $L\cong\pi_1(N)$
would be a non-trivial free product, a contradiction to the hypothesis
that $L$ has an infinite cyclic normal
subgroup. Hence after possibly re-indexing the $M_i$ we may assume
that $\pi_1(M_i)$ is trivial for every $i>1$ and hence that the
fundamental group of $M\doteq M_1$ is isomorphic to
$L$. Now, since $M$ is prime, it follows from \cite[Lemma
3.13]{hempel} that either
$M$ admits an $S^2$-fibration over $S^1$---in which case the conclusion
of the present lemma holds---or $M$ is irreducible. 

In the
case where $M$ is irreducible, we apply \cite[Theorem 1.1]{cj} or
\cite[Corollary 8.8]{gabai},
according to which a compact, irreducible $3$-manifold whose
fundamental group contains an infinite cyclic normal subgroup
admits a
Seifert fibration. In the present situation, the simplicity of $Q$
implies that $\pi_1(M)\cong L\le\pi_1(Q)$ has no rank-two free abelian
subgroup. Note also that $\pi_1(M)\cong L$ is in particular infinite. But any
irreducible, orientable, compact Seifert fibered space $M$ whose
fundamental group is infinite and has no rank-two free abelian
subgroup is homeomorphic to $D^2\times S^1$. Thus the conclusion holds
in this case as well. 
\EndProof

\Number\label{tgdef}
If $F$ is a closed surface whose components are orientable, the {\it
  total genus} of $F$ is defined to be the sum of the genera of the
components of $F$, and will be denoted $\tg(F)$.

\EndNumber

\Lemma\label{insurance}
Let $Q$ be a compact, connected, orientable $3$-manifold. Then for any
compact $3$-manifold $T\subset Q$ and any prime $p$, we have
$\rk_p(T)\le\rk_p(T,Q)+\tg(\partial T)$.
\EndLemma

\Proof
All homology groups in this argument will be understood to have
coefficients in $\FF_p$. We may assume without loss of generality that $T\subset\inter Q$. We
set $X=\overline{Q-T}$. Thus $X\cup T=Q$, while $F\doteq X\cap T=\partial
T \subset\partial X$. We denote by $K$ and $L$ the respective
kernels of the inclusion homomorphisms $H_1(T)\to H_1(Q)$ and
$H_1(F)\to H_1(X)$. From the exactness of the fragment
$$H_1(F)\to H_1(X)\oplus H_1(T)\to H_1(Q)$$
of the Mayer-Vietoris sequence, we deduce that $K$ is contained in the
image of $L$ under the inclusion
homomorphism $H_1(F)\to H_1(T)$. Hence $\dim K\le\dim L$. 

It is a standard observation that if $X$ is any compact, orientable
$3$-manifold, the kernel of the inclusion homomorphism $H_1(\partial
X)\to H_1(X)$ is self-orthogonal in the antisymmetric intersection pairing on
$H_1(\partial X)$. 
In particular $L$ is self-orthogonal. But since the components of $F$
are closed surfaces, the intersection pairing on $H_1(F)$ is
non-degenerate. Hence $\dim L\le(\dim H_1(F))/2=\tg(F)$.

Since $K$ is the kernel of the inclusion homomorphism 
$H_1(T)\to H_1(Q)$, and $\rk_p(T,Q)$ is the dimension of its image, we
have
$$\rk_p(T)=\dim H_1(T)=\rk_p(T,Q)+\dim K\le \rk_p(T,Q)+\dim L\le \rk_p(T,Q)+\tg(F).$$
\EndProof

\Definition\label{stand with us}
A group will be termed {\it freely indecomposable} if it is neither
trivial nor infinite cyclic, and cannot be written as a free product
of two non-trivial subgroups.
\EndDefinition

\Lemma\label{from DS} (Cf. \cite[Lemma 8.5]{moreDS}).
Let $M$ be a closed $3$-manifold, and let $p$ be a prime. Suppose that
either $M$ is orientable or $p=2$.
Suppose
that 
$h:\pi_1(M)\to\ZZ/p\ZZ\times \ZZ/p\ZZ$ is a surjective homomorphism.
Then 
$$\rk_p(\ker h)\ge2\rk_p(M)-3.$$
\EndLemma

\Proof

Let $\Gamma_1$ denote the normal subgroup of $\pi_1(M)$ generated by
all 
commutators 
and $p$-th powers. 
According to  \cite[Lemma 1.5]{sw}, if $n$ is any integer less than or
equal to $(\dim V)-2$,  if $E$ is any subgroup generated by $n$
elements of $\pi_1(M)$, and if $D$ denotes the subgroup $ E\Gamma_1$
of $\pi_1(M)$, we have $\dim H_1(D;\FF_p) \ge 2n+1$. To apply this, we
let $P$ denote the kernel of the homomorphism from $H_1(M;\FF_p)$
to $\FF_p^2$ induced by $h$; we set $n=\dim P=\rk_p(M)-2$;
we choose elements $x_1,\ldots,x_n$ of $\pi_1(M)$ whose images in $H_1(M;\FF_p)$
form a basis of $P$; and we take $E$ to be the subgroup of
$\pi_1(M)$ generated by $x_1,\ldots,x_n$. Then $D\doteq E\Gamma_1=\ker
h$, and hence
$$\rk_p(\ker h)=\dim H_1(D;\FF_p) \ge 2n+1=2\rk_p(M)-3.$$


\EndProof


The following result was stated in the introduction (in slightly
different language) as Proposition B.
As mentioned in the introduction, it is an improvement over
\cite[Prop. 1.1]{sw}, and is proved by basically the same method.

\Proposition\label{real Proposition B}
Let $p$ be a prime, and let $M$ be a closed $3$-manifold. Assume that
either $M$ is orientable or $p=2$. Let $E$ be
a finitely generated subgroup of $\pi_1(M)$, and set
$\eta=\rk_p(E)$. If $\rk_p (M) \ge\max(3,\eta+2)$, then $E$ has infinite index in $\pi_1(M)$. In fact, $E$ is contained in infinitely many
distinct subgroups of $\pi_1(M)$ having finite index.
\EndProposition

\Proof
Set $\eta'=\max(\eta,1)$. By induction on $n\ge0$ we shall show that
there is a subgroup $L_n$ having index $4^n$ in $\pi_1(M)$ such that
$E\le L_n$ and $\rk_p(L_n)\ge\eta'+2$. This immediately implies the
conclusion of the proposition.

For the base case $n=0$, we set $L_0=\pi_1(M)$. The inequality $\rk_p
(L_0) \ge\eta'+2$ is a restatement of the hypothesis  $\rk_p (M)
\ge\max(3,\eta+2)$, and the other asserted properties of $L_0$ are
trivial. Now suppose that we are given and $n\ge0$ and a subgroup
$L_n$ having the asserted properties. Since $E\le L_n$,
$\rk_p(E)=\eta\le\eta'$, and $\rk_p(L_n)\ge\eta'+2$, there is a
surjective homomorphism $h:L_n\to\ZZ/2\ZZ\times \ZZ/2\ZZ$ such that
$L_n\doteq\ker h$
contains $E$. 

Since $L_n$ has finite index in $\pi_1(Y)$, the group $L_n$ is
isomorphic to the
fundamental group of a closed $3$-manifold $M$, and if $p$ is odd then
$M$ is orientable. It therefore follows from
Lemma \ref{from DS}
that 
$$\rk_p(L_{n+1})=\rk_p(\ker
h)\ge2\rk_p(L_n)-3\ge2(\eta'+2)-3=2\eta'+1.$$ 
Since $\eta'\ge1$, it 
follows that $\rk_p(L_{n+1})\ge\eta'+2$. Furthermore, by construction
we have $|\pi_1(M):L_{n+1}|=4|\pi_1(M):L_n|=4^{n+1}$, and $E\le
L_{n+1}$. Thus the induction is complete.
\EndProof


The following result was stated in the introduction (in slightly
different language) as Proposition A:

\Proposition \label{real Proposition A}
Let $Y$ be a closed, simple $3$-manifold, and
let $A$ be a
finitely generated, 
freely indecomposable subgroup of $\pi_1(Y)$. 
Set $\redeta =\rk_2(A)$, and assume that
$\rk_2(Y)\ge\max(3\redeta -1,6)$.
Then $Y$ contains a connected closed incompressible surface of genus
at most $\redeta $.
\EndProposition


\Proof
The hypothesis implies that $\rk_2(Y)\ge \max(1,\redeta +2)$, which with
Proposition \ref{real Proposition B} implies that $A$ has infinite index in
$\pi_1(Y)$.
Hence the covering space $\tY$ of $Y$ determined by the subgroup $A$
is non-compact.

Since $A$ is finitely generated, the main theorem of
\cite{scott-core}
gives a compact, connected, three-dimensional submanifold $K$ of $\tY$
such that the inclusion homomorphism $\pi_1(K)\to\pi_1(\tY)$ is an
isomorphism. In particular, $\pi_1(K)$ is isomorphic to $A$.

Suppose that every component of $\partial K$ is a sphere (which is in
particular the case if $\partial K=\emptyset$). Since the simple
manifold $Y$ is aspherical by
\ref{keep it simple}, $\tY$ is also aspherical, and hence each
component $F$ of $\partial Y$ bounds a homotopy $3$-ball $B_F$ in
$\tY$. If $B_F\supset K$ for some component $F$ of $\partial K$, then
$\pi_1(K)\cong A$ is trivial, a contradiction to the free
indecomposability of $A$. If $B_F\cap K=F$ for every component $F$ of
$\partial K$, then $\tY$ is compact, and again we have a
contradiction. Hence some component $F_0$ of $\partial K$ has strictly
positive genus.

Suppose that the inclusion homomorphism $\pi_1(F_0)\to\pi_1(K)$ has
non-trivial kernel. Then by the loop theorem, there is a properly
embedded disk $D$ in $K$ such that $C\doteq\partial D$ is a
homotopically non-trivial simple closed curve in $F_0$. If $D$
separates $K$, the non-triviality of $C$ implies that both components
of $K-D$ are non-simply connected, and hence $A\cong\pi_1(K)$ is a
free product of two non-trivial groups, a contradiction to free
indecomposability. If $D$ does not
separate $K$, then $A\cong\pi_1(K)$ is a
free product of an infinite cyclic group with some (possibly
trivial) group, and again we have a contradiction to free
indecomposability. Hence the inclusion homomorphism
$\pi_1(F_0)\to\pi_1(K)$ is injective. Since the inclusion homomorphism
$\pi_1(K)\to\pi_1(\tY)$ is an isomorphism, and the covering projection
$\tY\to Y$ is $\pi_1$-injective,  $\pi_1(Y)$ has a subgroup isomorphic to a
genus-$g$ (orientable) surface group, where $g>0$ denotes the genus of
$F_0$. 

It is a standard consequence of Poincar\'e-Lefschetz duality that the
total genus of the boundary of the compact, orientable $3$-manifold
$K$ is at most $\dim H_1(K;\QQ)=\dim H_1(A;\QQ)\le \redeta $. Hence $0<g\le
\redeta $. With the hypothesis of the present proposition, we now deduce that
$\rk_2(Y)\ge\max(3g - 1, 6)$. It now follows from 
Theorem 1.1 of \cite{singular-two} (the statement of which was reviewed
in the introduction) that
$Y$ contains a connected incompressible closed surface whose (positive) genus is
at
most $g\le \redeta $.
\EndProof

\Lemma\label{surgery lemma}
Let $\redQ $ be a simple $3$-manifold, let $K$ be a compact, connected PL space 
with freely indecomposable fundamental group, and let
$\redgnought:K\to 
\inter
\redQ $ be a $\pi_1$-injective \redPL map. Let $T_0$ be a 
compact, connected 
three-dimensional
submanifold of $\inter\redQ $ such that $\redgnought(K)\subset T_0$. Suppose that no component of $\partial T_0$ is a sphere. Then
there is a 
compact, connected three-dimensional 
submanifold  $T$  of 
$\inter\redQ $ 
such that 
\begin{enumerate}[(a)]
\item$\partial T$ is incompressible in $\redQ $;
\item$\chibar(T)\le\chibar(T_0)$;
\item  $\rk_p(T)\le \rk_p(T_0)$ for every prime $p$;
\item the image of the inclusion homomorphism   $H_1(T;\ZZ)\to
  H_1(\redQ ;\ZZ)$ is contained in the image of the inclusion homomorphism   $H_1(T_0;\ZZ)\to
  H_1(\redQ ;\ZZ)$, so that in particular $\rk_p(T,Q)\le\rk_p(T_0,Q)$
  for every prime $p$; and 
\item$\redgnought$ is homotopic to a map 
$\redg  $ 
such that $\redg (K)\subset T$. 
\end{enumerate}
Furthermore, if we assume in addition
$\partial T_0$ is compressible in $\redQ $, and that no component of $\partial T_0$ is a torus, then $T$ may be chosen so that
\begin{enumerate}
\item[(b')]$\chibar(T)<\chibar(T_0)$.
\end{enumerate}
\EndLemma

\Proof
Let us say that a closed, orientable (possibly disconnected)
$2$-manifold $G\subset \inter Q$ is obtained by an {\it essential surgery}
from another closed, orientable (possibly disconnected)
$2$-manifold $F\subset \inter Q$ if there exist an annulus $A\subset F$ and
disks $D_+, D_-\subset \inter Q$ such that $\partial(D_+\cup D_-)=\partial
A=(D_+\cup D_-)\cap F$, $G=(F-A)\cup D_+\cup D_-$, and neither of the
curves $\partial D_+,\partial D_-$ bounds a disk in $F$. The annulus
$A$ will be called the {\it domain annulus} of the essential surgery. 

The following
observation will be useful:

\Claim\label{new rokhel}
Let $F\subset \inter Q$ be a closed orientable surface, and suppose that
$G\subset \inter Q$ is obtained from $F$ by an essential surgery. Then
$\chibar(G)=\chibar(F)-2$.

Furthermore, if we assume that $F$ has no
$2$-sphere component, then either (i) $G$ has no $2$-sphere component,
or (ii) the domain annulus of the essential surgery is a
non-separating annulus in a
torus component $Y$ of $F$, and $G$ is the disjoint union of $F-Y$ with a
$2$-sphere.
\EndClaim

Let $V$ be a compact three-dimensional submanifold of $\inter Q$. We define an {\it improving
  ball} for $V$ to be a
$3$-ball $E\subset \inter Q$ 
such that 
$A\doteq E\cap\partial V$ is an annulus in $\partial
E$, and no component of $\overline{(\partial V)-A}$ is a
disk. We observe that in this situation we have
either $E\subset V$ or $E\cap V=A$; in these respective cases
we say that $E$ is an {\it internal improving
  ball} or an {\it external improving
  ball.} 
A compact
submanifold $W$ of $\inter Q$ is said to be {\it obtained from $V$ by an internal
  improvement} if there is an internal improving ball $E$ for $V$ such
that $W=\overline{V-E}$; and $W$  is said to be {\it obtained from $V$ by an external
  improvement} if there is an external improving ball $E$ for $V$ such
that $W=E\cup V$.

If $V$ is a compact three-dimensional submanifold of $\inter Q$, and $D$ is a
compressing disk for $\partial V$, then by enlarging $D$ slightly we
obtain an improving ball for $V$. Hence:

\Claim\label{well, compress mah boundary!}
If $V$ is a compact three-dimensional submanifold of $\inter Q$, and
$\partial V$ admits a
compressing disk, then $V$ admits an improving ball.
\EndClaim

We also observe:
\Claim\label{how high}
If $V$ is a compact, connected three-dimensional submanifold of $\inter Q$,
and if
$W$ is a submanifold of $\inter Q$ obtained from $V$ by an internal or external improvement, then:
\begin{itemize}
\item $\partial W$ is obtained from $\partial V$ by an essential surgery;
\item $W$ has either one or two components;
\item if $W$ has two components, then it is obtained from $V$ by an
  internal improvement, and $\partial W$ is obtained from $\partial V$
  by an essential surgery whose domain annulus separates the component of
  $\partial V$ which contains it.
\end{itemize}
\EndClaim

Now we claim:
\Claim\label{beito shel gvir}
If $V$ is a compact, connected three-dimensional submanifold of $\inter Q$,
and if $W$ is a submanifold obtained from $V$ by an internal or
external improvement, then $\chibar(W)=\chibar(V)-1$. If in addition
we assume that no component of $\partial V$ is a $2$-sphere, then for each
component $X$ of $W$ we have $\chibar(X)<\chibar(V)$.
\EndClaim

To prove \ref{beito shel gvir}, first note that, according to \ref{how
  high}, $\partial W$ is obtained from $\partial V$ by an essential
surgery. We let $A$ denote the domain annulus of the surgery.  By \ref{new rokhel}, we have $\chibar(\partial
W)=\chibar(\partial V)-2$. This gives $\chibar(
W)=\chibar(\partial
W)/2=\chibar(\partial V)/2-1=\chibar( V)-1$, and the first assertion
of \ref{beito shel gvir}
is established. To prove the second assertion, we add the assumption
that $\partial V$ has no $2$-sphere component. According to \ref{how high}, $W$
has either one or two components. If $W$ has one component, the
conclusion of the second assertion of \ref{beito shel gvir},
is immediate from the first
assertion. Now suppose that $W$ has two  components $X_1$ and
$X_2$. According to the final assertion of \ref{new rokhel}, either (i) $\partial W$ has no $2$-sphere component,
or (ii) $A$ is a
non-separating annulus in a
torus component of $F$. But according to the final assertion of
\ref{how high}, since $W$ has two components, $A$ must separate the component of
  $\partial V$ which contains it; thus (ii) cannot occur. Hence
  $\partial W$ has no $2$-sphere component, and therefore
$\chibar(X_i)=\chibar(\partial X_i)/2\ge0$ for $i=1,2$. Since we have
$\chibar(X_1)+\chibar(X_2)=\chibar(W)=\chibar(V)-1$ by the first
assertion of 
\ref{beito shel gvir},
it now follows
that $\chibar(X_i)<\chibar(V)$ for $i=1,2$, as required for the second
assertion. This
completes the proof 
of \ref{beito shel gvir}.

Now let   $J$ denote the image of the inclusion homomorphism   $H_1(T_0;\ZZ)\to
  H_1(\redQ ;\ZZ)$. 
Fix a triangulation of $K$, and for every positive integer $d$ let $K^{(d)}\subset K$ denote the
underlying space of its $d$-skeleton.
Let $\scrV$ denote the set of all 
compact, connected submanifolds  $V$  of $\inter Q $ such that 
\begin{enumerate}
\item$\chibar(V)\le\chibar(T_0)$;
\item  $\rk_p(V)\le \rk_p(T_0)$ for every prime $p$;
\item the image of the inclusion homomorphism   $H_1(V;\ZZ)\to
  H_1(\redQ ;\ZZ)$ is contained in $J$; and 
\item
the map $\redg_0|K^{(2)}:K^{(2)}\to \inter Q$ is homotopic to a \redPL map
$\redg_V :K^{(2)}\to \inter Q $ such that
$\redg_V(K^{(2)})
\subset V$.
\end{enumerate}
Let $\scrV_1$ denote the subset of $\scrV$ consisting of all elements
$V$ such that $\partial V$ has no $2$-sphere components. Note that we
have $T_0\in\scrV_1$, and hence $\scrV_1\ne\emptyset$.

We claim:
\Claim\label{new outside}
If $V$ is an element of $\scrV$, and $W\subset \inter Q$ is a manifold
obtained from $V$ by an external improvement, then $W\in\scrV$.
\EndClaim

To prove \ref{new outside}, we first note that since $W$ is obtained by an
external improvement from the connected manifold $V$, it is itself connected.
We must show that
Conditions (1)--(4) of the definition of $\scrV$ hold with  $W\doteq E\cup
V$ in place of $V$. 
Each of these conditions on $W$ will be
established using the corresponding condition on $V$. Using the first
assertion of \ref{beito shel gvir}, we find that
$\chibar(W)<\chibar(V)\le\chibar (T_0)$, which gives Condition
(1). 
The definition of an external improvement implies that the inclusion homomorphism 
$\iota:H_1(V;\ZZ)\to
  H_1(W ;\ZZ)$ is surjective, which gives $\rk_p(W)\le\rk_p(V)\le\rk_p(T_0)$ for
  any prime $p$; this establishes Condition (2). The surjectivity of $\iota$
  also shows that the image of the inclusion homomorphism
  $H_1(W;\ZZ)\to H_1(Q;\ZZ)$ coincides with the image of the inclusion homomorphism
  $H_1(V;\ZZ)\to H_1(Q;\ZZ)$, which is in turn contained in $J$, and (3) is established. Finally,
  since $W\supset V$, Condition (4) for $V$ implies Condition (4) for
  $W$ upon setting 
$g_{W}=g_V $. 
This completes the proof of
  \ref{new outside}.

  Next we claim: 
\Claim\label{new inside}
If $V$ is an element of $\scrV_1$, and $W\subset \inter Q$ is a manifold
obtained from $V$ by an internal improvement, then some component of
$W$ belongs to $\scrV$.
\EndClaim

To prove
\ref{new inside},
we write $W=\overline{V-E}$, where $E$ is an internal improving ball
for $V$. 
Note that $W$ is $\pi_1$-injective in $V$.
We set $k=\frakc(W)$; as observed in \ref{how high}, $k$ is
either $1$ or $2$. The distinction between the cases $k=1$ and $k=2$
will appear at several points in the proof of
\ref{new inside}. In the case $k=2$, we will denote by $X_1$ and $X_2$
the components of $W$, and we will set $X_i^*=X_i\cup E$ for $i=1,2$.

We fix a base
point
$e\in V$;
in the 
case $k=2$
we take
$e\in E$, and 
in the case $k=1$ we take $e\in W$.
In the case $k=2$,
for $i=1,2$, we denote by $L_i$ the image of the
inclusion homomorphism $\pi_1(X_i^\star,e)\to\pi_1(V,e)$. Then
$\pi_1(V,e)$ is the free product of its subgroups $L_1$ and $L_2$.
In the case $k=1$, 
we denote by $L$ the image of the
inclusion homomorphism $\pi_1(W,e)\to\pi_1(V,e)$. Then
$\pi_1(V,e)$ is the free product of $L$ with an infinite cyclic
subgroup.

Since $V\in\scrV$, 
Condition (4) of the definition of $\scrV$ 
asserts
that
the map 
$\redgnought|K^{(2)}$ 
is homotopic to a \redPL map
$\redg_V :K^{(2)}\to \inter Q $ such that $\redg_V
(K^{(2)})\subset V$.
We may write $g_V=j_V\circ h_V$, where $h_V$ is a \redPL map from
$K^{(2)}$ to $V$ and $j_V:V\to \inter Q$ denotes the inclusion. Since
$g_V$ is homotopic to $g_0|K^{(2)}$, and $g_0$ is
$\pi_1$-injective by hypothesis, the map $g_V$ is
$\pi_1$-injective, and in particular $h_V$ is $\pi_1$-injective.

We  may suppose $h_V$ to be chosen
within its homotopy class in
such a way that $h_V(q)=e$ for some base point $q\in K$.
Let $D$ denote
the image of the homomorphism $(h_V)_\sharp:\pi_1(K,q)\to\pi_1(V,
e)$. Since $h_V$
is $\pi_1$-injective, 
$D$ is isomorphic to
$\pi_1(K)$ and is therefore freely indecomposable. It follows from the
Kurosh subgroup theorem that a freely indecomposable subgroup of a
free product is  conjugate to a subgroup (necessarily not infinite
cyclic) of a factor. Hence $D$ is conjugate in $\pi_1(V,e)$ to a
subgroup of $L_1$ or $L_2$ in 
the case $k=2$, 
and to a subgroup  of $L$ in
the case $k=1$. This implies that in 
the case $k=2$,  the map 
$h_V^{(1)}\doteq h_V|K^{(1)}:K^{(1)}\to V$ is freely homotopic to a \redPL map 
whose image is contained in $X_1^\star$ or $X_2^\star$; and that in 
the case $k=1$, the map $h_V^{(1)}:K^{(1)}\to V$
is freely homotopic to a \redPL map 
whose image is contained in $W$. 
But in the case $k=2$, the set $X_i$ is a
deformation retract of $X_i^*$ for $i=1,2$. Hence in all cases, there
exist a component $\redWbang $ of $W$ such that
$h_V^{(1)}:K^{(1)}\to V$ is freely homotopic to the composition of a \redPL
map $h_{\redWbang }^{(1)}:K^{(1)}\to {\redWbang }$ with the inclusion $\kappa:{\redWbang }\to V$.

Since $\kappa\circ h_{\redWbang }^{(1)}:K^{(1)}\to V$ is homotopic to 
$h_V^{(1)}\doteq h_V|K^{(1)}:K^{(1)}\to V$,
the homotopy extension property of PL spaces implies that $\kappa\circ
h_{\redWbang }^{(1)}$ can be extended to a \redPL map $\lambda: K^{(2)}\to V$. Since
$\lambda(K^{(1)})=h_{\redWbang }^{(1)}(K^{(1)})\subset {\redWbang }$, and the
component $\redWbang $ of $W$ is $\pi_1$-injective in $V$, there is a \redPL map
$h_{\redWbang }:K^{(2)}\to \redWbang $ such that $h_{\redWbang }|K^{(1)}=\lambda|K^{(1)}=\kappa\circ
h_{\redWbang }^{(1)}$. Hence if $j_{\redWbang }:\redWbang \to \inter Q$ denotes the inclusion map,
$j_{\redWbang }\circ h_{\redWbang }|K^{(1)}:K^{(1)}\to \inter Q$ is homotopic to
$j_{V^!}\circ h_{V^!}|K^{(1)}:K^{(1)}\to \inter Q$ and hence to
$g_0|K^{(1)}$. Now set $g_{W^!}=j_{\redWbang }\circ h_{\redWbang }:K^{(2)}\to \inter Q$. Since the  simple manifold $Q$  is aspherical by \ref{keep
  it simple}, and the maps $g_{\redWbang },g_0|K^{(2)}
:K^{(2)}\to \inter Q$ restrict to homotopic maps from $K^{(1)}$ to $\inter Q$, the
maps $g_{\redWbang }$ and $g_0|K^{(2)}
$ must themselves be homotopic. Since
$g_{\redWbang }(K^{(2)})\subset \redWbang $, this shows that Condition (4) of the definition of
$\scrV$ holds with $\redWbang $ in place of $V$.

To complete the proof of 
\ref{new inside},
it 
now suffices 
to show
that Conditions (1)--(3) of the definition of
$\scrV$ also hold with $\redWbang $ in place of $V$.

Since $V\in\scrV_1$, the boundary of $V$ has no $2$-sphere component. 
Since $W$ is a submanifold obtained from $V$ by an internal
improvement, it now follows from the final assertion of \ref{beito shel gvir} that the component
$X$ of $W$ satisfies $\chibar(X)<\chibar(V)$. But since in particular
we have $V\in\scrV$, Condition (1) of the
definition of $\scrV$ gives 
$\chibar(V)\le\chibar(T_0)$. Hence $\chibar(\redWbang )<\chibar(T_0)$, which establishes
 Condition (1) of the definition of
$\scrV$ for $\redWbang $.

To establish Condition (2) of the definition of $\scrV$ for $\redWbang $, we
note that in the case where $k=2$  we have $\rk_p(V)=\rk_p(X_1)+\rk_p(X_2)$, and
that in the case where $k=1$ we have $\rk_p(V)=\rk_p(W)+1$; since $\redWbang $ is one of
the $X_i$ in the case $k=2$ and is equal to $W$ in the case $k=1$, we always have
$\rk_p(\redWbang )\le\rk_p(V)$. Since $V$ satisfies Condition (2), it follows
that $\rk_p(\redWbang )\le\rk_p(T_0)$, which is Condition (2) for $\redWbang $. To
establish Condition (3) for $\redWbang $, we
note that since $\redWbang \subset V$,
the image of the inclusion homomorphism   $H_1(X;\ZZ)\to
  H_1(\redQ ;\ZZ)$ is contained in the image of the inclusion homomorphism   $H_1(V;\ZZ)\to
  H_1(\redQ ;\ZZ)$, which is in turn contained in $J$ since  $V$
  satisfies Condition (3). This completes the proof of \ref{new inside}.

We have observed that $\scrV_1\ne\emptyset$. Furthermore, for every
$V\in\scrV_1$, since $\partial V$ has no $2$-sphere component, we have
$\chibar(V)=\chibar(\partial V)/2\ge0$. We may therefore define a
non-negative integer $m$ by setting $m=\min_{V\in\scrV_1}\chibar(V)$. 
Since we have observed that $T_0\in\scrV_1$, we have
\Equation\label{ivy twine}
\chibar(T_0)\ge m.
\EndEquation

We
claim:
\Claim\label{omnibus}
If $V$ is any element of $\scrV_1$ with $\chibar(V)=m$, then either
\begin{itemize}
\item $V$ is incompressible, or
\item there exist a torus component $Y$ of
$\partial V$ and an element $U$ of $\scrV_1$ such that $\partial
U=(\partial V)-Y$.
\end{itemize}
\EndClaim

To prove \ref{omnibus}, we will assume that $\partial V$ is not
incompressible, and show that the second alternative of the conclusion
holds. First 
note that since $\partial V$  has no
$2$-sphere components by the definition of  $\scrV_1$, 
the assumption that $\partial V$ is not
incompressible implies that $\partial V$  has a compressing disk
  (see \ref{don't squeeze-a da fruit}).
Hence by \ref{well, compress mah boundary!}, $V$ admits an improving
ball; there therefore exists a submanifold $W$ of $\inter Q$ which is
obtained from $V$ by an internal or external improvement.
It follows from \ref{new inside} in the case where $W$  is
obtained from $V$ by an internal  improvement, and from \ref{new outside} in the case where $W$ is
obtained from $V$ by an external  improvement, that some component $X$ of
$W$ belongs to $\scrV$. Since no boundary component of $V\in\scrV_1$
is a $2$-sphere, it follows from \ref{beito shel gvir} that
$\chibar(X)<\chibar(V)=m$. According to the definition of $m$, it
follows that $X\notin\scrV_1$. As $X\in\scrV$, the definition of
$\scrV_1$ now implies that $\partial X$ has a $2$-sphere component. In
particular, $\partial W$ has a $2$-sphere component.

According to \ref{how high}, $\partial W$ is obtained from $\partial
V$ by an essential surgery. Since $\partial W$ has a $2$-sphere
component and $\partial X$ does not, it follows from 
the final assertion of \ref{new rokhel} that the domain annulus $A$ of the essential surgery is a
non-separating annulus in a
torus component $Y$ of $\partial V$, and $\partial W$ is the disjoint
union of $(\partial V)-Y$ with a
$2$-sphere $S$. 
Since in particular $A$ does not separate the component
of $\partial V$ containing it, it follows from the final assertion of 
\ref{how high} that $W$ is connected, i.e. $X=W$.

Since the simple manifold
$Q$ is in particular irreducible, the $2$-sphere $S$ bounds a ball
$B\subset \inter Q$. Hence either $W\subset B$ or $W\cap B=S$.

Consider the case in which $W\subset B$. Let
$\beta:\pi_1(V)\to\pi_1(\inter Q)$ denote the inclusion homomorphism. 
In the subcase where $W$ is obtained from $V$ by an external
improvement, we have $V\subset W\subset B$, so that $\beta$ is the
trivial homomorphism. In the subcase where $W$ is obtained from $V$ by an internal
improvement, there is an internal improving ball $E$ for $V$ such that
$W=\overline{V-E}$. Then $V$ is the union $W$ with the ball $E$, and $E\cap W$ is the
disjoint union of two disks; since $W\subset B$, it follows that the
image of $\beta$ has rank at most $1$ in this subcase. Thus in both
cases,
the image of $\beta$ is cyclic. Now
since $Q$ is aspherical by \ref{keep it simple},  $\pi_1(Q)$ is
torsion-free; hence
the image of $\beta$ is trivial or infinite
cyclic.
Now if $\redg_V :K^{(2)}\to \inter Q $ is the map given by Condition (4) of the definition of
$\scrV$, then since
$g_V$ is homotopic to $g_0|K^{(2)}$, and $g_0$ is
$\pi_1$-injective by hypothesis, the map $g_V$ is
$\pi_1$-injective. 
Since 
$\redg_V (K^{(2)})\subset V$, and $\beta$ has trivial or infinite
cyclic image,
$\pi_1(K)\cong\pi_1(K^{(2)})$ 
is either
trivial or infinite cyclic; this contradicts the free indecomposability
of $\pi_1(K)$. 
Hence 
the case  $W\subset B$ cannot occur, and 
we must have $W\cap B=S=\partial B$.

It now follows that  
$U\doteq W\cup B$ 
is a submanifold of $\inter Q$ with
$\partial U=(\partial \redV)-Y$. In particular we have $\chibar(\partial
U)=\chibar(\partial \redV)$, and hence
$\chibar(U)=\chibar(\partial U)/2=\chibar(\partial \redV)/2=\chibar(\redV)=m\le\chibar(T_0)$,
where the last inequality follows from (\ref{ivy twine}).
Thus Condition (1) of the definition of $\scrV$ holds with $U$ in
place of $V$.
On the other hand, since $W=X\in\scrV$, Conditions (2)--(4) of the definition of $\scrV$ hold with $W$ in
place of $V$. Now since $U\doteq W\cup B$, where $W\cap B$ is the boundary
of the ball $B$, we have $\rk_p(U)=\rk_p(W)$, and the image of the inclusion homomorphism   $H_1(U;\ZZ)\to
  H_1(\redQ ;\ZZ)$ coincides with the image of the inclusion homomorphism   $H_1(W;\ZZ)\to
  H_1(\redQ ;\ZZ)$; since $W$ satisfies Conditions (2) and (3), it
  follows that $U$ also satisfies these conditions. Likewise, since
  $U\subset W$, and $W$ satisfies Condition (4), it
  follows that $U$ also satisfies this condition. This shows that
  $U\in\scrV$.
Since $\partial U=(\partial W)-Y$, the $2$-manifold $\partial W$ has
no $2$-sphere component, and hence $U\in\scrV_1$. Thus $U$ has the
properties asserted in the second alternative of the conclusion of
\ref{omnibus}, and the proof of \ref{omnibus} is complete.


  Let us now set $\scrV_2=\{V\in\scrV_1:\chibar(V)=m\}\ne\emptyset$, and
$n=\min_{V\in\scrV_2}\frakc(\partial V)$. We fix an element $T$ of
 $\scrV_2$ with $\frakc(\partial T)=n$. 
Since $T\in\scrV_2\subset\scrV_1\subset\scrV$, it follows from the
definition of $\scrV$ that $T$ is a compact, connected
three-dimensional submanifold of $\inter Q$.
Since $T\in\scrV_2$, it
 follows from \ref{omnibus} that either $\partial T$ is
 incompressible, or there exist a torus component $Y$ of $\partial T$
 and an element $Z$ of $\scrV_2$ such that $\partial Z=(\partial
 T)-Y$. 
But the latter alternative gives $\frakc(\partial Z)=\frakc(\partial T)-1=n-1$,
 contradicting the minimality of $n$. Hence $\partial T$ is incompressible,
i.e.  Condition (a)
of the  conclusion of the lemma
 holds with this choice
of $T$. Since  $T\in\scrV$, 
Conditions (b)--(d) also hold. To show
that Condition (e) holds, we note that by Condition (4) in the
definition of $\scrV$, the map $\redg_0|K^{(2)}:K^{(2)}\to \inter Q$ is
homotopic to a \redPL map
$\redg_T :K^{(2)}\to \inter Q $ such that
$\redg_T(K^{(2)})
\subset T$. But since $\partial T$ is incompressible and $Q$ is
simple, it follows from Remark \ref{more fruit} that $T$ is simple,
and is therefore aspherical by \ref{keep it simple}. Hence $g_T$
extends to a \redPL map $g$, with domain $K$, such that $g(K)\subset T$. Since
$g_0$ and $g$, regarded as maps from $K$ to $\inter Q$, restrict to homotopic
maps from $K^{(2)}$ to $\inter Q$, the asphericity of $Q$ implies that $g_0$
and $g$ are homotopic. Thus Condition (e) holds, and 
the
first conclusion 
of the lemma
 is established.

To prove the second conclusion of the lemma, it suffices to show that if $\partial
T_0$ is not incompressible in $\redQ $, and if no component of
$\partial T_0$ is a torus, then with the choice of $T$ described above
we have $\chibar(T)<\chibar(T_0)$; equivalently, we need to show that
$\chibar(T_0)>m$. Since we have observed that $T_0\in\scrV_1$, the
definition of $m$ guarantees that $\chibar(T_0)\ge m$. Assume that
$\chibar(T_0)=m$. Then according to \ref{omnibus}, either
$\partial
T_0$ is  incompressible in $\redQ$, or some component of
$\partial T_0$ is a torus. This contradiction completes the proof of
the second conclusion.
\EndProof

\DefinitionsRemarks\label{spanish onion}
Let $F_0$ and $F_1$ be compact two-dimensional submanifolds of a
$3$-manifold $M$. Suppose that each $F_i$ is either properly embedded
in $M$ or contained in $\partial M$. A $3$-dimensional submanifold $P$
of $M$ will be called a {\it parallelism} between $F_0$ and $F_1$ if
there is a homeomorphism $h:F_0\times[0,1]\to P$ such that
$h_i(F_0\times\{i\})=F_i$ for $i=0,1$, and $(\partial F)\times[0,1]\subset\partial M$. We shall say that $F_0$ and
$F_1$ are {\it parallel} if there is a parallelism between $F_0$ and
$F_1$. Note that parallel surfaces are in particular disjoint and homeomorphic.

A submanifold of $M$ will be called a {\it thickened surface} if it is
a parallelism between two properly embedded surfaces in $M$ or,
equivalently, if it is a regular neighborhood of a properly embedded
surface in $M$. We define a {\it core} of a thickened surface
$N\subset M$ to be a properly embedded surface in $M$ having $N$ as a
regular neighborhood. The uniqueness of a core of a given thickened
surface $N$, up to isotopy in $N$, is guaranteed by regular
neighborhood theory.

We will frequently need to consider a thickened surface having an
annulus as a core. Such a thickened surface will be called a {\it \ta}.

A compact, properly embedded $2$-manifold in $M$ will be termed {\it
  boundary-parallel} if it is parallel to some $2$-manifold in
$\partial M$.
\EndDefinitionsRemarks

If a $3$-manifold $N$ is exhibited as an $I$-bundle over a
$2$-manifold $F$, we define a {\it vertical} subset of $N$ to be one which
is a union of fibers. If $q:N\to F$ denotes the $I$-fibration of $N$,
then $V\doteq q^{-1}(\partial N)$ is a two-dimensional vertical submanifold of $N$ contained
in $\partial N$; it will be referred to as the {\it  vertical
  boundary} of $N$. We define the {\it horizontal boundary} of $N$ to
be the $2$-manifold $\overline{(\partial N)-V}$.

We will make frequent use of
 Proposition \ref{holiday-proposition}
below, and its Corollary \ref{holiday-corollary}. The proof of
 Proposition \ref{holiday-proposition} depends on the following lemma.

\Lemma\label{erev chag}
Let $\tau$ be a fixed-point-free, orientation-reversing 
involution of a
compact, orientable $2$-manifold $S$. 
Let $C\subset\inter S$ be a simple
closed curve, and suppose that $C$  is homotopic in $S$ to a
simple closed curve whose image under  $\tau$ is
disjoint from $C$. Then $C$ is isotopic in $S$, by an isotopy constant
on $\partial S$, to a simple closed curve $C'$ such
that $\tau(C')\cap C'=\emptyset$.  
\EndLemma

\Proof
 Like all results about manifolds in this paper, the lemma is to be
 interpreted in the PL category. However, it is convenient to deduce it
 from the following nearly analogous statement in the smooth category: 
 \Claim\label{analogous}
 If $\tau$ is a fixed-point-free, orientation-reversing 
 smooth involution of a
 compact, orientable, smooth $2$-manifold $S$,  if
  $C\subset\inter S$ is a smooth simple
 closed curve, and if $C$  is homotopic in $S$ to a smooth
 simple closed curve whose image under  $\tau$ is
 disjoint from $C$, then $C$ is homotopic in $S$
 to a smooth simple closed curve $C'\subset\inter S$ such
 that $\tau(C')\cap C'=\emptyset$.  
 \EndClaim

 To prove \ref{analogous}, we first observe that the assertion is
 trivial if $C$ is homotopically trivial, and we therefore assume that
 it is not. We fix 
 a $\tau$-invariant Riemannian metric on $S$.

Since the
involution $\tau$ is fixed-point-free, we may arrange that $S$ has a strictly
convex boundary, by choosing the metric on $S$ to be the pull-back of
some metric on $S/\langle\tau\rangle$ in which
$\partial(S/\langle\tau\rangle)$ is strictly convex. We now apply Theorem 2.1
of \cite{FHS}. The statement of  that theorem in Section 2 of \cite{FHS} includes the hypothesis that the surface $S$ is closed,
but it is pointed out in the first paragraph of Section 4 of
\cite{FHS} that the result remains true if $S$ is a compact surface
which has
strictly convex boundary in a given metric. The theorem implies that there is a smooth simple
closed curve $C^*\subset\inter F$ which is homotopic to $C$ and has least length
among all closed curves homotopic to $C$. (This is one of two
alternative conclusions 
of \cite[Theorem 2.1]{FHS}; the other alternative would contradict the
orientability of $S$.) 

By hypothesis there is a simple closed curve $C_1$, homotopic to $C$,
such that  $\tau(C_1)$ is
disjoint from $C$. Now $C^*$ is a shortest curve in the homotopy class
of $C$, and $\tau(C^*)$ a shortest curve in the homotopy class
of $\tau(C_1)$. Since $C$ and $\tau(C_1)$ are disjoint, 
it follows from \cite[Corollary 3.4]{FHS}
that $C^*$ and $\tau(C^*)$ are either disjoint or equal. (Again the
discussion in the first paragraph of Section 4 of
\cite{FHS} justifies the application of \cite[Corollary 3.4]{FHS} to a
surface which need not be closed, but is compact and has strictly
convex boundary.) 
If $C^*$ and $\tau(C^*)$ are  disjoint, the conclusion of \ref{analogous}
follows upon setting $C'=C^*$. If
$C^*$ and $\tau(C^*)$ are  equal, then $\tau$ leaves a small metric
tubular neighborhood $A$ of $C^*$ invariant. Since $S$ is orientable,
$A$ is an annulus; and since $\tau$ is fixed-point-free and reverses
orientation, it must interchange the boundary curves of $A$.  
The
conclusion of
\ref{analogous} therefore holds if we take $C'$ to be one of
the boundary curves of $A$. (Each of these boundary curves is homotopic in $S$ to $C^*$ and
therefore to $C$.) Thus \ref{analogous} is proved.

 Now suppose that we are given a fixed-point-free,
 orientation-reversing, PL
 involution of a
 compact, orientable PL $2$-manifold $S$, and a PL simple
 closed curve $C\subset\inter S$ such that the hypothesis of the lemma
 holds. We
 fix a triangulation of the  surface $S/\langle\tau\rangle$ that realizes the PL
 structure inherited from $S$. We fix a smooth structure on $S/\langle\tau\rangle$
 which is  compatible with its triangulation, in the sense that every
 closed $1$-simplex is a smooth closed arc. Then  $S$ inherits a smooth
 structure which is $\tau$-invariant and is compatible with a
 triangulation realizing its given PL structure. Now by hypothesis, there is a PL simple closed curve $C_1$, homotopic to $C$,
 such that  $\tau(C_1)$ is
 disjoint from $C$. The PL simple closed curves $C$ and $C_1$
 may be arbitrarily well approximated by smooth simple closed curves
 $\Csmooth$ and $\Csmooth_1$. If these approximations are good enough, the hypotheses
 of \ref{analogous} will hold with $\Csmooth$ and $\Csmooth_1$ in place
 of $C$ and $C_1$. Hence $\Csmooth$ is homotopic in $\inter S$
 to a smooth simple closed curve $\Csmoothprime$ such
 that $\tau(\Csmoothprime)\cap \Csmoothprime=\emptyset$.  We may
 approximate $\Csmoothprime$ arbitrarily well by a PL simple closed
 curve $C'$.
 By making our approximations good enough we can guarantee that $C'$ is
 homotopic to $C$ and that $\tau(C')\cap C'=\emptyset$.  Finally,
 since $C,C'\subset\inter
 S$ are homotopic in $S$, it follows from \cite[Theorem 2.1]{Epstein}
 that they are isotopic in $S$ by an isotopy constant
 on $\partial S$.
\EndProof



\Proposition\label{holiday-proposition}
Let $N$ be an orientable $3$-manifold which is an $I$-bundle over a
compact, connected $2$-manifold. Let $V$ denote the vertical boundary
of $N$. Let $\scrA$ 
be a
compact, $\pi_1$-injective, properly embedded $2$-manifold  in
$N-V$, each component of which is an annulus. Then either $\scrA$ has a component which is
boundary-parallel in $N-V$, or $\scrA$ is isotopic, by an isotopy of $N$
which is constant on $V$, to a vertical submanifold of $N$.
\EndProposition

\Proof
We first consider the case in which $\scrA$ is connected. Thus $\scrA$
is an annulus.

The existence of the $\pi_1$-injective annulus $\scrA$ in $N$ implies
that $\pi_1(N)$ is infinite, so that the $2$-manifold which is the
base of the $I$-bundle $N$
is aspherical. Hence:
\Claim\label{it's still Thursday}
$N$ is aspherical.
\EndClaim

Let $H$ denote the horizontal boundary of $N$. We have
$\partial\scrA\subset\inter H$.

We define a covering space $\tN$ of $N$ as follows. If $N$ is a
trivial $I$-bundle we take $\tN=N$ to be the trivial covering. If $N$
is a twisted $I$-bundle we take $\tN$ to be the unique two-sheeted
covering such that the $I$-fibration inherited by $\tN$ is a trivial
$I$-fibration. We let $p:\tN\to N$ denote the covering map. In both
cases we identify $\tN$ with $F\times[0,1]$, where $F$ is a compact,
connected, orientable surface, in such a way that $\tH\doteq
p^{-1}(H)=F\times\{0,1\}$. We set $F_j=F\times\{j\}$ for $j=0,1$.

We fix an embedding $\eta:S^1\times[0,1]\to N$ with
$\eta(S^1\times[0,1])=\scrA$. For $j=0,1$ we set $C_j
=\eta(S^1\times\{j\})$. We choose a
lift $\teta:S^1\times[0,1]\to\tN$, and set $\tC_j
=\teta(S^1\times\{j\})$. By symmetry we may assume
that 
$\tC_0\subset\inter F_0$.

Consider the subcase in which 
$\tC_1\subset\inter F_0$. 
In this subcase, the inclusion $\scrA\hookrightarrow N$ is homotopic
rel $\partial\scrA$ to a map whose image is contained in $\inter H$.
It then follows from \cite[Lemma 5.3]{waldhausen} that $\scrA$ is
parallel in $N$ to an annulus in $\inter F_0$,
and the conclusion holds in this subcase.

We now turn to the subcase in which 
$\tC_1\subset\inter
F_1$. 
For the proof in this subcase, we first observe that if
$\theta$ is any embedding of $S^1$ in $F$, then the map
$(x,t)\mapsto(\theta(x),t)$ is an
embedding of $S^1\times[0,1]$ in $F\times[0,1]=\tN$. 
We shall denote
this embedding by $\teta^\dagger_\theta$, and we shall denote the
immersion $p\circ\teta^\dagger_\theta$ by 
$\eta^\dagger_\theta$.
 We
 claim:

\Claim\label{in this subcase}
There exists an embedding $\theta$  of $S^1$ in $F$ such that
$\eta^\dagger_\theta: S^1\times[0,1]\to N$ is an embedding, and such
that the simple closed curve $p(\theta(S^1)\times\{0\})=
\eta^\dagger_\theta(S^1\times\{0\})$ is isotopic to $C_0$
by an isotopy of $ H$ which is constant on $\partial H$.
\EndClaim

To prove \ref{in this subcase}, first note that if $N$ is a trivial
$I$-bundle then
$\eta^\dagger_\theta: S^1\times[0,1]\to N$ is an embedding
for {\it every}  embedding $\theta$  of $S^1$ in $F$; thus in this
case we need only choose an embedding $\theta$ such that
$\theta(S^1)\times\{0\}=C_0$. Now suppose that $N$ is a twisted $I$-bundle.
Then the non-trivial deck transformation $\sigma$ of $\tN=F\times[0,1]$
is given by $\sigma(x,t)=(\tau(x),1-t)$, where $\tau$ is a
fixed-point-free, orientation-reversing involution of $F$. Since we
are in the subcase where $\tC_1\subset F_1$, there are simple closed
curves $E_0,E_1\subset\inter F$ such that $\tC_j=E_j\times\{j\}$ for
$j=0,1$. The $\pi_1$-injectivity of $\scrA$ implies that the $E_j$ are
homotopically non-trivial, and since $\scrA$ is an annulus, $E_0$ and
$E_1$ are homotopic. On the other hand, $\teta$ maps $S^1\times\{0\}$
onto $E_0\times\{1\}$, while $\sigma\circ\teta$ maps $S^1\times\{1\}$
onto $\tau(E_1)\times\{0\}$; since $\eta$  is an embedding, the images
of $\teta$ and $\sigma\circ\teta$ must be mutually disjoint, and hence
$E_0\cap\tau(E_1)=\emptyset$. Thus $E_0$ is homotopic in $F$ to a
simple closed curve whose image under the involution $\tau$ is
disjoint from $E_0$. It then follows from 
Lemma \ref{erev chag}, applied with $F$ and $E_0$ playing the
respective roles of $S$ and $C$,
that $E_0$ is isotopic in $F$,
by an isotopy constant on $\partial F$,
 to a simple closed curve $E_0'$ such
that $\tau(E_0')\cap E_0'=\emptyset$. The conclusion of \ref{in this
  subcase} then follows upon choosing $\theta$ to be an embedding
whose image is $E_0'$.

It follows from \ref{in this subcase} that we may assume, after
possibly modifying $\scrA$ 
by an isotopy that is constant on
$V$, 
that there is an
embedding 
$\theta$  of $S^1$ in $F$ such that
$\eta^\dagger\doteq\eta^\dagger_\theta: S^1\times[0,1]\to N$ is an
embedding, and agrees on $S^1\times\{0\}$ with the embedding $\eta$. 
We set $\teta^\dagger=\teta^\dagger_\theta$.

Note that the annulus
$\scrA^\dagger \doteq\eta^\dagger(S^1\times[0,1])$ is vertical and
shares the boundary component $C_0$ with $\scrA$. We shall complete the
proof in this subcase by showing that $\scrA$ and $\scrA^\dagger $ are
isotopic by an ambient isotopy of $N$ which is constant on $V$ and on
$C_0$.

We let $o$ denote the standard base point $1\in S^1$. We define a path
$\alpha^\std$ in  $S^1\times[0,1] $ by $\alpha^\std(t)=(o,t)$. (The
subscript `$\std$' stands for ``standard.'') We denote the paths
$\eta\circ\alpha^\std$ and $\eta^\dagger\circ\alpha^\std$
in $N$ by
$\alpha$ and $\alpha^\dagger$ respectively. It follows from the definition
of $\eta^\dagger=\eta^\dagger_\theta$ that $\alpha^\dagger$ is a
vertical path. 
We define lifts $\talpha$ and $\talpha^\dagger$ of $\alpha$ and
$\alpha^\dagger$ respectively by setting $\talpha(t)=\teta(o,t)$ and
$\talpha^\dagger(t)=\teta^\dagger (o,t)$ for each $t\in[0,1]$.
Since we are in the subcase where $\tC_1\subset F_1$, we have
$\talpha(j)\in F_j$ for $j=0,1$. The definition of
$\teta^\dagger=\teta^\dagger_\theta$ implies that $\talpha^\dagger(j)\in F_j$ for $j=0,1$.
Hence if we denote by $P:\tN=F\times[0,1]\to F$ the
projection to the first factor, and define a path $\tbeta$ in
$F_1\subset\tN$ by $\tbeta(t)=(P(\alpha(1-t)),1)$, then
$\talpha(1)=\tbeta(0)$, and hence $\talpha\star\tbeta$ is a
well-defined path in $\tN$. Note that $\talpha\star\tbeta$ is
fixed-endpoint homotopic to the vertical path $\talpha^\dagger$. Hence
if we define a path $\beta$ in $H$ by $\beta=p\circ\tbeta$, then 
$\delta\doteq\alpha\star\beta$ is a
well-defined path in $N$, and
\Claim\label{upstairs downstairs}
$\delta$ is
fixed-endpoint homotopic in $N$ to  $\alpha^\dagger$.
\EndClaim

Now let us fix a regular neighborhood $R$ of $H$ in $N$. There is a
homeomorphism $\rho:H\times[0,2]\to R$ such that $\rho(x,2)=x$ for
every $x\in\partial H$, and $\rho(H\times\{0\})=\Fr_NR$. Set $R^\maybeast =\rho(H\times[1,2])$
and $N^\maybeast=N-\overline{R^\maybeast}$, and denote by $\ell:N^\maybeast\to N$ the
homeomorphism which is the identity on $\overline{N-R}$ and maps
$\rho(x,t)$ to $\rho(x,2t)$ for $(x,t)\in H\times[0,1]$.

We define parametrized arcs $\alpha^\maybeast$ and $\beta^\maybeast$ in $N^\maybeast$ and $R^\maybeast$
respectively by setting $\alpha^\maybeast=\ell^{-1}\circ\alpha$ and setting
$\beta^\maybeast(t)=\rho(\beta(t),t+1)$ for $0\le t\le1$. Then
$\delta^\maybeast\doteq\alpha^\maybeast\star\beta^\maybeast$ is a well-defined parametrized arc
in $N$, and it follows from \ref{upstairs downstairs} that
\Claim\label{mostly downstairs}
The parametrized arc $\delta^\maybeast$ is
fixed-endpoint homotopic in $N$ to  $\alpha^\dagger$.
\EndClaim

Since $\alpha([0,1])\subset N^*$ and $\beta([0,1])\subset R^*$, and since the expression defining $\beta^\maybeast(t)$ has a monotone increasing
function in its second coordinate, we have:
\Claim\label{new isotopy part}
The parametrized arc $\delta^\maybeast$ is isotopic to $\alpha$ by an ambient
isotopy of $N$ which is constant on $C_0$ and on $V$.
\EndClaim

According to \ref{new isotopy part}, there is a homeomorphism $h:N\to N$,
isotopic to the identity by an 
isotopy  which is constant on $C_0$ and on $V$, such that
$h\circ\alpha=\delta^\maybeast$. 
Then $\eta^\maybeast\doteq h\circ
\eta:S^1\times[0,1] \to N$ is an embedding, isotopic to $\eta$.
We have
$\eta^\maybeast\circ\alpha^\std =\alpha^\maybeast$. 

For $j=0,1$, let $\gamma_j^\std$ denote the parametrization
$s\mapsto(e^{2\pi i s},j)$ of $S^1\times\{j\}\subset S^1\times[0,1]$
by the unit interval. We observe that the composition
$\overline{\alpha^\std }\star\gamma_0^\std\star\alpha^\std $ (where
as usual $\overline{\alpha^\std }$ denotes the path defined by
$\overline{\alpha^\std }(t)=\alpha^\std(1-t)$) is a loop based at
$(o,1)$, well defined up to  orientation-preserving
reparametrization of the domain
interval, and is based homotopic to $\gamma_1^\std$. Hence if
we set $\gamma_j^\maybeast=\eta^\maybeast\circ\gamma_j^\std$ and
$\gamma_j^\dagger=\eta^\dagger\circ\gamma_j^\std$ for $j=0,1$, then
$\gamma_1^\maybeast$ and $\gamma_1^\dagger$ are respectively based homotopic to
$\overline{\delta^\maybeast }\star\gamma_0^\maybeast\star\delta^\maybeast$ 
and to
$\overline{\alpha^\dagger }\star\gamma_0^\maybeast\star\alpha^\dagger$. But it
follows from \ref{mostly downstairs} that
$\overline{\delta^\maybeast }\star\gamma_0^\maybeast\star\delta^\maybeast$ 
and 
$\overline{\alpha^\dagger }\star\gamma_0^\maybeast\star\alpha^\dagger$  are
based homotopic to each other. 
Hence 
$\gamma_1^\maybeast$ is based homotopic to $\gamma_1^\dagger$ in $N$. 
Now since $\eta^\maybeast$ and $\eta^\dagger$ are boundary-preserving embeddings, and both map
$S^1\times\{0\}$ to the simple closed curve $C_0$, both $\gamma_1^\maybeast$
and $\gamma_1^\dagger$ are simple loops in $(\inter H)-C_0$.  The
$\pi_1$-injectivity of $\scrA$ implies that $C_0$ is $\pi_1$-injective
in $H$;  hence $(\inter H)-C_0$ is $\pi_1$-injective in the horizontal boundary
$H$ of the $I$-bundle $N$, and therefore in $N$ itself. It now follows
that $\gamma_1^\maybeast$ is based homotopic to $\gamma_1^\dagger$ in $(\inter
H)-C_0$. 


We now apply \cite[Theorem 4.1]{Epstein}, which is a fixed-basepoint version of a theorem due to
Baer (see \cite[Theorem 2.1]{Epstein}). Theorem 4.1 of \cite{Epstein}
implies
that if two homotopically non-trivial simple loops in the
interior
 of an
orientable surface are
based homotopic, then they are isotopic
by an ambient isotopy of the surface which is constant on the common base point of
the loops, and is also constant outside a compact subset of the
surface. 
In the
present context, this implies that there is a self-homeomorphism
$h^\maybeast_1$ of $(\inter H)-C_0$ which is isotopic to the identity by an ambient isotopy that is
constant on the point $(o,1)$, is also constant outside a compact
subset of $(\inter H)-C_0$,  and satisfies $h_1^\maybeast\circ\gamma_1^\maybeast=\gamma_1^\dagger$.
We may extend $h_1^\maybeast$ to a homeomorphism $h^\maybeast:N\to N$ which
is isotopic to the identity by an isotopy that is constant on $V\cup
C_0\cup \{(o,1)\}$.

We set 
$\eta\maybeprime=h^\maybeast\circ h\circ\eta$, so that $\eta\maybeprime$ is an embedding of
$S^1\times[0,1]$ in $N$ and is isotopic to $\eta$ by an isotopy which
is constant on $V$. Hence
$\scrA\maybeprime\doteq\eta\maybeprime(S^1\times[0,1])$ is an annulus isotopic to
$\scrA$, again by an isotopy which
is constant on $V$. It follows from the construction that $\eta\maybeprime$ agrees with
$\eta^\dagger$ on $S^1\times\{0\}$ and on  $S^1\times\{1\}$. 
The
construction also gives that $\eta\maybeprime\circ\alpha^\std$ is
fixed-endpoint isotopic to $\delta^\maybeast$, which by \ref{mostly downstairs}
is in turn
fixed-endpoint homotopic to  $\alpha^\dagger$. Hence 
the restriction of $\eta\maybeprime$ to
$K\doteq (S^1\times\{0,1\})\cup (\{o\}\times[0,1])\subset
S^1\times[0,1]$ 
%
is homotopic rel
$S^1\times\{0,1\} $ to $\eta^\dagger|K$.




But $S^1\times[0,1]$ is obtained from its subpolyhedron
$K$ by attaching a $2$-cell, and $N$ is
aspherical by \ref{it's still Thursday}. Hence 
$\eta\maybeprime$ is
homotopic rel $S^1\times\{0,1\} $ 
to $\eta^\dagger$. Hence there is a homeomorphism of $\scrA\maybeprime$ onto
$\scrA^\dagger$ which, regarded as a map of $\scrA\maybeprime$ into $N$, is
homotopic rel $\partial\scrA\maybeprime$ to the inclusion map.
As
$\scrA^\dagger$ and $\scrA\maybeprime$ are in particular $\pi_1$-injective,
properly embedded, compact, connected orientable $2$-manifolds in the
orientable $3$-manifold $N$, it now follows from \cite[Corollary
5.5]{waldhausen} that $\scrA\maybeprime$ is isotopic to $\scrA^\dagger$ by an
ambient isotopy which is constant on $\partial N$. As we
have observed that $\scrA$ and $\scrA\maybeprime$ are isotopic by an isotopy
 constant on $V$, we deduce that
$\scrA$ is isotopic to $\scrA^\dagger$ by an
 isotopy constant on $V$, as required for the proof in
this subcase. This completes the proof of the proposition in the case where $\scrA$
is connected.

We now turn to the general case. We use induction on the number $n$ of
components of $\scrA$. The case $n=0$ is trivial, and the case $n=1$
is the  case where $\scrA$ is connected, which has already been proved.

Now suppose that $n>1$, and that the proposition is true under the
assumption that $\scrA$ has $n-1$ components. We shall show
that if the hypotheses of the proposition hold,  if $\scrA$ has $n$
components, and if none of these components is
boundary-parallel in $N-V$,  then $\scrA$ is isotopic, by an isotopy of $N$
which is constant on $V$, to a vertical submanifold of $N$. 

Choose a component $A$ of $\scrA$, and set $\scrA_-=\scrA-A$. Then the
hypotheses of the proposition hold when $\scrA$ is replaced by
$\scrA_-$; furthermore, $\scrA_-$ has $n-1$ components, and none of these components is
boundary-parallel in $N-V$. Hence, in view of the induction
hypothesis, we may assume after an isotopy which is constant on $V$ that $\scrA_-$ is vertical
in $N$.

Let $N^{\rm spl}$ denote the manifold obtained by splitting $N$ along
$\scrA_-$, and let $\iota:N^{\rm spl}\to N$ denote the canonical surjection. Since $\scrA_-$ is vertical, $N^{\rm spl}$ inherits an $I$-bundle
structure from $N$; the vertical boundary of $N^{\rm spl}$ is $V^{\rm spl}\doteq
\iota^{-1}(V\cup \scrA_-)$. The manifold $A^{\rm spl}\doteq\iota^{-1}(A)$ is an
annulus which is mapped homeomorphically onto $A$ by $\iota$. Let $N^{\rm spl}_0$
denote the component of $N^{\rm spl}$ containing $A^{\rm spl}$. Then $N^{\rm spl}_0$, as a
component of $N^{\rm spl}$, acquires an $I$-bundle structure; the vertical
boundary $V^{\rm spl}_0$ of $N^{\rm spl}$ is a union of components of $V^{\rm spl}$. The
hypotheses of the proposition hold with $N^{\rm spl}_0$, $V^{\rm spl}_0$ and $A^{\rm spl}$
playing the roles of $N$, $V$ and $\scrA$. Since $A^{\rm spl}$ is connected,
it follows from the case of the proposition which has already been proved
that $A^{\rm spl}$ is isotopic, by an isotopy of $N^{\rm spl}_0$ which is constant on
$V^{\rm spl}_0$, to a vertical annulus $B$ in $N^{\rm spl}_0$. Now
$\iota(B)\cup\scrA_-$ is the required vertical submanifold of $N$ which is
isotopic to
$\scrA$  by an isotopy of $N$ that is constant on $V$. 
\EndProof


\Corollary\label{holiday-corollary}
Let $N$ be an orientable $3$-manifold which is an $I$-bundle over a
compact, connected $2$-manifold. Let $V$ denote the vertical boundary
of $N$. Suppose that $N'$ is a compact three-dimensional submanifold of $N$ such that
every component of $\Fr_NN'$ is a $\pi_1$-injective, properly embedded annulus in
$N-V$. (Note that $N'$ may contain one or more components
of $V$.) Then either $\Fr_NN'$ has a component which is
boundary-parallel in $N-V$, or $N'$ is isotopic, by an isotopy of $N$
which is constant on $V$, to a vertical submanifold of $N$.
\EndCorollary

\Proof
After possibly modifying $N'$ by a small non-ambient isotopy, we may
assume that $N'\cap
V=\emptyset$. Under this assumption, the result follows upon applying
Proposition 
\ref{holiday-proposition}
to $\scrA\doteq\Fr_NN'$, and observing that a three-dimensional submanifold
of $N$, whose frontier is a properly embedded $2$-manifold, is
vertical if and only if its frontier is vertical.
\EndProof

\Number\label{SFS}
The theory of Seifert fibered spaces is presented, for example, in 
Chapter 12 of \cite{hempel}.
In this paper, ``Seifert fibered space'' means a compact $3$-manifold,
not necessarily
connected, equipped with a specific Seifert fibration. By analogy with
the terminology for $I$-bundles, we shall say that
a submanifold of a Seifert fibered space is {\it vertical}
if it is a union of fibers.



\EndNumber

\Proposition\label{curse of the drinking class} 
Let $K$ be a solid torus equipped with a Seifert fibration,
 and let $F$ be a compact, non-empty, two-dimensional  submanifold of
 $\partial K$ which is vertical in $K$. Let $\scrA$ be a properly embedded $2$-manifold in $K$ whose
components are $\pi_1$-injective annuli, and suppose that $\scrA\cap
F=\emptyset$. Then $\scrA$ is isotopic, by an isotopy of $K$
which is constant on $F$, to a vertical submanifold of $K$.
\EndProposition

\Proof
Observe that for any component $A$ of $\scrA$, since
$A$ is a $\pi_1$-injective annulus in the solid torus $K$, the inclusion map
$H_1(A,\QQ)\to H_1(K,\QQ)$ is an isomorphism. It follows that each
component of $\scrA$
separates $K$.

We first consider the case in which $\scrA$ is a single annulus. Since
the compact $2$-manifold $F\subset\partial K$ is vertical in $K$ and
is disjoint from the non-empty, $\pi_1$-injective $1$-manifold
$\scrA\subset\partial K$, the components of $F$ are annuli whose
cores are parallel in $\partial K$ to the components of $\partial
\scrA$. Hence, after modifying $\scrA$ by an isotopy which is constant
on $F$, we may assume that $\scrA$ is vertical.

Since $\scrA$ separates $K$, we may write $K=X_1\cup X_2$, where the $X_i$
are submanifolds of $K$ such that $X_1\cap X_2=\scrA$. Choosing a base
point in $\scrA$, we can then exhibit $\pi_1(K)$ as a free product
with amalgamation $G_1\star_C G_2$, where $G_i$ and $C$ represent the
respective images of $\pi_1(X_i)$ and $\pi_1(\scrA)$ under
inclusion. Since $\pi_1(K)$ is infinite cyclic, one of the $G_i$ must
be equal to $C$; after re-labeling we may assume $G_1=C$, so that the
inclusion homomorphism $\pi_1(\scrA)\to\pi_1(X_1)$ is an
isomorphism. It then follows from
\cite[Theorem 10.2]{hempel}
that:
\Claim\label{parallel thing}
The pair $(X_1,\scrA)$ is homeomorphic to $(D^2\times
S^1,\alpha\times S^1)$, where $\alpha\subset\partial D^2$ is an arc.
\EndClaim

For $i=1,2$, set $L_i= X_i\cap\partial K$ and $F_i=F\cap X_i$. Then the
$X_i$ and $F_i$ are vertical submanifolds of $\partial K$, and
$F_i\subset\inter L_i$ for $i=1,2$. (Either of the $F_i$ may be empty,
although they cannot both be empty since $F\ne\emptyset$.) Let us fix a vertical
annulus $L_1'\subset\inter L_1$ which contains $F_1$. It follows from \ref{parallel thing}
that $X_1$ is a regular neighborhood of $L_1'$ in $K$. But since $L_1'$ is
vertical it has a regular neighborhood $N$ in $K$ which is vertical in
$K$ and which, like $X_1$, is disjoint from $F_2$. By the
uniqueness of regular neighborhoods, $X_1$ is carried onto $N$ by an
isotopy of $K$ which is constant on $F=F_1\cup F_2$. This isotopy in particular
carries $\scrA=\Fr_KX_1$ onto the vertical annulus $\Fr_KN$, which
gives the conclusion of the proposition in the case where $\scrA$ is connected.

To prove the proposition in general, we use induction on the number of
components of $\scrA$; if this number is $0$ the result is trivial,
and if it is $1$ then $\scrA$ is connected, in which case the result
has been proved. Now suppose that $\scrA$ has $n$ components, where
$n>1$, and that the result is true for the case where the number of
components is less than $n$. Choose a component $A$ of $\scrA$. By the
connected case of the proposition, we may assume after an isotopy that
$A$ is vertical. Since we have seen that $A$ separates $K$, we may write $K=X_1\cup X_2$, where the $X_i$
are submanifolds of $K$ such that $X_1\cap X_2=A$. Since $A$ is
vertical, the $X_i$ are vertical solid tori and inherit Seifert
fibrations from $K$. For each $i\in\{1,2\}$, set $\Phi_i=(F\cap
X_i)\cup A$ and $\scrA_i=(\scrA\cap X_i)-A$. Then the hypotheses of the
proposition  hold with $K_i, \Phi_i$ and $\scrA_i$ in place of $X$,
$F$ and $\scrA$. Since $\scrA_i$ has fewer than $n$ components,
$\scrA_i$ is isotopic, by an isotopy of $X_i$ that is constant on
$\Phi_i$, to a vertical submanifold of $X_i$. Since this holds for
$i=1,2$, the conclusion follows.
\EndProof

The following result on $3$-manifolds will be needed in Section
\ref{pushing section}:

\Proposition\label{it goes down good}
Any orientable $3$-manifold which is covered by a solid torus is
itself a solid torus.
\EndProposition

\Proof
Suppose that a $3$-manifold $M$ has a covering space $\tM$ which is a
solid torus. Since $\tM$ is compact and irreducible and has non-empty
boundary, $M$ has these properties as well. Since $M$ is irreducible
and has non-empty boundary, it follows from the Sphere Theorem that
$M$ is aspherical, and hence that $\pi_1(M)$ is torsion-free. But
$\pi_1(M)$ is virtually cyclic since $\pi_1(\tM)$ is cyclic, and every
virtually cyclic, torsion-free group is infinite cyclic. It is a
special case of \cite[Theorem 5.2]{hempel} that a compact, orientable, irreducible
$3$-manifold with non-empty boundary and infinite cyclic fundamental group is a solid torus.
\EndProof

The following fact from homotopy theory will be needed in Section
\ref{pushing section}.

\Proposition\label{from-brown-book}
Let $(X,A)$ be a topological pair with the homotopy extension
property, and suppose that $X$ and $A$ are compact. 
Let $Q$ be a Hausdorff space, and let $B$ be a closed subspace  of
$Q$. 
Suppose that $h : X \to Q$ is a
continuous map such that
$h^{-1}(B)=A$, that $h|A:A\to B$ is a homotopy equivalence, and
that $h|(X-A):X-A\to Q-B$ 
is a
homeomorphism. Then 
$h : X \to Q$ is a homotopy equivalence.
\EndProposition

\Proof
This is Assertion
7.5.8 on page 298 of \cite{brown-book}.
\EndProof


\section{The characteristic submanifold}\label{suction}

\Number\label{real pairs}
By a {\it PL pair} we mean an ordered pair $(X,Y)$, where $X$ is a PL
space and $Y$ is a closed PL subset of $X$. 
A PL pair $(X,Y)$ is said to
be {\it compact} or {\it connected} if $X$ is compact or connected,
respectively.
A {\it component} of a PL pair $(X,Y)$ is a PL pair of the form
$(L,Y\cap L)$ where $L$ is a component of $X$.
Two PL pairs $(X,Y)$ and
$(X',Y')$ are said to be {\it PL homeomorphic} if there is a PL
homeomorphism from $X$ to $X'$ mapping $Y$ onto $Y'$; such a
homeomorphism is called a {\it PL homeomorphism} between the given pairs.

We define a {\it $3$-manifold pair} to be a PL pair $(X,F)$,
where $X$ is a 
(PL)
 $3$-manifold and $F$ is a $2$-manifold which is a
closed 
(PL)
subset of $\partial X$. A $3$-manifold pair $(X,F)$ is 
said to
be
{\it compact} or
  {\it orientable} if $X$ is 
compact or orientable, respectively.

We observe that in a situation where $X$ is a compact $3$-manifold and $Y$ is a compact
three-dimensional submanifold of $X$ such that $Y\cap\partial X$ is a
$2$-manifold,  there are two $3$-manifold pairs that arise naturally:
the pair $(Y,Y\cap\partial X)$ and the pair $(Y,\Fr_XY)$. Both of
these constructions of pairs will appear in this section; for example,
the first appears in Subsection \ref{from MS, GS version}, and the
second appears in Proposition \ref{who knows}.
\EndNumber

In this section we will use the following result, which  will be proved, in a slightly more general form, in \cite{squeeze}:

\begin{squeezingtheorem}\label{one and only}
Let $X$ be a compact, orientable, irreducible, boundary-irreducible
$3$-manifold, with $\partial X\ne\emptyset$. Let $Z$ and $Z'$
be compact (but possibly disconnected) $3$-dimensional submanifolds
of $X$ such that $\Fr_XZ$ and $\Fr_XZ'$ are properly embedded,
$\pi_1$-injective $2$-manifolds in $X$. Suppose that the inclusion map
$(Z,Z\cap\partial X)\hookrightarrow(X,\partial X)$ is homotopic as a
map of pairs to a map $j$ such that $j(Z)\subset Z'$ (and therefore
$j(Z\cap\partial X)\subset Z'\cap\partial X$). Then there is a
homeomorphism $J:X\to X$, isotopic to the identity, such that
$J(Z)\subset Z'$.

If in addition we have $Z\cap Z'=\emptyset$ and $Z$ is connected, then $Z$ is a thickened
surface, and a core of $Z$ is parallel to a frontier component of $Z'$.
\end{squeezingtheorem}


\Definition\label{redundant}
Let $Y$ be a 
three-dimensional 
submanifold of a
boundary-irreducible, orientable $3$-manifold $X$ such that each component of $\Fr_XY$
is a properly embedded
surface
in $X$. We will say that a component
$K$ of $Y$ is {\it redundant} if $K$ is a thickened surface in $X$ and
has a core which is parallel in $X$ to a component of $\Fr_X(Y-K)$.
\EndDefinition

\Lemma\label{about redundancy}
Let $Y$ be a
compact, 
three-dimensional
submanifold of a
boundary-irreducible, orientable $3$-manifold $X$ such that 
 $\Fr_XY$
is a properly embedded, $\pi_1$-injective $2$-manifold in $X$. 
Then for any component $K$ of $Y$, the
following conditions are equivalent.
\begin{enumerate}
\item $K$ is a redundant component of $Y$.
\item $K$ is a thickened surface in $X$, and some
component of $\overline{X-Y}$ is a parallelism between a component of
$\Fr_XK$ and a component of $\Fr_X(Y-K)$.
\item The inclusion map $(K,K\cap\partial X)\hookrightarrow
  (X,\partial X)$ is homotopic as a map of pairs to a map
  whose image is contained in $Y-K$.
\end{enumerate}
\EndLemma

\Proof
The implications (2) $\Rightarrow$ (1) $\Rightarrow$ (3) are
obvious. The implication (3) $\Rightarrow$ (2) follows from
  the second assertion of
Squeezing
Theorem \ref{one and only}.
\EndProof

\Number\label{from MS, GS version}
The theory of characteristic submanifolds was developed by Johansson
\cite{Jo} and Jaco-Shalen~\cite{JS}, In this subsection we shall
review the essentials of the theory from the point of view of
\cite{JS}. In the following discussion, the terms ``well-embedded,''
``perfectly embedded,'' ``$S^1$-pair'' and ``$I$-pair,'' which are used
in \cite{JS}, are replaced by their definitions.

A {\it Seifert pair} 
is a compact, orientable $3$-manifold pair $(\boldW, \boldE)$ such
that for each component $\reddestK$
of $\boldW$, either
\begin{enumerate}[(a)]
\item $\reddestK$ is Seifert-fibered and $\reddestK\cap\boldE$ is a vertical
  $2$-manifold in $\partial \reddestK$, or
\item $\reddestK$ is an $I$-bundle over a surface and $\reddestK\cap\boldE$ is the associated
$\partial I$-bundle.
\end{enumerate}
A Seifert pair $(\boldW, \boldE)$,  is {\it non-degenerate} if $\boldW$ has no 
component $\reddestK$ such that $\pi_1(\reddestK)=\{1\}$, or $\pi_1(\reddestK)\cong\ZZ$ and
$\reddestK\cap\boldE=\emptyset$.

Let $\reddestX$ be a compact, orientable $3$-manifold which is irreducible and
boundary-irreducible, and has non-empty boundary. (This is a slightly
more restrictive hypothesis than the one used in \cite{JS}.)
 A map $f$ from a Seifert pair $ (\boldW, \boldE)$ to  $(\reddestX,\partial
 \reddestX)$ is said to be {\it essential} if for every
component $\reddestK$ of $\boldW$, the map $(f|\reddestK)_\sharp:
\pi_1(\reddestK)\to\pi_1(\reddestX)$  is injective, and $f|\reddestK:(\reddestK,\reddestK\cap\boldE)\to
(\reddestX, \partial \reddestX)$ 
cannot be homotoped as a map of pairs
to a map into $(\partial X,\partial X)$. 
The  Characteristic Submanifold Theorem, in a somewhat more special
form than is proved in \cite{JS}, asserts
that
up to isotopy
there is a unique
submanifold $\boldV$ of $\reddestX$
satisfying the following conditions:
\begin{enumerate}
\item $(\boldV,\boldV\cap\partial \reddestX)$ is a non-degenerate Seifert
  pair, and $\Fr_X\boldV$ is $\pi_1$-injective in $X$ and has no
  boundary-parallel component;
\item every essential map of an arbitrary non-degenerate  Seifert pair into
  $(\reddestX,\partial \reddestX)$ is
homotopic (as a map of pairs) to a map into $(\boldV, \boldV\cap\partial
\reddestX)$; and 
\item $\boldV$ has no component $\reddestK$ such that the inclusion
$(\reddestK,\reddestK\cap\partial \reddestX)\to (\reddestX, \partial \reddestX)$ is
homotopic (as a map of pairs) to a map into $(\boldV-\reddestK, (\boldV-\reddestK)\cap\partial \reddestX)$.
\end{enumerate}

Note that in view of Lemma \ref{about redundancy},
Condition (3) above is equivalent to the condition that
$\boldV$ has no redundant component. 

In
this paper, we shall define a {\it characteristic submanifold} of $\reddestX$
to be a submanifold $\boldV$ satisfying Conditions (1)--(3). (In the language of
\cite{JS}, the pair $(\boldV,\boldV\cap\partial \reddestX)$ is a
characteristic pair of the pair $(X,\partial X)$.)
\EndNumber

\Number\label{our thick and thin}
Let $\boldU$ be a compact submanifold of a $3$-manifold $\reddestX$ such
that each component of $\Fr_\reddestX\boldU$ is a properly embedded surface in $\reddestX$. 
We shall denote by $\Theta_\boldU$ the
union of all 
components of $\boldU$ that are not thickened surfaces in $X$.
\EndNumber

\Lemma\label{hare today}
Let $\reddestX$ be a compact, orientable $3$-manifold which is irreducible and
boundary-irreducible, and has non-empty boundary. Let $\boldV$ denote the
characteristic submanifold of $\reddestX$, and suppose that $\boldU$ is a submanifold of
$\reddestX$ such that Conditions  (1) and (2) in the definition of a characteristic submanifold given in
\ref{from MS, GS version} hold with  $\boldU$ in place of $\boldV$.
Then
$\Theta_\boldU$ is isotopic to $\Theta_{\boldV}$.
\EndLemma

\Proof
If, in addition to Conditions  (1) and (2) in the definition of a characteristic submanifold, Condition
(3) also holds, then $\boldU$  is a
characteristic submanifold of $\reddestX$, 
and is therefore  isotopic to 
$\boldV$. If Condition (3) does not hold, then
there is a component $R$ of
$\boldU$ such that the inclusion $(\boldU,\boldU\cap\partial \reddestX)\hookrightarrow(\reddestX,\partial \reddestX)$ is
homotopic to a map into $(\boldU- R, (\boldU- R)\cap\partial
\reddestX)$. 
  It is immediate that Conditions (1) and (2) hold with
  $\boldU_1\doteq\boldU- R$ in place of $\boldV$.
Furthermore, it
follows from 
Lemma \ref{about redundancy} that $R$ is a redundant component of
$\boldU$. In particular
$R$ is a thickened annulus or torus, so that $\Theta_{\boldU_1}=\Theta_\boldU$. Since in addition
$\boldU_1$ has fewer components than $\boldU$, it follows by induction that
there is a 
submanifold
 $\boldU'$ of $\reddestX$ such that 
$\Theta_{\boldU'}=\Theta_\boldU$, and such that Conditions (1)---(3)
hold with $\boldU'$ in place of $\boldV$. These conditions say that
$\boldU'$ is a characteristic submanifold, 
and hence $\boldU'$
is isotopic to $\boldV$. This immediately implies the conclusion.
\EndProof


\DefinitionRemarks\label{acylindrical}
We define an {\it acylindrical pair} to be a compact, orientable
$3$-manifold pair $(C,\redcale)$ 
satisfying the following conditions: 
\begin{enumerate}
\item each component of $C$ is irreducible, no
  component is a ball, and no component has a fundamental group with a
  rank-$2$ free abelian subgroup;
\item each component of $\redcale$ is an annulus and is $\pi_1$-injective
  in $C$;
\item no component of $(C,\redcale)$ is 
homeomorphic to a pair of the form $(Y,V)$, where $Y$ is an $I$-bundle
over a (possibly non-orientable) surface and $V$ is the vertical
boundary of the $I$-bundle $Y$;
\item no component of $C$ is a solid torus;
and
\item for every properly embedded, $\pi_1$-injective annulus
  $\Delta$ in $ C-\cale$, either
(a) $\Delta$ is boundary-parallel in $C-\cale$, or
(b) there is a component $E$ of $\cale$ such that  $\Delta$ is
parallel to $E$ in $(C-\cale)\cup E$ (see \ref{spanish onion}).
\end{enumerate}

It follows from Condition (1) of this definition that if $(C,\redcale)$
is an acylindrical pair, then no component of $\partial C$ is a
sphere. In addition, we assert that no component of $\partial C$ is a
torus. Indeed, if $\partial C$ has a torus component, it  follows
from Condition (1) that some component of $C$ is a solid torus; 
this contradicts Condition (4).

If $(C,\redcale)$ is an acylindrical pair, it follows from Condition (2) of the definition that no component of $(\partial
C)-\inter\redcale$ is a disk.
We also assert that  
no component of $(\partial
C)-\inter\redcale$ is an annulus. Indeed, suppose that $F$ is an annulus
component of $(\partial
C)-\inter\redcale$. 
Since we have seen that no component of $\partial C$
is a torus, the boundary curves of $F$ are contained in distinct
components $A$ and $A'$ of $\redcale$. Now $A\cup F\cup A'$ is an annulus
in $\partial C$. If $N_0$ denotes  the frontier of a small regular neighborhood of $A\cup
F\cup A'$ in $C$, then $\Delta\doteq \Fr_CN_0$  is a properly embedded annulus in $C$, disjoint
from $\redcale$. Condition (4), applied with this choice of $\Delta$,
gives a submanifold $N$ of $C$ which must be either $N_0$ or
$N_1\doteq \overline{C-N_0}$. But $N_1$ is not a solid torus since
$\partial C$ has no torus components; and the intersection of the solid
torus $N_0$ with $\partial C$ contains two distinct components of
$\redcale$. This contradiction establishes our assertion.

Thus for every
acylindrical pair
$(C,\redcale)$, every component of $(\partial
C)-\inter\redcale$ has strictly negative Euler characteristic.

\EndDefinitionRemarks

\Definition\label{essential annulus}
We define an {\it essential annulus} in a boundary-irreducible
$3$-manifold $M$ to be a properly embedded annulus in $M$ which is
$\pi_1$-injective and is not boundary-parallel (see \ref{spanish
  onion}).
\EndDefinition

\DefinitionsRemarks\label{characteristic stuff}
A submanifold $K$ of a boundary-irreducible $3$-manifold $X$ is called
{\it \bindinglike} if 
$K$ is a solid torus, 
$K\cap\partial X\ne\emptyset$,
and each component of $K\cap\partial X$ is an
annulus which is $\pi_1$-injective in $K$.

A submanifold $K$ of an orientable $3$-manifold $X$ is called
{\it \pagelike} if 
$K$ can be given the
structure of an $I$-bundle over a 
compact, connected
$2$-manifold in such a way that
$K\cap\partial X$ is the horizontal boundary of $K$.

Note that if $K$ is a \pagelike\ or \bindinglike\ submanifold of a
boundary-irreducible, orientable $3$-manifold $X$ then each component of $\Fr_XK$
is a properly embedded annulus in $X$. We will say that $K$ has an
{\it essential frontier} (relative to $X$) if each component of
$\Fr_XK$ is an essential annulus in $X$.

\EndDefinitionsRemarks

\DefinitionsRemarks\label{purely}
Let $X$ be a boundary-irreducible $3$-manifold, and let $K$ be a
\bindinglike\ submanifold of $X$. Then $K$ is \pagelike\ if and only if $K$ can be given the
structure of an $I$-bundle over an annulus or M\"obius band in such a way that
$K\cap\partial X$ is the horizontal boundary of $K$. Note also that a
\pagelike\ submanifold of $X$ has Euler characteristic $0$ if and only
if it is \bindinglike.

We define a {\it purely \bindinglike} submanifold of  $X$ to be a \bindinglike\ submanifold
which is not \pagelike.
We define a {\it purely \pagelike} submanifold of  $X$ to be a \pagelike\ submanifold
which is not bindinglike. 
\EndDefinitionsRemarks

The rest of this section will be devoted to the proof of the following result.

\Proposition\label{who knows}
Let $X$ be a simple $3$-manifold with non-empty boundary, and let $\Sigma$ be a
compact three-dimensional submanifold of $X$. Then $\Sigma$ is
isotopic to the characteristic submanifold of
$X$ if and only if it has the following properties.
\begin{enumerate}
\item 
Each component
of $\Sigma$ 
is either \bindinglike\ or \pagelike, and has essential frontier
relative to $X$.
\item
For each component $\maybeR$
of $\overline{X-\Sigma}$, 
either the pair $(\maybeR,\Fr_X\maybeR)$  is acylindrical,
or $\maybeR$ is a \ta;
and in the latter case,  the components of $\Fr_XR$ are contained in
distinct components of $\Sigma$, of which one is purely
\bindinglike\ and the other is purely \pagelike.
\end{enumerate}
\EndProposition

Lemmas 
\ref{if not why not}--\ref{thursday lemma} below are needed for the
proof of Proposition \ref{who knows}.

\Lemma\label{if not why not}
Let $X$ be a simple manifold, let $\Sigma_0$ denote its characteristic
submanifold, and let $Q$ be a submanifold of $X$ which is either
\pagelike\ or \bindinglike, has essential frontier, and is disjoint from $\Sigma_0$. Then $Q$
is a \ta\ in
$X-\Sigma_0$, and some component of $\overline{X-(\Sigma_0\cup Q)}$ is a
parallelism between a component of $\Fr_XQ$ and a component of
$\Fr_X\Sigma_0$.
\EndLemma

\Proof
Since $Q$ is pagelike or bindinglike and has essential frontier, the
pair $(Q,Q\cap\partial X)$ is a Seifert pair, and each component of
$\Fr_XQ$ is $\pi_1$-injective in $X$ and not boundary-parallel. Hence
Condition (1) of the definition of a characteristic submanifold (see
\ref{from MS, GS version}) continues to hold with
the disjoint union
$\Sigma_1\doteq\Sigma_0\cup Q$ in place of $\Sigma_0$.
Since 
$\Sigma_1\supset\Sigma_0$, 
Condition (2) of the definition of a characteristic submanifold
also continues to hold when $\Sigma_0$ is replaced
by $\Sigma_1$. If Condition (3) of the definition
were to hold with $\Sigma_0$  replaced
by $\Sigma_1$, then by the Characteristic Submanifold Theorem (see
\ref{from MS, GS version}),
$\Sigma_1$ and $\Sigma_0$ would be isotopic; this is impossible since
$\Sigma_1$ has one more component than $\Sigma_0$. Hence $\Sigma_1$
has at least one redundant component.

If $Q$ is a redundant component of
$\Sigma_1$, then Condition (2) of Lemma \ref{about redundancy} holds
with $\Sigma_1$ and $Q$ playing the respective roles of $Y$ and $K$; 
this immediately implies the conclusion of the present lemma. We shall
complete the proof by assuming that $Q$ is not a redundant component
of $\Sigma_1$, and obtaining a contradiction.

As we have seen that $\Sigma_1$ has a redundant component, we may
label the redundant components of $\Sigma_1$ as $L_1,\ldots,L_n$ for
some $n\ge1$. As we have assumed that $Q$ is not a redundant component
of $\Sigma_1$, the $L_i$ are all components of $\Sigma_0$. In view of
Lemma  \ref{about redundancy}, for $i=1,\ldots,n$, the manifold
$L_i$ is a \ta\ in $X$,
and we may choose a component $P_i$ of $\overline{X-\Sigma_1}$ which is a parallelism between a
component of $\Fr_XL_i$ and a component $F_i$ of $\Fr_X\Sigma_1$;
furthermore, we have $F_i\not\subset L_i$. For each $i$, we have a priori that either
$F_i\subset Q$ or $F_i\subset\Sigma_0$. We distinguish several cases
and subcases.

{\bf Case I.} $F_i\subset\Sigma_0$ for some $i$.

In this case, $P_i$ is a parallelism between a frontier component of
$L_i$ and a frontier component of some other component of $\Sigma_0$;
as $L_i$ is a \ta\ in
$X$, it then follows from the definition that $L_i$, regarded as a
component of $\Sigma_0$, is  redundant. This 
is a contradiction, since we observed in \ref{from MS, GS version}
that $\Sigma_0$ has no redundant component.

{\bf Case II.} $Q$ is a \ta\ in $X$, and for each $i\in\{1,\ldots,n\}$ we have
  $F_i\subset Q$.

In this case, we have in particular that   $F_1\subset Q$; the
existence of the parallelism $P_1$ then implies that a core of $Q$
is parallel in $X$ to some  component of 
  $\Fr_XL_i\subset\Fr_X\Sigma_0=\Fr_X(\Sigma_1-Q)$. By definition
  this implies that $Q$ is a redundant component of $\Sigma_1$, a
  contradiction.

{\bf Case III.} $Q$ is {\it not} a \ta\ in $X$; and for each $i\in\{1,\ldots,n\}$ we have
  $F_i\subset Q$.

  In this case, for $i=1,\ldots,n$, the parallelism $P_i$ meets $Q$ in
  the single common frontier component $F_i$, and meets
  $\Sigma_0=\Sigma_1-Q$ in its other frontier component, which is
  contained in the component $L_i$ of $\Sigma_0$. Since the components
  $L_1,\ldots,L_n$ of $\Sigma_0$ are distinct, the components
  $P_1,\ldots,P_n$ of $\overline{X-\Sigma_1}$ are also distinct. Since
    each $L_i$ is a \ta\ in $X$, and
    each $P_i$ is a parallelism, it follows that
 $Q^+\doteq Q\cup\bigcup_{i=1}^n(L_i\cup P_i)$ is a regular neighborhood of
 $Q$. Since $Q$ is either \pagelike\ or \bindinglike, 
has essential frontier,
and is not a
\ta\ in
 $X$, the manifold $Q^+$ is likewise \pagelike\ or \bindinglike,
has essential frontier,
and
is not a \ta\ in $X$.

 Set $\Sigma_2=\Sigma_1\cup(P_1\cup\cdots\cup P_n)$. Then $Q^+$ is a  component of $\Sigma_2$,
and the remaining components of $\Sigma_2$ are   the components of $\Sigma_1$
 that are distinct from $Q$ and are not among the $L_i$. The
 components of $\overline{X-\Sigma_2}$ are the
 components of $\overline{X-\Sigma_1}$ that are not among the
 $P_i$. In particular, $\Fr_X\Sigma_2$ is a union of components of $\Fr_X\Sigma_1$.

 Consider the subcase in which $\Sigma_2$ has a redundant component, say $K$.
 Then
 $K$ is a \ta\ in $X$,
 and a core of $K$ is
 parallel to
 a component $G$ of $\Fr_X(\overline{\Sigma_2-K})$. Since $Q^+$ is not a \ta\ in $X$, we have $K\ne Q^+$; hence $K$ is a
component of $\Sigma_1$
 that is distinct from $Q$ and is not among the $L_i$. Since the
 component $G$ of $\Fr_X\Sigma_2$ is in particular a component of
 $\Fr_X\Sigma_1$, it now follows that $K$ is a redundant component of
 $\Sigma_1$; by definition this means that $K$ is one of the $L_i$, a
 contradiction.

 There remains the subcase in which $\Sigma_2$ has no redundant
 component. Thus 
Condition (1) of the definition of a characteristic submanifold
 holds with $\Sigma_2$ in place of $\Sigma_0$. 
Since $Q^+$ is
 \bindinglike\ or \pagelike\
and has essential frontier,
and every component of $\Sigma_2$ other
 than $Q^+$ is a component of $\Sigma_0$, 
Condition (1) of the definition of a  characteristic submanifold 
also
 holds with $\Sigma_2$ in place of $\Sigma_0$; and so does Condition
(2), 
since $\Sigma_2\supset\Sigma_0$. Hence $\Sigma_2$ is isotopic to $\Sigma_0$.

For $i=0,2$, let $n_i$ denote the number of components of $\Sigma_i$
that are not thickened annuli in 
$X$. The components of $\Sigma_2$ are $Q^+$, which is not a \ta, and the components of
$\Sigma$ distinct from $L_1,\ldots,L_n$; each of the $L_i$ is a \ta. Hence $n_2=n_0+1$. But
since $\Sigma_2$ is isotopic to $\Sigma_0$ we have $n_2=n_0$, and we
have a contradiction in this subcase as well.
\EndProof

\Lemma\label{may 22}
Let $U$
be a compact
 submanifold of  a simple $3$-manifold $X$, and suppose that each component
of $U$ 
is either \bindinglike\ or \pagelike, and has essential frontier
relative to $X$. Let $\Sigma_0$ be a characteristic submanifold of
$X$. Then $U$ is isotopic to a submanifold of $\Sigma_0$.
\EndLemma

\Proof
As a preliminary observation, note that if $\Fr_XU$ were empty, then
$X$ would be either a solid torus or an $I$-bundle over a closed
surface; in either case we would have a contradiction to simplicity. Hence
$\Fr_XU\ne \emptyset$.

According to Squeezing Theorem \ref{one and only}, since
$\Fr_X\Sigma_0$ and $\Fr_XU$ are $\pi_1$-injective, the
conclusion of the lemma will hold provided that the inclusion map $(U,U\cap \partial
X)\hookrightarrow(X,\partial X)$ is homotopic to a map into
$(\Sigma_0,\Sigma_0\cap\partial X)$. The hypotheses of the lemma imply
that $(U,U\cap \partial
X)$ is a non-degenerate Seifert pair. Since $\Sigma_0$ is a
characteristic submanifold of $X$, Condition (2) of the definition
given in \ref{from MS, GS version} holds with $\Sigma_0$ in place of
$\Sigma$; hence we need only verify that the inclusion map $(U,U\cap \partial
X)\hookrightarrow(X,\partial X)$ is an essential map of pairs.

The $\pi_1$-injectivity of the inclusion $U\hookrightarrow X$
follows from the essentiality of $\Fr_XU$. Now suppose that $(U,U\cap \partial
X)\hookrightarrow(X,\partial X)$ is homotopic, as a map of pairs, to a map into $(\partial
X,\partial X)$. Since we have observed that $\Fr_XU\ne
\emptyset$, we may choose a component $F$ of $\Fr_XU$. The
inclusion map  $(U,U\cap \partial
X)\hookrightarrow(X,\partial X)$ is then in particular homotopic, as a
map of pairs, to a map into $(\partial
X,\partial X)$. As $X$ is irreducible and boundary-irreducible, and $F$ is
$\pi_1$-injective, it now follows by  \cite[Lemma V.3.2]{JS} that $F$
is boundary-parallel in $X$, a contradiction to the essentiality of $\Fr_XU$.
\EndProof

\Lemma\label{need to add}
Let $X$ be an orientable $3$-manifold, and let $P$ be a parallelism
between a properly embedded $2$-manifold in $X$ and a $2$-manifold in
$\partial X$. Then no essential annulus in $X$ can be a subset of $P$.
\EndLemma

\Proof
We may regard $P$ as a trivial $I$-bundle over a $2$-manifold in such
a way that the horizontal boundary of $P$ is the disjoint union of 
a surface $B_0\subset\partial X$ and a properly embedded surface
$B_1\subset X$, and the vertical boundary of $P$ is contained in
$\partial X$. If $A$ is an essential annulus in $X$ and $A\subset P$, then
after an isotopy we may assume that $\partial A\subset\inter B_0$. 
Applying
Proposition  \ref{holiday-proposition} taking  $N=P$ and 
$\scrA=A$,
we deduce that either   $A$  is
boundary-parallel in $X$, contradicting the essentiality of $A$, or
that $A$ is isotopic, by an isotopy of $P$
which is constant on the vertical boundary of $P$, to a vertical
annulus in $P$. The latter alternative is ruled out by the observation
that $A$ meets only one component of the horizontal boundary of $P$,
namely $B_0$, whereas a vertical annulus must meet both components of
the horizontal boundary.
\EndProof





\Lemma\label{red red rose}
Suppose that $K$ and $K'$ are compact, connected submanifolds of a
connected, boundary-irreducible $3$-manifold $X$. Suppose that 
$K'\subset K$, and that
every component
of $\Fr_XK'$ is an essential annulus in $X$. Then:
\begin{itemize}
\item if $K$ is a \bindinglike\ submanifold of $X$, then $K'$ is also a
  \bindinglike\ submanifold of $X$; and
\item if $K$ is a \pagelike\ submanifold of $X$, then $K'$ is also a
  \pagelike\ submanifold of $X$.
\end{itemize}
\EndLemma

\Proof
To prove the first assertion, note that since the components of
$\Fr_XK'$ are essential annuli in $X$, they are in particular
$\pi_1$-injective in $X$; hence $K'$ is $\pi_1$-injective in $X$, and
in particular it is $\pi_1$-injective in $K$. If $K$ is \bindinglike,
then $\pi_1(K)$ is infinite cyclic, and hence $\pi_1(K')$ is either
infinite cyclic or trivial. But since $X$ is boundary-irreducible and
$K$ is a solid torus, we have $K\ne X$ and hence $K'\ne X$; as $X$ is
connected, it follows that $\Fr_XK'$ is non-empty, and therefore
contains at least one $\pi_1$-injective annulus. Hence $\pi_1(K')$
is non-trivial, and hence infinite cyclic. 
But $K'$ is irreducible, since the solid torus $K$ is irreducible and
each component of $\Fr_KK'$ meets $\partial K$. Hence
by 
\cite[Theorem 5.2]{hempel},
$K'$ is a solid torus. 
Since each component of
$\Fr_XK'$ is a $\pi_1$-injective annulus in $K$', each component of
$K'\cap\partial X$ is also a $\pi_1$-injective annulus in $K'$. Hence  $K'$ is a \bindinglike\ submanifold of $X$.

To prove the second assertion, suppose that $K$ is \pagelike, and fix
an $I$-bundle structure on $K$ such that $K\cap\partial X$ is the
horizontal boundary of $K$. Let $V$ denote the vertical boundary of
$K$. 
Since $K'\subset K$ and $\Fr_XK=V$, we may assume
after a small isotopy 
that $K'\cap V=\emptyset$. Since the components of $\Fr_KK'$ are
essential annuli in $X$, they are in particular $\pi_1$-injective
annuli in $K$.
Hence we may apply 
Corollary \ref{holiday-corollary}, 
taking $N=K$
and $N'=K'$ and choosing $V$ as above, to deduce that either (a) $\Fr_KK'$ has a component which is
boundary-parallel in $K-V$, or (b) $K'$ is isotopic, by an isotopy of $K$
which is constant on $V$, to a vertical submanifold of $K$. But (a)
would imply that some component of $\Fr_XK'=\Fr_KK'$ is
boundary-parallel in $X$; this is impossible, since by hypothesis the
components of $\Fr_XK'$ are essential annuli in $X$. Hence (b) holds,
and it follows that $K'$ is a \pagelike\ submanifold
of $X$.
\EndProof

\Lemma\label{thursday lemma}
Let $X$ be a boundary-irreducible, orientable $3$-manifold, and let  $H$ be 
a 
submanifold of $X$ which is either \bindinglike\ or \pagelike, and has
essential frontier in $X$.
Let $\maybeR$ be a submanifold of $X$ such that 
 every component of $\Fr_XR$ is an essential
annulus in $X$. Suppose  that $\Fr_XR\cap\Fr_XH=\emptyset$.
Suppose also that either $(\maybeR ,\Fr_XR )$ is an
acylindrical pair, or $R$ is a \ta\ in $X$. 
Then every component of 
$R\cap H$ is a \ta\
in $X$. Furthermore, any core annulus of a component of $R\cap H$ is
parallel in $R$ to a  component of $\Fr_XR$.
%
\EndLemma


\Proof
Let a component $L$ of $\maybeR \cap H$ be given. 
We may write $\Fr_XL$ as a disjoint union  $\reallyredcala\discup\cale$,
where $\reallyredcala
$ is a union of components of $\Fr_XH$ and
$\cale
$ is a union of components of $\Fr_X\maybeR $. 
(Either $\reallyredcala$ or $\cale$ may be empty.)
In
particular, it follows from the hypothesis that each component of $\Fr_XL$
is an essential annulus in $X$. It then follows from Lemma
\ref{red red rose} that $L$ is either a 
\bindinglike\ or a \pagelike\ 
submanifold of $X$.

We shall show:
\Claim\label{cold}
There is a \ta\ 
 $P$ in $X$ with the following properties: (a) $L\subset
 P\subset H$; (b) every component of $\Fr_PL$ is a component of
 either $\cale$ or $\reallyredcala$; and (c) at least
 one component of $\Fr_XP$ is a component of $\cale$.
\EndClaim

To prove \ref{cold}, we first note that by the hypothesis of the
lemma, either  $R$ is a \ta, or
$(R,\Fr_XR)$ is an acylindrical pair. 
To establish \ref{cold} when $R$ is a \ta,
we need only set $P=H$. In this case, (b) holds vacuously, and
(c) holds for the stronger reason that each of the two
components of $\Fr_XP$ is a component of $\cale$. 

To prove \ref{cold} when $(R,\Fr_XR)$ is an acylindrical pair, we need a preliminary construction.
Suppose that a component $\redA $ of 
$\reallyredcala$ is given. Since $\redA $ is an essential annulus in $X$, it is in particular
$\pi_1$-injective  in ${\maybeR }-\Fr_XR$;
it therefore follows from Condition (5)
of the definition of an acylindrical pair (see \ref{acylindrical})
that either $\redA$ is boundary-parallel in ${\maybeR }-\Fr_XR$, or
there is a component $E$ of $\Fr_XR$ such that  $\redA$ is
parallel to $E$ in $({\maybeR }-\Fr_XR)\cup E$ (see \ref{spanish onion}).
The first of these alternatives is ruled out by the essentiality of
$\redA$ in $X$.
Hence for each component $\redA $ of 
$\reallyredcala$, we may fix a component $E_\redA$ of $\Fr_XR$ and a parallelism
$P_\redA$ in  $({\maybeR }-\Fr_XR)\cup E_\redA$ between $\redA$ and $E_\redA$.

For each component $\redA$ of $\reallyredcala$ we must have either $L\cap
P_\redA =\redA $, or   $P_\redA \supset L$. 
Consider the 
subcase 
in which $L\cap
P_\redA =\redA $  for
every component $\redA$ of $\reallyredcala$. 
(This includes the 
subcase 
in which $\reallyredcala=\emptyset$.)
Then the pair $({\maybeR },\Fr_XR)$ is
homeomorphic to $(L,\reallyredcala)$.
Since $L$ is either
\pagelike\ or \bindinglike, it then follows that either (i) $\maybeR $ is a solid torus, or
(ii) $\maybeR $ may be given the structure of an $I$-bundle in such a way
  that $\Fr_XR $ is the vertical boundary of $\maybeR $. In view of Conditions (3) and (4) of the definition
  of an acylindrical pair, either of the alternatives (i) and (ii)
  gives a contradiction. 
Hence there must be a component
  $\redA_0$ of $\reallyredcala$ such that $P_{\redA_0} \supset L$. 
(In
  particular we have $\reallyredcala\ne\emptyset$ in the case where
$(R,\Fr_XR)$ is an acylindrical pair.)

Now set $P=P_{\redA_0}$. By construction, $P$ is a \ta, 
and it has the properties (a) and (b) stated in \ref{cold}. Furthermore, the component
$E_{A_0}$ of $\Fr_XP$ is a component of $\cale$, so that $P$ has
Property (c). Thus \ref{cold} is
established in all cases.

Now fix a \ta\ $P$ having the properties (a)--(c) stated in \ref{cold}, and fix
a component $E$ of $\Fr_XP$ which is a component of $\cale$. We may 
equip $P$ with the
structure of a trivial $I$-bundle over an annulus $G$ in such a way that
$\Fr_XP$ is
the vertical
boundary of $P$. 
Property (b) of the \ta\ $P$ implies that each component of $\Fr_PL$
is
an essential
annulus in $X$; hence any such component is $\pi_1$-injective in $P$,
and is not parallel in $P$ to a component of the horizontal boundary
of $P$. It therefore follows from 
Corollary \ref{holiday-corollary} 
that the
$I$-bundle structure on $P$ may be chosen in such a way that $L$ is vertical.
Thus $L$ is an induced $I$-bundle over a connected subsurface $G'$ of
$G$. 
Since the components of  $\Fr_PL$ are  $\pi_1$-injective annuli in $P$,
the components of $\Fr_GG'$ are  $\pi_1$-injective simple closed
curves in $G$; hence $G'$ is an annulus sharing a core curve with
$G$. Finally, Property (c) of $P$ implies that $G'$ has a boundary
curve whose preimage under the $I$-fibration of $P$ is a component of
$\cale$; this implies that a core
annulus   of $L$ is
parallel in $P$, and hence in $R\supset P$,
 to a  component of $\cale$.
%
%
\EndProof

\Proof[Proof of Proposition \ref{who knows}]
Let $\Sigma_0$ denote the
characteristic submanifold of $X$ (which we fix within its isotopy class).

We first prove the ``only if'' assertion of the 
proposition. 
We must
show that $\Sigma=\Sigma_0$ has Properties
(1) and (2) 
from the statement of the proposition.

To establish Property (1), we consider an arbitrary component
$\reddestK$ of $\boldV$. According to Condition  (1) in the definition
of the characteristic submanifold (see \ref{from MS, GS version}),
$\Fr_X\reddestK$ is $\pi_1$-injective and has no boundary-parallel
component. Thus we need only show that $\reddestK$ is either
pagelike or bindinglike. Condition  (1) in the definition of the
characteristic submanifold also implies that
$(\reddestK,\reddestK\cap\partial X)$ is a nondegenerate Seifert
pair. Hence one of the alternatives (a), (b) of the
definition of a Seifert pair holds with $\reddestK\cap\partial X$ in place of 
$\reddestK\cap\boldE$. In the case in which (a) holds, it follows from
the definitions that $\reddestK$ is pagelike.

In the case in which (b) holds, 
$K$ admits a Seifert fibration. 
In particular, each component of $\partial K$ is a torus. Since $X$ is connected and 
$\partial X\ne\emptyset$, we have
$\partial K\ne\emptyset$. If we fix a torus $T\subset\partial K$, the
simplicity of $X$ implies that $T$ is not $\pi_1$-injective in $X$. But
since $\Fr_X\reddestK$ 
is $\pi_1$-injective
in $X$, the submanifold $K$ of $X$ is also $\pi_1$-injective. Hence
$T$ is not $\pi_1$-injective in $K$. The  $\pi_1$-injectivity of
$\Fr_X\reddestK$, together with the irreducibility of the simple
manifold $X$, also implies that $K$ is irreducible. But  it is a standard
consequence of the Loop Theorem that a compact, orientable,
irreducible $3$-manifold, whose boundary contains a torus which is not
$\pi_1$-injective, is a solid torus.


If the vertical
$2$-manifold  $\reddestK\cap\partial X$ were the full boundary of
the solid torus
 $\reddestK$, the
connectedness of $X$ would imply that $X=K$; this is a contradiction,
since a solid torus is not simple. This shows 
that the components of $\reddestK\cap\partial X$
are $\pi_1$-injective annuli in $\reddestK$. Thus $\reddestK$ is
bindinglike in this case. This completes the proof that $\boldV$ has
Property (1) from the statement of the proposition.

To prove
that $\boldV$ has
 Property 
(2),
we consider an arbitrary component $\maybeR $
of $\overline{X-\Sigma_0}$. 
We will first show that 
either the pair $(\maybeR,\Fr_X\maybeR)$  is acylindrical,
or $\maybeR$ is a \ta.
To show this, we assume that $R$ is not a \ta,
and show that
Conditions (1)--(5) of Definition \ref{acylindrical}
hold with $(\maybeR,\Fr_X\maybeR) $ playing the role of $(C,\cale)$. 

To establish
Conditions (3) and (4) of Definition \ref{acylindrical}, assume that 
either $\maybeR $ is a solid torus, or $(\maybeR ,\Fr_X\maybeR )$ is 
homeomorphic to a pair of the form $(Y,V)$, where $Y$ is an $I$-bundle
over a (possibly non-orientable) surface and $V$ is the vertical
boundary of the $I$-bundle $Y$. Then $\maybeR $ is a \pagelike\ or
\bindinglike\ submanifold of $X$, and has essential frontier (since
$\Sigma_0$ has essential frontier).
Hence $\maybeR $ is a regular neighborhood in $X$ of a
submanifold $Q$ of $X$ which is either \pagelike\ or \bindinglike,
and has essential frontier. Since $Q$ is disjoint from $\Sigma_0$,
it follows from Lemma \ref{if not why not} that  $Q$
is a 
thickened annulus in
$X$,
a contradiction. 

To establish Condition
(5) of Definition \ref{acylindrical}, suppose that
$\Delta$ is a
properly embedded, $\pi_1$-injective annulus
  in $ \maybeR -\Fr_X\maybeR $. If $\Delta$ is essential in $X$, then a regular
  neighborhood $Q$ of $\Delta$, regarded as a submanifold of $X$, is
  both \pagelike\ and \bindinglike, has essential frontier, and is
  disjoint from $\Sigma_0$. 
Applying Lemma \ref{if not why not} with this choice of $Q$, we deduce that
some component of $\overline{X-(\Sigma_0\cup Q)}$ is a
parallelism between a component of $\Fr_XQ$ and a component of
$\Fr_X\Sigma_0$. This implies Alternative (b) of Condition (5). 
If $\Delta$ is not essential in $X$, we may fix a parallelism
$P\subset X$
between $\Delta$ and some annulus $J\subset\partial X$. Since every
component of $\Fr_X\maybeR \subset\Fr_X\Sigma_0$ is an essential annulus in
$X$, it follows from Lemma \ref{need to add} that no component of
$\Fr_X\maybeR $ is contained in $J$; hence $J\cap\Fr_X\maybeR =\emptyset$, and
$\Delta$ is boundary-parallel in $\maybeR -\Fr_X\maybeR $. This is Alternative (a) of Condition (5).

The pair $(\maybeR ,\Fr_X\maybeR )$ satisfies
Condition (2) of Definition \ref{acylindrical} because the components
of $\Fr_X\maybeR$ are components of
of $\Fr_X\Sigma_0$, and are therefore essential annuli in $X$.
 To verify Condition (1), first note that $\maybeR $ is
irreducible because the simple manifold $X$ is irreducible, and every
component of $\Fr_X\Sigma_0$ meets $\partial X$. 
The assertion that no
component of $\maybeR $ is a ball is immediate from the definition of
simplicity if $\Sigma_0=\emptyset$, i.e. $\maybeR =X$; if
$\Sigma_0\ne\emptyset$ then $\Fr_X\maybeR $ is non-empty, so that $\maybeR $
contains at least one $\pi_1$-injective annulus and is therefore not a
ball. Finally, since $\pi_1(X)$ has no rank-$2$ free abelian subgroup
by simplicity, and $\maybeR $ is $\pi_1$-injective in $X$ by the
essentiality of the components of $\Fr_X\Sigma_0$, the group
$\pi_1(\maybeR )$ also has no rank-$2$ free abelian subgroup. 
This
completes the proof that either $(\maybeR,\Fr_X\maybeR)$  is acylindrical
or $\maybeR$ is a \ta.

To complete the proof that $\boldV$ has Property (2) from the
statement of the proposition, we shall assume that $R$ is a \ta, and
show that
the two components $F_0$ and $F_1$ of $\Fr_X\Sigma_0$ 
are contained in
distinct components of $\Sigma_0$, of 
which one is purely
\bindinglike\ and the other is purely \pagelike. 
Assume that this is false,
and let us label the components $\Fr_XR$
as $A_1$ and $A_2$. For
$i=1,2$, let $\reddestK_i$ denote the component of $\boldV$ containing $A_i$; a
priori, $\reddestK_1$ and $\reddestK_2$ may or may not be distinct. 
By our assumption, 
$\reddestK_1$ and $\reddestK_2$ are either both pagelike or
both bindinglike.

Set $K^*=\reddestK_1\cup
R\cup \reddestK_2$ (so that $K^*=\reddestK_1\cup
R$ if $\reddestK_1= \reddestK_2$).  We claim that the connected
$3$-manifold pair $(K^*,K^*\cap\partial X)$ is a
nondegenerate Seifert pair.

If $\reddestK_1$ and $\reddestK_2$ are both bindinglike, then
each $\reddestK_i$ is a solid torus having non-empty intersection with
$\partial X$, and each component of $\reddestK_i\cap\partial
  X$ is an annulus in $\partial \reddestK_i$ which is
  $\pi_1$-injective in $\reddestK_i$. Hence $A_i$ is a
  $\pi_1$-injective annulus in $K_i$ for $i=1,2$. Since $R$ is a
  thickened annulus, it follows that $K^*$ may be given the structure
  of a Seifert fibered space in such a way that $K^*\cap\partial X\ne\emptyset$ is
  vertical; thus in this case $(K^*,K^*\cap\partial X)$ is a
nondegenerate Seifert pair.

If $\reddestK_1$ and $\reddestK_2$ are  both pagelike, then
the union of the manifolds $\reddestK_1$ and $\reddestK_2$ (which are
either disjoint or identical)
may be given the structure of an $I$-bundle over a surface
 in such a way that
  $(\reddestK_1\cup\reddestK_2)\cap\partial X$ is the horizontal
  boundary of $\reddestK_1\cup\reddestK_2$. Since $R$ is a
  thickened annulus, the $I$-bundle structure on
  $\reddestK_1\cup\reddestK_2$ extends to an $I$-bundle structure on
  $K^*$, and  $K^*\cap\partial X$ is its horizontal boundary. Hence  $(K^*,K^*\cap\partial X)$ is a
 Seifert pair. To establish 
nondegeneracy in this case we need only verify that $K^*$ is
non-simply-connected; this follows from the essentiality of
$\Fr_XR$. Thus our claim is established in all cases.

If we set 
$\boldU^*=\boldV\cup R$, it now follows that 
$(\boldU^*,\boldU^*\cap\partial X)$ is a nondegenerate Seifert
pair. Furthermore, the components of $\Fr_X\boldU^*$ are components of
$\Fr_X\boldV$ and are therefore essential annuli in $X$.
Thus Condition (1) of the definition of a characteristic submanifold
(see \ref{from MS, GS version})
holds with $\boldU^*$ in place of $\boldV$.
Condition (2) of that definition is immediate
  since $\boldU^*\supset \boldV$. Hence by Lemma \ref{hare  today}, $\Theta_{\boldU^*}$ is
  isotopic to $\Theta_{\boldV}$.

If $\reddestK_1$ and $\reddestK_2$ are distinct and neither of them is
a thickened annulus, then $\boldU^*$ is
obtained from $\boldV$ by replacing two  components, neither of
which is a thickened surface, by a single component; this implies that $\Theta_{\boldU^*}$ has strictly fewer components
than $\Theta_{\boldV}$, a contradiction. If $\reddestK_1$ and $\reddestK_2$ are distinct and
at least one of them, say $\reddestK_1$, is a thickened annulus, then since $R$ is a
thickened annulus, the component $K_1$ of $\boldV$ is redundant; this
is a contradiction, since we observed in \ref{from MS, GS version}
that $\boldV$ has no redundant component.

Now suppose that
$\reddestK_1=\reddestK_2$.
Then there is a simple closed curve $\gamma$ in $\inter
K^*\subset\inter \boldU^*$ whose homological intersection number with $A_1$ in
$X$ is equal to $1$. In particular the homology class of $\gamma$ is
not in the image of the inclusion homomorphism $H_1(\Fr_XK^*)\to
H_1(K^*)$, and hence $K^*$ is not a thickened annulus. If
$\reddestK_1=\reddestK_2$ is a thickened annulus, this implies
that $\Theta_{\boldU^*}$ has one more component than $\Theta_{\boldV}$, a
contradiction. Finally, if $\reddestK_1=\reddestK_2$ is not a
thickened annulus, then $A_1$ is a component of $\Fr_X\Theta_{\boldV}$ having homological intersection number
$1$ with $\gamma$; but since $\gamma$ is contained in
$\inter \Theta_{\boldU^*}$, it has  homological intersection number $0$ with every
component of $\Fr_X\Theta_{\boldU^*}$, which is in turn isotopic to $\Fr_X\Theta_{\boldV}$. Again we have a contradiction.

This completes the proof of the ``only if'' assertion of the proposition.


We now turn to the proof of the ``if''  assertion. Suppose that
$\Sigma$ is a compact three-dimensional submanifold of $X$  which has
Properties 
(1) and (2) from the statement of the present proposition.

According to Property (1) 
from the statement of the  proposition, 
each component
of $\Sigma$ 
is either \bindinglike\ or \pagelike, and has essential frontier
relative to $X$.
Hence by Lemma \ref{may
  22}, applied with $U=\Sigma$,
$\Sigma$ is  isotopic to a submanifold of
$\Sigma_0$.
We may therefore assume without loss of generality that
$\Sigma\subset\Sigma_0-\Fr_X\Sigma_0$. 

We shall complete the proof
by showing that $\Sigma_0$ is a regular neighborhood of
$\Sigma$. Equivalently, we must show, for each component $L$ of
$\overline{\Sigma_0-\Sigma}$, that
\begin{whoknowsclaim}\label{i swear}
$L$ is a parallelism between a component of
$\Fr_X\Sigma$ and a component of $\Fr_X\Sigma_0$.
\end{whoknowsclaim}


Consider an arbitrary component $L$ of
$\overline{\Sigma_0-\Sigma}$. 
Let $\maybeR $ denote the component of $\overline{X-\Sigma}$
containing $L$, and let $H$ denote the component of $\Sigma_0$
containing $L$. Then $L$ is a component of $\maybeR \cap H$. 

Since $\Sigma$ has Property 
(2), 
either 
the pair $(\maybeR ,\Fr_X\maybeR )$  is acylindrical, or
the manifold $\maybeR $ is a 
\ta.

It therefore follows from  Lemma \ref{thursday lemma} that $L$ is a
\ta\ in $X$,
and that a core annulus of $L$ is
parallel in $R$ to a  component of $\Fr_XR$.

 Each frontier component of the 
\ta\
$L$ is contained either in $\Fr_XR\subset\Fr_X\Sigma$ or in
$\Fr_XH\subset\Fr_X\Sigma_0$. In order to establish \ref{who knows}.1, we need to show that $L$ has one
frontier component in $\Fr_XR$ and one in $\Fr_XH$.

Assume that this is false. Then either (i) $\Fr_XL\subset\Fr_XR$ or
(ii) $\Fr_XL\subset\Fr_XH$. 

If (i) holds, then $L=R\subset H$, and in particular $R$ is a 
\ta\ in $X$.
Furthermore, since $\Sigma$ has Property (2), the components of $\Fr_XR$ are contained in
distinct components $K_1$ and $K_2$  of
$\Sigma$, where $K_1$  is purely
\bindinglike\ and $K_2$ is purely \pagelike. Since $\Sigma_0$ contains
both $\Sigma$ and $R$, the connected manifold $K_1\cup R\cup K_2$ must
be contained in a component $K$ of $\Sigma_0$. By Property (1) of
$\Sigma_0$, the submanifold $K$ of $X$ is either \pagelike\ or
\bindinglike. But if $K$ is \pagelike, then by Lemma \ref{red red rose},
$K_1$ is \pagelike; this is a contradiction since 
$K_1$  is purely
\bindinglike. Likewise, if
$K$ is \bindinglike, then by Lemma \ref{red red rose},
$K_2$ is \bindinglike; this is a contradiction since 
$K_2$  is purely
\pagelike.

If (ii) holds, then $L=H\subset R$.  In particular $H$ is a 
\ta\ in $X$,
and a core annulus  $A$ of $H$  
 is
parallel in $R$ to a  component $E$ of $\Fr_XR$. Since $R$ is a
component of $\overline{X-\Sigma}$, we have
$E\subset\Fr_X\Sigma$; in particular, $E$ is contained in a component
$Y$ of $\Sigma$, which in turn is contained in a component $H'$  of
$\Sigma_0$. On the other hand, since the component $H$ of $\Sigma_0$ is
contained in $R$, no component of $\Sigma$ can be contained in $H$. Hence the
components $H$ and $H'$ of $\Sigma_0$ are distinct.
Now since $E$ and $A$ are parallel, the inclusion map $H\hookrightarrow X$ is homotopic in $X$ to a map
  whose image is contained in $H'$. It therefore follows from Lemma
  \ref{about redundancy} that $H$ is a redundant component of
  $\Sigma_0$. 
This is a contradiction, since we observed in \ref{from MS, GS version} that the characteristic
  submanifold $\boldV$ of $X$ has no redundant component. 
\EndProof

\section{
Books of $I$-bundles}\label{book section}

\DefinitionRemarks\label{books and such}
A {\it generalized book of $I$-bundles}
 is a triple $\calw =
(W,\calb,\calp)$, where $W$ is a (possibly empty and possibly disconnected) compact, orientable
$3$-manifold, and $\calb , \calp  \subset W$ are compact (possibly
disconnected) submanifolds such
that
\begin{itemize}
\item $\calb$  is equipped with a Seifert fibration over a compact (possibly disconnected) $2$-manifold;
\item $\calp$  is equipped with an $I$-fibration over a (possibly disconnected) 2-manifold, and every component of $\calp$  has non-positive Euler characteristic;
\item $W =\calb \cup\calp$;
\item $\calb  \cap \calp$  is the vertical boundary of $\calp$; and
\item $\calb  \cap \calp$  is vertical in the Seifert fibration of $\calb$.
\end{itemize}
We shall denote $W$, $\calb$  and $\calp$  by $|\calw|$, $\calb_\calw$ and $\calp_\calw$, respectively. The components of $\calb_\calw$ will be called bindings of $\calw$, and the components of $\calp_\calw$ will be called its pages. The submanifold $\calb\cap \calp$, whose components are properly embedded annuli in $W=|\calw|$, will be denoted $\cala_\calw$.

A generalized book of I-bundles $\calw$ will be termed {\it connected}
if the manifold $|\calw|$ is connected. 
A {\it component} of a generalized book of $I$-bundles $\calw
=
(W,\calb,\calp)$
 is a
(connected) generalized book of $I$-bundles of the form
$(X,X\cap\calb,X\cap\calp)$ where $X$ is a component of $W$.





The definition given above is essentially identical with the one given
in \cite[Definition 2.1]{acs-singular}. The condition that $\calb$ and
$\calp$ be equipped with a {\it specific} Seifert fibration and a
specific $I$-fibration was not emphasized in \cite{acs-singular}, but
it will be important in the present paper.


%


If 
$B$ is a binding of
a generalized book of $I$-bundles $\calw$, 
the non-singular
fibers of the Seifert-fibered $3$-manifold $B$ are simple closed curves that  all represent the same
isotopy class in $B$, 
and hence in $|\calw|$. 
We denote by
 $\calf_\calw$ the set of all simple closed curves in $|\calw|$ that are
isotopic in $|\calw|$ to non-singular fibers of bindings of $\calw$. Thus  $\calf_\calw$ is a
finite union of
isotopy classes of simple closed curves in $|\calw|$.

\Claim\label{E or A is injective}
\textnormal{
In an $I$-bundle over a $2$-manifold whose components are all of
non-positive Euler characteristic, the vertical boundary is
$\pi_1$-injective. Furthermore, in a Seifert fibered space, any vertical
annulus is $\pi_1$-injective. It follows that if $\calw$ is any
generalized book of $I$-bundles, $\redcalA_\calw$ is $\pi_1$-injective in $\calp_\calw$ and in $\calb_\calw$. Since
$\calp_\calw\cup\calb_\calw=|\calw|$ and
$\calp_\calw\cap\calb_\calw=\redcalA_\calw$, it follows that
$\calp_\calw$, $\calb_\calw$ and $\redcalA_\calw$ are all
$\pi_1$-injective in $|\calw|$.
}
\EndClaim

\Claim\label{derry down down}
\textnormal{
If $\calw$ is a generalized book of $I$-bundles, then for each
component $B$ of $\calB_\calw$, the $2$-manifold
$B\cap\partial|\calw|=\overline{(\partial B)-(\redcalA_\calw\cap\partial
  B)}$ is a disjoint union of annuli and tori, while for each
component $P$ of $\calp_\calw$, the $2$-manifold
$P\cap\partial|\calw|$ is the horizontal boundary of an $I$-bundle
over a surface of non-positive Euler characteristic, and hence
$\chi(P\cap\partial|\calw|)\le0$. It follows that each component of
$\partial|\calw|$ has non-positive Euler characteristic, i.e. no
component of $\partial|\calw|$ is a sphere.
}

{\rm
Since the Euler characteristic of a compact, orientable $3$-manifold
is one-half the Euler characteristic of its boundary, it now follows
that if $\calw$ is a generalized book of $I$-bundles, each component
of $|\calw|$ has non-positive Euler characteristic.}
\EndClaim

\Claim\label{vertical in book}
\textnormal{
If $\calw$ is a generalized book of $I$-bundles, we 
define a
{\it $\calw$-vertical submanifold
}
 of $|\calw|$ to be a submanifold, each of
whose components is
contained in  $\calb_\calw$ or $\calp_\calw$ and is vertical in the 
Seifert fibration or
$I$-fibration 
(respectively) of $\calb_\calw$ or
$\calp_\calw$. 
}

\textnormal{\noindent
We observe that if a $\calw$-vertical submanifold of $|\calw|$ is
contained in $\partial|\calw|$, then it is
contained in  $\calb_\calw$ and is vertical in the
Seifert fibration of $\calb_\calw$.
}

\textnormal{\noindent
We define a submanifold of $|\calw|$ to be {\it strongly
  $\calw$-vertical} if it is  $\calw$-vertical and is disjoint from
$\redcalA_\calw$. 
}
\EndClaim

\Claim\label{subsubmarine}
\textnormal{
Note that if $\calw$ is a generalized book of $I$-bundles, any
$\calw$-vertical submanifold of $|\calw|$ is isotopic by
a small isotopy to a strongly $\calw$-vertical submanifold. 
}

\EndClaim

\EndDefinitionRemarks

\DefinitionRemarks\label{real books}
In this paper, we define a {\it book of $I$-bundles} to be a
generalized book of $I$-bundles in which each 
binding
is a solid torus.

This definition differs from the one given in \cite{acs-singular} in
that we do not require the underlying manifold of a book of $I$-bundles to be non-empty, and
we do not require each binding to meet a page.
These differences will become crucial in Proposition \ref{new
  characterization}, which indirectly involves the notion of a book of
$I$-bundles, and would not be true with the stronger
definition given in \cite{acs-singular}.

\Claim\label{confetti}
\textnormal{
If $\calw$ is a book of $I$-bundles, each component of $|\calw|$
contains either a binding or a page of $\calw$; and if $Y$ is a
binding or page,
the intersection of $Y$ with $\inter|\calw|$ consists of the interior
of a disjoint union of closed annuli in the closed surface $\partial Y$. Hence each
component of $|\calw|$ has non-empty boundary.
}
\EndClaim

\Claim\label{substrict}
\textnormal{
A book of $I$-bundles $\calw$ is said to be {\it strict} if each of its pages
has strictly negative Euler characteristic.
}
\EndClaim

\Claim\label{derry ferry}
\textnormal{
If $\calw$ is a strict book of $I$-bundles, the argument of \ref{derry
  down down} gives stronger information. In this case, for each
component $B$ of $\calB_\calw$, the $2$-manifold
$B\cap\partial|\calw|=\overline{(\partial B)-(\redcalA_\calw\cap\partial
  B)}$ is a disjoint union of annuli and tori, while for each
component $P$ of $\calp_\calw$, the $2$-manifold
$P\cap\partial|\calw|$ is the horizontal boundary of an $I$-bundle
over a surface of strictly negative Euler characteristic, and hence
$\chi(P\cap\partial|\calw|)<0$. It follows that each component of
$\partial|\calw|$ either has strictly negative Euler characteristic or
is contained in a component of $|\calw|$ consisting of a single
binding; such a component is a solid torus.
}

{\rm
Now let $\calw$ be a connected strict book of $I$-bundles such that
$|\calw|$ is not a solid torus. By \ref{confetti} we have
$\partial|\calw|\ne\emptyset$. Since
$\chi(|\calw|)=\chi(\partial|\calw|)/2$, and each component of
$\partial|\calw|)$ has strictly negative Euler characteristic, we have
$\chi(|\calw|)<0$. 
}
\EndClaim

\EndDefinitionRemarks


\Number\label{new-gbp}
\textnormal{
A $3$-manifold pair $(W,\redcale)$ is called a {\it generalized bookish pair} if
there exists a generalized book
of $I$-bundles $\calw$ such that (a) $|\calw|=W$, and (b) each component of
$\redcale\subset W$ is a $\calw$-vertical annulus in $\partial W$. If
a generalized book of $I$-bundles $\calw$ satisfies (a) and (b), we
must have $\cale\subset\calb_\calw$; and after possibly modifying
$\calw$ by a small isotopy we can arrange that (c)  $\redcalA_\calw\cap\redcale=\emptyset$.
A generalized book of $I$-bundles $\calw$  will be said to {\it \realize}
the generalized bookish pair $(W,\redcale)$ if it satisfies (a), (b)
and (c). 
Note that Condition (c) implies that $\redcale$ is strongly $\calw$-vertical.
}

{\rm
Note that, since the definition of a generalized book of $I$-bundles
includes the condition that the underlying manifold be compact and
orientable, a generalized bookish pair is a compact, orientable
$3$-manifold pair.
}

\Claim\label{on accouna}
\textnormal{
If $(W,\redcale)$ is a generalized bookish pair, then $\redcale$ is
$\pi_1$-injective in $W$. To show this, we  
note that if
$\calw$ realizes  $(W,\redcale)$ then $\redcale$ is a disjoint union of
annuli in 
$W$ that are vertical in $\calw$. By an observation made in
\ref{vertical in book}, these are vertical annuli in 
the boundary of the Seifert fibered space
$\calb_\calw$, and are therefore $\pi_1$-injective in $\calb_\calw$;
and $\calb_\calw$ is in turn $\pi_1$-injective in $W$ by \ref{E
  or A is injective}. 
}
\EndClaim

\Claim\label{degenerates}
\textnormal{
A $3$-manifold pair will be termed {\it degenerate} if it is
 homeomorphic to $(D^2\times S^1,\alpha\times S^1)$, where $\alpha$ is
 an arc in $\partial D^2$. Note that a degenerate pair is in
 particular a generalized bookish pair, and is \realiz ed by a strict
 book of $I$-bundles.
}
\EndClaim

\Claim\label{excelsior}
\textnormal{
We define a {\it strict bookish
  pair} to be 
a generalized bookish pair $(\calw,\redcale)$ such that (1) the pair $(\calw,\redcale)$ is \realize d by some
strict book of $I$-bundles, and (2) 
no component of $(\calw,\redcale)$ is degenerate.
}
\EndClaim

\Claim\label{verticacklepuss}
\textnormal{
Let $(W,\redcale)$ be a  strict bookish pair. A submanifold of $W$
will be termed {\it
  $(W,\redcale)$-vertical } if it is $
\calw$-vertical for some strict book of $I$-bundles $\calw$ \realiz ing 
$(W,\redcale)$.
It follows from 
\ref{subsubmarine}
that if $U$ is a $(W,\redcale)$-vertical submanifold  of $W$, then $U$ is
strongly
  $\calw$-vertical for some strict book of $I$-bundles $\calw$ \realiz ing 
$(W,\redcale)$. (Indeed, any realization $\calw_0$ such that $U$ is
  $\calw_0$-vertical may be modified by a small isotopy to obtain a realization $\calw$ such that
$U$ is
strongly
  $\calw$-vertical.)
}

\EndClaim

Note that we have not defined ``bookish pair'' without the modifier
``generalized'' or ``strict,'' and we have not defined 
  $(W,\redcale)$-verticality  for a generalized
  bookish pair $(W,\redcale)$. 


\EndNumber


\Proposition\label{make 'em negative}
Every  generalized bookish pair is realized by a generalized
book of $I$-bundles each of whose pages has strictly negative
Euler characteristic. Furthermore, if a generalized book of
$I$-bundles 
$\calw$ \realiz es a given generalized bookish pair
$(W,\redcale)$, then $(W,\redcale)$ is \realiz ed by some generalized book of $I$-bundles
$\calw_0$ such that
(i)  each  page of $\calw_0$ has strictly negative
Euler characteristic, 
(ii) each page of $\calw_0$ is a page of $\calw$,
and (iii)
$\calf_{\calw_0}=\calf_{\calw}$.
\EndProposition

\Proof
This follows from the proof of
\cite[Lemma
2.3]{acs-singular}. If $\calw$ \realiz es the generalized bookish pair
$(W,\redcale)$, then in particular $\calw$ is a generalized book of
$I$-bundles with $|\calw|=W$. The  proof of \cite[Lemma
2.3]{acs-singular} then gives a generalized book of
$I$-bundles $\calw_0$ such that $|\calw_0|=W$, and such that every
page of $\calw_0$ has strictly negative Euler characteristic. The
proof also shows that $\calb_0\doteq \calb_{\calw_0}$ contains
$\calb\doteq \calb_\calw$, that each component of $\calb_0$ contains at
least one
component of $\calb$, 
that each page of $\calw_0$ is a page of $\calw$,
and that the Seifert fibration of $\calb_0$
extends the Seifert fibration of $\calb$. These conditions imply both
that  $\calf_{\calw_0}=\calf_{\calw}$, and that any vertical
submanifold of $\calb$ is also a vertical submanifold of $\calb_0$. In
particular, since $\redcale$, which is a disjoint union of annuli in
$\partial W$, is contained in $\calb$ and is vertical in $\calb$, it
follows that $\redcale$ is contained in $\calb_0$ and is vertical in
$\calb_0$. Hence, by an observation made in Subsection \ref{new-gbp},
we may assume, after possibly modifying  $\calw_0$ by a small isotopy,
that it is a realization of $(W,\redcale)$.
\EndProof

The  following lemma is a variant of \cite[Lemma 2.5]{acs-singular}.

\Lemma\label{even more like 2.5}
Suppose that $Q$ is a simple $3$-manifold and that $\calu$ is a
connected generalized book
of $I$-bundles such that $U\doteq|\calu|\subset Q$. Suppose that
every page of $\calu$ has strictly negative Euler characteristic and
is $\pi_1$-injective in
$Q$, and  that every fiber of every binding  of $\calu$ is also
$\pi_1$-injective in $Q$. Then there is a connected book of
$I$-bundles $\calu^+$ with $U^+\doteq |\calu^+|\subset Q$ 
such that
\begin{enumerate}[(a)]
\item $U^+\supset U$;
\item  every page of $\calu^+$ is a page of $\calu$;
\item $\partial U^+$ is a union of components of $\partial
  U$; 
\item every component of $\overline{U^+-U}$ is a solid torus;
and
\item for every binding $B$ of $\calu$ there is a binding $\hatB$ of
  $\calu^+$ such that $\hatB\supset B$, and the Seifert fibration of
  $\hatB$ extends the Seifert fibration of $B$.
\end{enumerate}
\EndLemma

\Proof
This follows from the  proof of  \cite[Lemma 2.5]{acs-singular},
specifically its third to sixth paragraphs. Before discussing the
details, we must point out that the hypothesis that the manifold $M$
is closed, which appears in the statement of \cite[Lemma
2.5]{acs-singular}, is the result of an editing error; it is not used
in the proof, and does not hold in the applications given in
\cite{acs-singular}. The corrected version, with the word ``closed''
removed, will be cited in the discussion below.

If we let the objects
denoted in the statement of the present lemma by $Q$, $\calu$ and $U$
play the respective roles of the objects denoted by $M$, $\calw_0$ and
$W$ in the proof of \cite[Lemma 2.5]{acs-singular}, then the four
cited paragraphs apply without change. (The hypothesis of the present
lemma imply all the hypotheses of the corrected version of \cite[Lemma 2.5]{acs-singular}
except the hypothesis that $\chi(W)<0$; and the latter hypothesis is
used only to establish, in the second paragraph of the proof of
\cite[Lemma 2.5]{acs-singular}, the fact, used in the third paragraph
of that proof, that the fibers of the bindings of $\calw$ are
homotopically non-trivial. The homotopic non-triviality of the fibers
of the bindings of $\calu$ is included in the hypothesis of the
present lemma.)

Hence, under the hypotheses of the present lemma, the conclusion of \cite[Lemma 2.5]{acs-singular}  holds with $Q$, $\calu$ and $U$
playing the respective roles of  $M$, $\calw_0$ and
$W$; that is, there is a connected book of $I$-bundles $\calu^+$ with
$U^+\doteq|\calu^+|\subset Q$ such that Conditions (1)---(7) of
\cite[Lemma 2.5]{acs-singular}  hold with $\calu^+$ and $U^+$ playing
the roles of $V$ and $\calv$. (Since $\calu^+$ is a book of
$I$-bundles in the slightly stronger sense defined in
\cite{acs-singular}, it is in particular a book of $I$-bundles in the
sense of the present paper.) Now Conditions (1), (4), (5) and (6) of
\cite[Lemma 2.5]{acs-singular}  give, respectively, Conditions (a),
(c), (d) and (b) of the conclusion of the present lemma. Condition
(e)  of the conclusion of the present lemma is not included in the
statement of \cite[Lemma 2.5]{acs-singular}; however, it is immediate
from the construction of the bindings of $\calv$ given in the third
paragraph of the proof of \cite[Lemma 2.5]{acs-singular}.
\EndProof

\Lemma\label{like 2.5}
Let $Q$ be a simple $3$-manifold, and let $(W,\redcale)$ be a generalized
bookish pair with $W\subset\inter Q$. Suppose that $(W,\redcale)$ has a
realization $\calw$ such that $\calp_\calw $ is $\pi_1$-injective in
$Q$, and such that every element of $\calf_\calw$ is also
$\pi_1$-injective in $Q$. Then
either ($\redalpha$) 
there are 
distinct
components $V_1$ and $V_2$ of $W$ such that 
the
component of $\overline{Q-V_1}$ containing $V_2$ is a solid torus
whose boundary is contained in $V_1$, or ($\redbeta$)
 there is a $3$-manifold pair $(W^+,\redcale^-)$ having the following properties:
\begin{enumerate}
\item $W\subset W^+\subset Q$;
\item each component of $(W^+,\redcale^-)$ is either a strict bookish pair or a degenerate pair;
\item $\partial W^+$ is a union of components of $\partial W$;
\item every component of $\overline{W^+-W}$ is a solid torus; and
\item $\redcale^-$ is the union of
all components of $\redcale$ contained in $\partial W^+$.
\end{enumerate}
\EndLemma

\Proof
According to Proposition \ref{make 'em negative}, 
$(W,\redcale)$ is \realiz ed by some generalized book of $I$-bundles
$\calw_0$ such that
each  page of $\calw_0$ has strictly negative
Euler characteristic, 
each page of $\calw_0$ is a page of $\calw$,
and 
$\calf_{\calw_0}=\calf_{\calw}$. In particular, 
$\calp_{\calw_0}$ is $\pi_1$-injective in $Q$, and
every element of $\calf_{\calw_0}$ is also
$\pi_1$-injective in $Q$.

If $\calu_1,\ldots,\calu_k$ denote the components (see \ref{books and
  such}) of $\calw_0$, it now follows that for $i=1,\ldots,k$, the
hypotheses of Lemma \ref{even more like 2.5} hold with $\calu_i$
playing the role of $\calu$. Hence there exist connected books of
$I$-bundles $\calu^+_1,\ldots,\calu^+_k$ with $|\calu^+_i|\subset Q$
such that, for $i=1,\ldots,k$, Conditions (a)--(e)
of  Lemma \ref{even more like 2.5} hold with $\calu^+_i$ and $\calu_i$
playing the respective roles of $\calu^+$ and $\calu$. We set
$U^+_i=|\calu^+_i|$ and $U_i=|\calu_i|$ for each $i$.

If the submanifolds $U^+_1,\ldots,U^+_k$ are not pairwise disjoint, then
there are distinct indices $i,j\in\{1,\ldots,k\}$ such that $U^+_j$ is
contained in a component $\Xi $ of $\overline{U^+_i-U_i}$.
According to Condition (d)
of  Lemma \ref{even more like 2.5}, $\Xi $ is a solid torus. Hence
Alternative ($\redalpha$) 
of the conclusion of the present lemma holds.

There remains the case in which  $U^+_1,\ldots,U^+_k$ are  pairwise
disjoint. In this case, $\calu^+_1,\ldots,\calu^+_k$ are the components of
a book of $I$-bundles $\calw^+$ with $W^+\doteq|\calw^+|\subset
Q$. Let $\redcale^-$ denote the union of
all components of $\redcale$ contained in $\partial W^+$; thus
$(W^+,\cale^-)$ is a $3$-manifold pair. We shall show that Alternative
($\redbeta$) of the conclusion holds in this case by verifying Properties
(1)--(5) of that alternative for the pair $(W^+,\cale^-)$.
Properties (1), (3) and (4) 
follow respectively from Conditions (a), (c) and (d) of the conclusion
of Lemma \ref{even more like 2.5} concerning the $\calu_i$. Property (5) is the definition of $\cale^-$. 

To verify Property
(2) 
for the pair $(W^+,\cale^-)$, first note that since each  page of the generalized
book of $I$-bundles $\calw_0$ has strictly negative
Euler characteristic, Condition  (b) of the conclusion
of Lemma \ref{even more like 2.5} guarantees that every page of the 
book of $I$-bundles $\calw^+$ has strictly negative
Euler characteristic; that is, $\calw^+$ is a strict book of
$I$-bundles. Next note that, since $\calw_0$ is a realization of the
generalized bookish pair $(W,\cale)$, the components of $\cale$ are
vertical annuli in the boundary of the generalized book of $I$-bundles
$\calw_0$. As was observed in 
\ref{vertical in book},
each component $E$ of $\cale$ is a vertical annulus in the boundary of
some binding $B_E$ of $\calw_0$. By Condition (e) of the conclusion
of Lemma \ref{even more like 2.5}, there is a binding $\hatB_E$ of
  $\calw^+$ such that $\hatB_E\supset B_E$, and the Seifert fibration of
  $\hatB_E$ extends the Seifert fibration of $B_E$. Hence $E$ is vertical
  in $\hatB_E$, and is therefore a vertical annulus in $\calw^+$. As
  this is the case for every component $E$ of $\cale$, the pair
$(W^+,\cale^-)$ is a generalized book of $I$-bundles realized by the
strict book of $I$-bundles $\calw^+$. It now follows from the
definition \ref{excelsior} that $(W^+,\cale^-)$ is either a degenerate
pair or a strict bookish pair. This is Property (2).
\EndProof


\DefinitionsRemarks\label{triad defs}
We define a {\it splitting} of a
 $3$-manifold $X$ to be
an ordered pair $(Z_1,Z_2)$, where $Z_1$ and $Z_2$ are 
three-dimensional submanifolds of $X$, closed as subsets of $X$, such
that 
$Z_1\cap Z_2$ is a
$2$-manifold, and $Z_1\cup Z_2=X$. Note that this definition implies
that $F\doteq Z_1\cap Z_2$
is a properly embedded surface in $X$, and that
$F=\Fr_XZ_i\subset\partial Z_i$ for $i=1,2$. In
particular $(Z_i,F)$ is a $3$-manifold pair for $i=1,2$.

Two splittings $(Z_1,Z_2)$ and $(Z_1',Z_2')$ of $X$ are said to be
{\it isotopic} if there is a self-homeomorphism $h$ of $X$, isotopic
to the identity, such that
for some---and hence for each---index $i\in\{1,2\}$ we have
$h(Z_i)=Z_i'$.

We define a {\it $3$-manifold triad}
to
be an ordered triple $(X,Z_1,Z_2)$ where $X$ is a $3$-manifold and $(Z_1,Z_2)$ is a splitting
of $X$. A $3$-manifold triad $(X,Z_1,Z_2)$ is said to be 
{\it compact,
connected,} or {\it orientable}  if $X$ is respectively compact,
connected, or orientable.
Note that connectedness of the triad $(X,Z_1,Z_2)$ does not imply that
$Z_1$ or $Z_2$ is connected.

We define an {\it acylindrical-bookish splitting} of a compact, connected, orientable
$3$-manifold $X$, 
or more briefly an  {\it \abs} of $X$,
to be a splitting $(C,W)$ of $X$ such that $(C,C\cap W)$ is an
acylindrical pair and $(W,C\cap W)$ is a strict bookish pair.
We define a {\it generalized \abs} of $X$
to be a splitting $(C,W)$ of $X$ such that $(C,C\cap W)$ is an
acylindrical pair and $(W,C\cap W)$ is a generalized bookish
pair.

Note that if $(C,W)$ is a generalized \abs\ of $X$ (and in particular
if $(C,W)$ is \anabs\ of $X$) then the components of 
$\redcale\doteq C\cap W$ are
annuli. According to the definition of an acylindrical pair (see
\ref{acylindrical}), $\redcale$ is $\pi_1$-injective in $C$; and
according to \ref{on accouna},
$\redcale$ is also $\pi_1$-injective in $W$. 
It follows that $C$, $W$ and
$\redcale$ are $\pi_1$-injective in $X$.

We define an {\it \abt}
(or  a  {\it generalized \abt}) to
be 
a 
compact, connected, orientable 
$3$-manifold triad $(X,C,W)$ such that 
$(C,W)$ is \anabs\ (or, respectively, a
generalized \abs)
of $X$. 
\EndDefinitionsRemarks

\Proposition\label{no nothing}
Let $X$ be a compact, connected,
orientable $3$-manifold. Suppose that $X$ is not a solid torus, and
that $X$ admits \anabs. Then no component of $\partial X$ is a sphere
or torus.
\EndProposition

\Proof
Fix an \abs\ $(C,W)$ of $X$, and set $\redcale=C\cap W$.
Let $Y$ be any  component of $\partial X$. We first consider the case where
$Y\cap\redcale=\emptyset$. In this case $Y$ is a component of either $\partial C$ or
$\partial W$. If $Y$ is a component of $\partial C$, then by
the discussion in \ref{acylindrical}, $Y$ is not a sphere or torus. If $Y$ is a
component of $\partial W$, 
then by \ref{derry down down},
$Y$ is not a sphere.
By \ref{derry ferry}, 
 if $Y$ is a  torus, 
then the component of
$C$ containing $Y$ is a solid torus; since $Y\cap\redcale=\emptyset$
and
$X$ is connected, it then follows that $X$ is a solid torus, a
contradiction to the hypothesis. Thus $Y$ cannot be a sphere or torus
in this case.

Now consider the case in which
$Y\cap\redcale\ne\emptyset$. In this case, the $2$-manifolds $Y\cap C$
and $Y\cap W$ are both non-empty
and have non-empty boundaries; 
they meet precisely in their common
boundary, and their union is $Y$. 
According to the discussion in \ref{acylindrical}, every component of
$(\partial
C)-\inter\redcale$ has strictly negative Euler characeristic. In
particular, every component of $Y\cap C$ has strictly negative Euler characeristic.  Since $Y\cap
C$ is non-empty, it now follows that
$\chi(Y\cap C)<0$. On the other hand, since $\redcale$ is $\pi_1$-injective
in $X$ by \ref{triad defs}, 
no component of  $Y\cap W$ is a disk; hence  $\chi(Y\cap
W)\le0$. 
We therefore have 
$\chi(Y)=\chi(Y\cap C)+\chi(Y\cap W)<0$,
which shows that $Y$ is not a sphere or a torus.
\EndProof

\Proposition\label{yes they are}
Suppose that $(C,W)$ is an \abs\ of a boundary-irreducible, compact, connected, orientable
$3$-manifold $X$. Then every component of $C\cap W$ is an essential
annulus (see \ref{essential annulus})
in $X$.
\EndProposition

\Proof
Set $\cale=C\cap W$.
According to \ref{E or A is injective}, each component of $\cale$ is a
$\pi_1$-injective annulus in $X$. It remains to show that no component
of $\cale$ is boundary-parallel. If $E$ is a boundary-parallel
component of $\cale$ then $E$ is the frontier of a submanifold $J$ of
$X$ such that $(J,E)$ is a degenerate  pair.
The submanifold $J$ must be a union of
 components of $C$ and components of $W$. According to \ref{derry down
   down},  any component of $W$ contained in $J$ must have non-positive
 Euler characteristic. Any component of $C$ contained in $J$ must have
 non-empty boundary, and since $(C,\cale)$ is an acylindrical pair,
 $\partial C$ has no sphere or torus boundary components by
 \ref{acylindrical}; hence every component of $C$ contained in $J$
 must have strictly negative Euler characteristic. But since every
 component of $\cale$ is an annulus, the sum of the
 Euler characteristics of the components of $C$ and of $W$ contained
 in $J$ is $\chi(J)=0$. It follows that $J$ contains no component of
 $C$, and by connectedness $J$ is a component of $W$. It now follows
 that $(W,\cale)$ has a  degenerate component, a contradiction to the
 definition of a strict bookish pair.
\EndProof

\Proposition\label{new characterization}
Let $\redX$ be a 
simple 
$3$-manifold. Then
up to isotopy $X$ admits a unique \abs\ 
$(C_0,W_0)$.

Furthermore, 
if
$\partial X\ne\emptyset$, and if
$U$ is a compact three-dimensional
submanifold of  $X$ such that

\begin{itemize}
\item  
$\Fr_XU$ is a properly embedded, $\pi_1$-injective, two-dimensional
submanifold of $\redX $, 
\item
each
component of $U$ has Euler characteristic $0$, and
\item no component of $\Fr_\redX U$ is 
a disk, a sphere or a
boundary-parallel annulus in $\redX $,
\end{itemize} 
then $U$ is isotopic in $X$ to a
$(W_0,C_0\cap W_0)$-vertical
 submanifold
 of 
$W_0$ which is disjoint from $W_0\cap C_0$.
\EndProposition

\Proof


To prove the first assertion, we first consider the case in which $\partial X=\emptyset$. In this
case, the pair $(X,\emptyset)$ is acylindrical. (Indeed,
Condition (1)
of Definition \ref{acylindrical}  holds in view of the simplicity of
$X$; Condition (2) holds trivially because the second term of the pair is
empty; Condition (3) holds because 
a non-empty  $I$-bundle has a non-empty boundary; 
and
Conditions (4) and (5) follow trivially from the fact that $X$ is closed.)
Hence $(X,\emptyset)$ is an \abs\ of $X$. 
To prove uniqueness in this
case, suppose that $(C,W)$ is an \abs\ of $X$. 
Since the closed
manifold $X$ cannot contain a properly embedded annulus, we have 
$C\cap W=\emptyset$; since $X$ is connected we have either
$C=\emptyset$ or $W=\emptyset$. The former alternative would give a
book of $I$-bundles $\calw$ with $|\calw|=X$;
since $X$ is connected and
$\partial X=\emptyset$, this contradicts \ref{confetti}.
Hence we have $(C,W)=(X,\emptyset)$, and uniqueness is established.

For the remainder of the proof of the first assertion, we assume that $\partial X\ne\emptyset$.
Then according to 
\ref{from MS, GS version}
the 
characteristic submanifold of $X$ is defined up to isotopy; 
it will be denoted by
$\Sigma_0$. 
According to Proposition \ref{who
  knows}, $\Sigma_0$ has the  properties 
(1) and (2) that are stated in
that proposition. By Property (1), 
we may write
$\Sigma_0=\caln_0\discup\calp_0$, where 
$\caln_0$ denotes the union of
all \bindinglike\ components of $\Sigma_0$, and   $\calp_0$
denotes the union of
all purely \pagelike\ components of $\Sigma_0$.

In view of an observation made in Subsection \ref{purely}, we may
equivalently define
$\calp_0$ to be the union of
all components of $\Sigma_0$ having strictly negative Euler
characteristic. 
(Either or both of the 
manifolds
$\caln_0$ and $\calp_0$
may be empty.)

Let us  set
$\scrL_0=\overline{X-\Sigma_0}$.

Since $\Sigma_0$ has the property
(2) that is stated in
Proposition \ref{who
  knows},
we
may write $\scrL_0$
as a disjoint union
$\Lacyl_0\discup\Ltriv_0$, where for each component $\maybeR $ of $\Lacyl_0$ the pair $(\maybeR ,\Fr_X\maybeR )$  is acylindrical,
and  each component  of $\Ltriv_0$ is a parallelism between two
components of 
$\Fr_X\Sigma_0$. (Condition (3) of the definition
(\ref{acylindrical}) of an acylindrical pair
guarantees that $\Lacyl_0$ and $\Ltriv_0$ are disjoint.)

Set $W_0=\Sigma_0\cup\Ltriv_0$. Then $(C_0,W_0)$ is a splitting of
$X$. 

We shall establish the existence assertion of the proposition by
showing that $(C_0,W_0)$ is \anabs.

Property 
(2)
of $\Sigma_0$ 
(from Proposition \ref{who
  knows})
also gives that for each component $\maybeR $ of $\Ltriv_0$, the two components of $\Fr_X\maybeR $ are contained in
distinct components of $\Sigma_0$, of which one is a component of
$\caln_0$ while the other is a component of $\calp_0$.
This implies that each component of $\caln_0\cup\Ltriv_0$
is \bindinglike, and in particular is a solid torus.

According to Property (1)
of 
$\Sigma_0$,
the components of $\Fr_X\Sigma_0$ (which contains both
$\Fr_X(\caln_0\cup\Ltriv_0)$ and 
$\calp_0\cap\Lacyl_0$) are essential annuli in
$X$. 
It now follows that if we fix a regular neighborhood $\calg_0$ of $\calp_0\cap\Lacyl_0$ in
$\Lacyl_0$, then each component of $\calb_0\doteq 
(\caln_0\cup\Ltriv_0)\discup\calg_0$ is a solid torus, and
each
component of $\Fr_X\calb_0$ is an essential annulus in $X$. 
Since the components of $\Fr_X\calb_0$ are in particular
$\pi_1$-injective annuli in $\calb_0$, we may fix a Seifert fibration
of $\calb_0$ in which $\Fr_X\calb_0$ is vertical. We observe, for
future reference, that the Seifert fibration of $\calb_0$ may be
chosen so that, in addition to the verticality of $\Fr_X\calb_0$, we
have that
\Claim\label{use it later}
$\caln_0$ is vertical in the Seifert fibration of $\calb_0$.
\EndClaim

Set 
$W_0^+=W_0\cup\calg_0$. Note that $W_0^+$ is isotopic to $W_0$ in
$X$. On the other hand, we have
$W_0^+=\calb_0\cup\calp_0$. By construction,
$\cala_0\doteq \calb_0\cap\calp_0$ is a union of components of
$\Fr_X\calb_0$, and $\cale_0\doteq \Fr_XW_0^+$ is equal to
$(\Fr_X\calb_0)-\cala_0$. Hence the components of $\cala_0$ and
$\cale_0$ are essential annuli in $X$, and in particular are
$\pi_1$-injective in $\calb_0$. Furthermore, $\calp_0$ is equipped
with an $I$-fibration whose vertical boundary is $\cala_0$, and each
component of $\calp_0$ has strictly negative Euler characteristic. It
now follows from the definitions that 
$\calw_0^+\doteq (W_0^+,\calb_0,\calp_0)$ 
is a strict book of $I$-bundles, with the Seifert fibraion of
$\calb_0$ and the $I$-fibration of $\calp_0$ that we have defined,
and that $\cale_0$ is 
$\calw_0^+$-vertical. 
The essentiality of the
components of $\cale_0$ includes the condition that none of them be
boundary-parallel; hence no component of $(W_0^+,\cale_0)$ is a degenerate
pair. This shows that $(W_0^+,\cale_0)=(W_0^+,\Fr_XW_0^+)$ is a strict
bookish pair. Since $W_0$ is isotopic to $W_0^+$, it follows that $(W_0,\Fr_XW_0)$ is a strict
bookish pair. But the pair $(C_0,\Fr_XC_0)$ is acylindrical by its very
definition. This shows that 
$(C_0,W_0)$ is \anabs, as claimed, and
establishes the existence assertion of the proposition.

Now consider an arbitrary \abs\
$(C,W)$ of $X$. We shall show that $W$ is isotopic to $W_0$ and hence
that $(C,W)$ is isotopic to $(C_0,W_0)$, thereby establishing the
uniqueness assertion of the lemma.

Let $\scrS$ denote the set of all submanifolds $\Sigma$ of $W$ 
satisfying the 
following conditions:

\begin{enumerate}[(i)]
\item $\Sigma\supset\Fr_XW$;
\item 
each component 
of $\Sigma$ is a \pagelike\ or \bindinglike\ submanifold of $X$ and
has essential frontier;
\item
each component
of $\overline{W-\Sigma}$ is a parallelism between two components of
$\Fr_W\Sigma$, 
which are not both
contained in \bindinglike\ components of $\Sigma$.
\end{enumerate}

Set $\cale=C\cap W$, and fix a strict book of $I$-bundles $\calw$ realizing the strict
bookish pair $(W,\cale)$. Fix a regular neighborhood $V$ of
$\cala_\calw$ in $W$ with $V\cap\cale=\emptyset$, and set
$\Sigma_1=\overline{W-V}$. We claim that $\Sigma_1\in\scrS$.

Condition (i) above obviously holds with $\Sigma=\Sigma_1$. 
To show that Condition (ii) holds with $\Sigma=\Sigma_1$, first note
that 
if $K$ is any component of $\Sigma_1$,
then $K$
admits either a binding or a page of $\calw$ as a regular
neighborhood. Hence 
$K$ is either a \bindinglike\ or a \pagelike\ submanifold of $X$.
On the other hand, 
since
every component of $\Fr_X\Sigma_1$ is parallel to a component of
$\cala_\calw$, 
Proposition \ref{yes they are} implies that 
each component of 
$\Fr_X\Sigma_1$
is an essential annulus in $X$; hence the \bindinglike\ or 
\pagelike\ submanifold $K$ of $X$ has essential frontier.
This
establishes Condition (ii).

To establish Condition (iii), 
note that $\overline{W-\Sigma_1}=V$.  
Consider an arbitrary
component   $L$ of $V$.
It is clear that $L$
is a parallelism between two components of
$\Fr_W\Sigma_1$. These  components of $\Fr_W\Sigma_1$ are contained in
components $K$ and $K'$ of $\Sigma_1$, where $K$ admits a
binding of $\calw$ as a regular neighborhood, and $K'$ admits a page
of $\calw$ as a regular neighborhood. Since $\calw$ is a strict book
of $I$-bundles, we have $\chi(K')<0$, and hence $K'$ is not a solid
torus. This establishes Condition 
(iii) 
with $\Sigma=\Sigma_1$, and
completes the proof that $\Sigma_1\in\scrS$.
 
In particular we have $\cals\ne\emptyset$.
Hence there is an element $\Sigma_2$ of $\cals$ such that
$\frakc(\overline{W-\Sigma_2})\le \frakc(\overline{W-\Sigma})$ for every $\Sigma\in\cals$.
We shall use Proposition  \ref{who knows}, with $\Sigma_2$ playing the role of $\Sigma$, to show that $\Sigma_2$ is isotopic to the
characteristic submanifold $\Sigma_0$ of $X$.

Since $\Sigma_2\in\scrS$, it follows
from Condition (ii)  of the definition of
$\scrS$ given above that $\Sigma_2$ has 
the property
(1) stated in Proposition \ref{who knows}.
To verify that $\Sigma_2$ has Property (2), we consider an arbitrary component $L$ of
$\overline{X-\Sigma_2}$. It follows from Condition (i) of the
definition of $\scrS$ that $L$ is a component of either
$\overline{X-W}$ or $\overline{W-\Sigma_2}$. If 
$L$ is a component of 
$\overline{X-W}$, then since $(C,W)$ is an \abs\ of $X$, the pair
$(L,\Fr_XL)$ is acylindrical. If $L$ is a component of
$\overline{W-\Sigma_2}$, then Condition 
(iii) 
of the definition of
$\scrS$ implies that $\Fr_XL$ has two components
 $\reddJ$ and $\reddJ'$, 
and
that $L$ is a parallelism between $\reddJ$ and $\reddJ'$. Let $K$ and $K'$ denote  the components (which a priori
need not be distinct) of $\Sigma_2$ containing $\reddJ$ and $\reddJ'$
respectively. Condition 
(iii) 
of the definition of
$\scrS$ 
also implies that $K$ and $K'$ 
are not both \bindinglike\ submanifolds of $X$.

Now suppose that $K$ and $K'$ are both \pagelike\ submanifolds of $X$.
 Then 
$K\cup K'$ can be given the
structure of an $I$-bundle over a $2$-manifold in such a way that
$(K\cup K')\cap\partial X$ is the horizontal boundary of $K\cup
K'$. (Note that in the union $K\cup K'$ the terms are a priori either disjoint
or equal.) Since $L$ is a parallelism between $\reddJ$ and $\reddJ'$, the $I$-fibration of $K\cup K'$ extends to an
$I$-fibration of $K\cup L\cup K'$ such that 
$(K\cup L\cup K')\cap\partial X$ is the horizontal boundary of $K\cup
L\cup K'$. 
Furthermore, the components of $\Fr_X(K\cup L\cup K')$ are
in particular components of $\Fr_X(K\cup K')$, and are therefore
essential annuli in $X$ by Condition (ii) of the definition of $\scrS$. It now follows that $K\cup L\cup K'$ is a
\pagelike\ submanifold of $X$, and hence that 
Condition (iii) of the definition of
$\scrS$ holds with $\Sigma_2\cup L$ playing the role of $\Sigma$. The
other conditions of the definition of $\scrS$ are immediate with
$\Sigma_2\cup L$ playing the role of $\Sigma$, and hence $\Sigma_2\cup
L\in\scrS$. Since $\frakc(\overline{W-(\Sigma_2\cup
  L)})=\frakc(\overline{W-\Sigma_2})-1$, this contradicts our choice of
$\Sigma_2$. 
Hence $K$ and $K'$ cannot both be \pagelike.

Since $K$ and $K'$  are not both \pagelike\and are not both
\bindinglike, they are distinct; furthermore, one of the submanifolds
$K$ and $K'$ must be purely \pagelike, while the other must be purely \bindinglike.

Property (2) 
for $\Sigma_2$ is
now established, and it follows from Proposition \ref{who knows} that
$\Sigma_2$ is isotopic to $\Sigma_0$. After an isotopy we may assume $\Sigma_2=\Sigma_0$.

We have observed that each
component  of
$\overline{X-\Sigma_2}$ is a component of either
$\overline{X-W}$ or $\overline{W-\Sigma_2}$. By 
Condition (iii) 
in the
definition of $\scrS$, 
each component  of 
$\overline{W-\Sigma_2}$ is a parallelism between two components of
$\Fr_W\Sigma_2$. On the other hand, for each component $L$ of 
$\overline{X-W}$, the pair $(L,\Fr_XL)$ is a component of the
acylindrical pair $(C,\Fr_XC)$; it then follows from Condition (3) of Definition
\ref{acylindrical}
that $L$ cannot be
a parallelism between two components of   $\Fr_W\Sigma_2$. 
In view of the definition of $\Ltriv_0$ and the
equality $\Sigma_2=\Sigma_0$, it now
follows that $\Ltriv_0=\overline{W-\Sigma_2}$, and therefore
$W_0=\Sigma_2\cup\Ltriv_0=W$. Hence the pairs $(C,W)$ and $(C_0,W_0)$
coincide, and the uniqueness assertion is established.

To prove the final assertion of the proposition, suppose that
$\partial X\ne\emptyset$, let $(C_0,W_0)$ be
defined as above, and suppose that $U$ is a compact three-dimensional
submanifold of $X$ such that the bulleted conditions in the statement
of the proposition hold.

Consider an arbitrary component $Q$ of $U$. If some component $S$ of $\partial Q$ is a $2$-sphere,
then  the simplicity of $X$ implies that $S$ is not a component
of $\partial X$; and since $\Fr_XQ$ is $\pi_1$-injective, some
component of  $\Fr_XQ$ must be a disk or a $2$-sphere, a
contradiction. Hence every component of $\partial Q$ has non-negative
Euler characteristic. But we have  $\chi(\partial
Q)=2\chi(Q)=0$, and hence every component of $\partial Q$ has 
Euler characteristic $0$, i.e. is a torus. Now since $\partial
X\ne\emptyset$, we have $\partial Q\ne\emptyset$. If we choose a
component $T$ of $\partial Q$, then since $T$ is a torus and $X$ is
simple, $T$ cannot be $\pi_1$-injective in $X$. Since $\Fr_XQ$ is
$\pi_1$-injective, $T$ is not $\pi_1$-injective in $Q$. But the
simplicity of $X$ and the $\pi_1$-injectivity of $\Fr_XQ$ also imply
that $Q$ is irreducible; since $\partial Q$ has a torus component
which is not $\pi_1$-injective in $Q$, it follows that $Q$ is a solid
torus. 

Each component of $\Fr_XQ$ is a compact $2$-manifold which is
$\pi_1$-injective in the solid torus $Q$, and is not a disk; hence
each such component is an annulus. Since $\Fr_XQ$ is $\pi_1$-injective
in $X$, and none of its components is a boundary-parallel annulus, all
its components are essential annuli. Furthermore, since $X$ is simple,
it cannot coincide with the solid torus $Q$, and hence
$\Fr_XQ\ne\emptyset$. It now follows that $Q$ is a pagelike
submanifold of $X$ with essential frontier. Since this is the case for
every component $Q$ of $U$, 
it follows from Lemma \ref{may 22} 
that $U$ is
isotopic in $X$ to a submanifold 
$U'$
of $\Sigma_0$. 


Recall that in the existence argument we fixed $I$-fibrations and
Seifert fibrations on $\calp_0$ and
$\calb_0$ respectively. Since $\caln_0$ is vertical in the Seifert
fibration of $\calb_0$ by \ref{use it later}, $\caln_0$ inherits an
$I$-fibration from $\calb_0$. In particular, each component of
$\caln_0$ or $\calp_0$ is now equipped with an $I$-fibration
or Seifert fibration, respectively. We claim:
\Claim\label{verticality claim}
For each component $K$ of $\Sigma_0=\caln_0\discup\calp_0$, the manifold  $U'\cap K$ is
isotopic in $K$ to a submanifold of $K$ which is vertical in the $I$-fibration
or Seifert fibration of $K$.
\EndClaim

To prove \ref{verticality claim},
we first consider the case
in which $K\subset\calp_0$. Set $V=\Fr_XK$. Since the components
of $\Fr_X(U'\cap K)\subset \Fr_XU'$ are essential annuli, these components of
$\pi_1$-injective and none of them is boundary-parallel in $K-V$. It
therefore follows from Corollary \ref{holiday-corollary}, applied with
$K$ and $U'\cap K$ playing the respective roles of $N$ and $N'$, that
$K\cap U'$ is
isotopic, by an ambient isotopy of $\calp_0$ which is constant on $V$
to a vertical submanifold $U''$ of the $I$-bundle $K$; thus
\ref{verticality claim} is established in this case.

Now consider the case in which $K\subset\caln_0$. In this case, again
using that the components of $\Fr_X(U'\cap K)=\Fr_K(U'\cap K)$ are $\pi_1$-injective
annuli, we apply  Proposition \ref{curse of the
  drinking class}, with $K$ defined as above and letting $\Fr_XK$ and
$\Fr_K(U'\cap K)$ play the respective roles of $F$ and $\cala$. This gives an
isotopy of $K$ which is constant on $\Fr_XK$ and carries $\Fr_X(U'\cap
K)$
onto a vertical submanifold of $K$. It follows that this isotopy
carries $U'\cap K$ onto a vertical submanifold of $K$, and
\ref{verticality claim} is established in this case as well. 

In view of \ref{verticality claim}, we may assume $U'$ to have been
chosen within its isotopy class in $\Sigma_0$ in such a way that
each of its components is
 vertical in the 
Seifert fibration or
$I$-fibration 
(respectively) of $\caln_0$ or
$\calp_0$. Since $\caln_0$ is vertical in $\calb_0$ and
inherits its Seifert fibration from $\calb_0$, we deduce that
each  component of $U'$ is either contained in $\calb_0$ and vertical
in its Seifert fibration, or contained in $\calp_0$ and vertical in its
$I$-fibration. According to the definitions this means that $U'$ is vertical in the
strict book of $I$-bundles $\calw_0^+$, and is therefore 
$(W_0^+,\Fr_XW_0^+)$-vertical.
Since $W_0^+$ is isotopic to $W_0$, the submanifold
$U'$---and hence the given submanifold $U$---is isotopic in $X$ to a
$(W_0,\Fr_XW_0)$-vertical submanifold $U''$ of $W_0$; after a further
small isotopy we may arrange that $U''$ is disjoint from
$\Fr_XW_0$. Since
$\Fr_XW_0=C_0\cap W_0$, this is the required conclusion.

\EndProof


\Notation\label{see X}
If $X$ is any simple manifold, the \abs\ given by Proposition
\ref{new characterization}, which is well-defined up to isotopy, will be
denoted by $(C_X,W_X)$. Thus $(X,C_X,W_X)$ is \anabt.
 We will also set $\redcale_X=C_X\cap W_X$.
\EndNotation

\Lemma\label{here again}
If $X$ is a simple $3$-manifold with $\partial X\ne\emptyset$, 
then for every prime $p$ we have 
$$\rk_\redp(C_\redX)\le \rk_\redp(\redX)+\frakc(C_X) -1\le \rk_\redp(\redX)+\chibar(X)-1.$$
\EndLemma

\Proof
We set $C=C_X$, $W=W_X$ and $\redcale=\redcale_X$, so that $C\cup W=X$ and
$C\cap W=\redcale$. We consider the following exact fragment of the
Mayer-Vietoris sequence,
where coefficients are taken in $\FF_p$ for a given prime $p$:
\Equation\label{pre-MV}
\xymatrix@C=1em{
H_1(\redcale) \ar[r]& H_1( C)\oplus H_1(W)\ar[r] &H_1(X)\ar[r]
&H_0(\redcale)\ar[r]& H_0(C)\oplus H_0(W)\ar[r] &H_0(X)\ar[r] &0.
}
\EndEquation

The exactness of (\ref{pre-MV}) implies:
\Equation\label{cool chicks}
\rk_\redp(\redcale)-(\rk_\redp(C)+\rk_\redp(W))+\rk_\redp(X)-\frakc(\redcale)+(\frakc(C)+\frakc(W))-\frakc(X)\ge0.
\EndEquation
We have $\frakc(X)=1$, and since the components of $\redcale$ are annuli,
we have $\rk_\redp(\redcale)=\frakc(\redcale)$. Thus (\ref{cool chicks})
simplifies to
\Equation\label{grand coulis}
\rk_\redp(C)\le\rk_\redp(X)+\frakc(C)-(\rk_\redp(W)-\frakc(W))-1.
\EndEquation

Now if  $L$ is any component of $W$, we have $\chibar(L)\ge0$ by
\ref{derry down down}. But since $X$ is connected and $\partial
X\ne\emptyset$ we have $\partial L\ne\emptyset$, and hence
$\chibar(L)=\rk_p(L)-(1+\dim H_2(L;\FF_p))\le \rk_p(L)-1$. We therefore
have 
$\rk_\redp(L)\ge1$ for each component $L$ of $W$, and hence 
$\rk_\redp(W)\ge\frakc(W)$.
Thus 
(\ref{grand coulis}) implies that
$\rk_\redp(C)\le\rk_\redp(X)+\frakc(C)-1$, which is the first inequality in
the conclusion of the lemma. 

The second inequality is
equivalent to the assertion $\chibar(X)\ge\frakc(C)$. To prove this,
first note that since $\partial X\ne\emptyset$ and $X$ is connected,
we have $\partial L\ne\emptyset$ for each component $L$ of 
$C$. Since, according to the discussion  in \ref{acylindrical}, no component of $\partial
C$ is a sphere or a torus,
it follows that $\chibar(L)\ge1$ for each component $L$ of $C$, so
that $\chibar(C)\ge\frakc(C)$. But
$\chibar(X)=\chibar(C)+\chibar(W)$ since the components of
$\redcale=C\cap W$ are annuli, and $\chibar(W)\ge0$ 
by \ref{derry down down}.
The required inequality $\chibar(X)\ge\frakc(C)$ now follows.
\EndProof

\Lemma\label{new where the cylinders are}
Let $(V,C,W)$ be a
generalized \abt, set $\redcale=C\cap W$, and let
$\calw$ be a  generalized book of
$I$-bundles which \realize s (\ref{new-gbp})
the generalized bookish pair
$(W,\redcale)$. Assume that $\calw$ has been chosen in such a way that each page of $\calw$ has strictly negative
Euler characteristic (cf. Proposition \ref{make 'em negative}).
Let $\psi$ be a 
(continuous) 
map from $V$ to 
a connected PL  
space $Q$.
Suppose that for each binding $B$ of $\calw$, 
there exist a space $\redscrSB$ of the homotopy type of $S^1$, and (continuous) maps
$\delta^B:B\to\redscrSB$ and $\epsilon^B:\redscrSB\to Q$, such that
$\epsilon^B\circ\delta^B=\psi|B$.
Then
for any prime $p$, we have
$$
\rk_\redp(\psi;V,Q)
\le\rk_\redp(C)+2\chibar(W)+3\chibar(C)-\frakc(C)+1.
$$
Furthermore, we have $\rk_\redp(\psi;V,Q)\le\max(1,\rk_\redp(C)+3\chibar(V))$,
and if either $C\ne\emptyset$ or $\chi(W)<0$, we have $\rk_\redp(\psi;V,Q)\le\rk_\redp(C)+3\chibar(V)$.
\EndLemma

\Proof
All homology groups in this proof will be understood to have
coefficients in $\FF_\redp$,
where $p$ is a given prime.

According to \ref{triad defs},
we have $C\cup W=V$; each component of  $\redcale\doteq C\cap W$ is an
annulus; $(C,\redcale)$ is an acylindrical pair; and $\redcale$ is
$\calw$-vertical.
We set $\calp=\calp_\calw$, $\calB=\calB_\calw$, and
$\redcalA=\redcalA_\calw$.
According to the definition of a 
generalized book of $I$-bundles (\ref{books and such}), 
$\calp$ has the structure of an $I$-bundle over a compact (but possibly
disconnected, empty, or non-orientable) $2$-manifold $S$; the
hypothesis implies that each 
component of $S$ has strictly negative Euler characteristic. 
We let $q:\calp\to S$
denote the fiber projection, so that $\redcalA=q^{-1}(\partial
 S)$. 
According to the definition of a
generalized book of $I$-bundles (\ref{books and such}), 
we have
$\calp\cap\calB=\cala$; and
by \ref{new-gbp} we have 
$\cale\subset\calb$
and 
$\redcale\cap\redcalA=\emptyset$.
Hence $\calp\cap C=\calp\cap(C\cap
W)=\calp\cap\cale=(\calp\cap\calb)\cap\cale=\cala\cap\cale=\emptyset$. To summarize:
\Claim\label{nuts and bolts}
The manifold $V$ is a threefold union $C\cup\calp\cup\calB$. We
have $\calB\cap C=\redcale$, $\calB\cap\calp=\redcalA$, and
$\calp\cap C=\emptyset$.
\EndClaim

Consider an arbitrary binding $B$ of $\calw$. Set $f^B=\psi|B:B\to Q$.
According to the hypothesis we have $f^B=\epsilon^B\circ\delta^B$.
We denote by $M^B$ the mapping cylinder of 
$\reddeltaB:B\to\redscrSB$. By definition, $M^B$ is obtained from the disjoint
union $(B\times [0,1])\discup\redscrSB$ by identifying $(x,1)$ with
$\reddeltaB(x)\in\redscrSB$ for every $x\in B$. The space $\redscrSB$ is
canonically identified with a deformation retract of $M^B$; in
particular:
\Claim\label{I said Aristotle}
$M^B$ has the homotopy type of $S^1$ for each binding $B$ of $\calw$.
\EndClaim

For each $B$, we define a map $\phi^B:(B\times [0,1])\discup\redscrSB\to Q$ by
$(x,t)\mapsto f^B(x)$ for $x\in B$ and $u\mapsto \epsilon^B(u)$ for
$u\in\redscrSB$. 
Since $\epsilon^B\circ\delta^B=f^B$, the map $\phi^B$
factors through a map $\hatf^B$ from the quotient space $M^B$ to $Q$. 

Now consider the disjoint union $V\discup\amalg_BM^B$, where $B$
ranges over the bindings of $\calw$, and 
define a space $\hatV$ to be the quotient of this disjoint union
formed by gluing $B\times\{0\}\subset M_B$ to $B\subset W\subset V$,
for each binding $B$, via the homeomorphism $(x,0)\mapsto x$. Then $V$
is canonically identified with a subspace of $\hatV$, as is each $M^B$. Since $V$ is connected by the definition of a
generalized \abt, 
 and since each $M^B$ is connected, the space $\hatV$ is connected. 
The spaces $V$ and $\amalg_BM^B$ are identified with subspaces of
$\hatV$, and there is a unique map $\hatpsi:\hatV\to Q$ which
restricts to $\psi$ on $V$ and to $\hatf^B$ on each $M^B$. If we denote by
$\iota$ the inclusion map $V\to\hatV$, then since
$\psi=\hatpsi\circ\iota$, it follows from
an observation made in \ref{redirects
  here}
that 
\Equation\label{yes i do}
\rk_\redp(\psi;V,Q)\le\rk_\redp(\hatV).
\EndEquation

Let $\hatcalB$ denote the union of the sets $M^B\subset\hatV$, where
$B$ ranges over the bindings of $\calw$.
From \ref{nuts and bolts} we deduce:
\Claim\label{hats and bats}
The space $\hatV$ is a threefold union $C\cup\calp\cup\hatcalB$. We
have $\hatcalB\cap C=\redcale$, $\hatcalB\cap\calp=\redcale'$, and
$\calp\cap C=\emptyset$.
\EndClaim

According to \ref{hats and bats}, we may regard $\hatV$ as the union
of the closed 
PL
sets 
$\hatcalB$ and $C\discup\calp$, and the intersection
of these sets is $\redcale^+\doteq \redcale\discup\redcale'$. From the
Mayer-Vietoris theorem, we have an exact sequence
\Equation\label{MV}
\xymatrix@C=1em{
H_1(\redcale^+) \ar[r]^-\alpha& H_1(\hatcalB)\oplus H_1(\calp \cup C)\ar[r] &H_1(\hatV)\ar[r]
&H_0(\redcale^+)\ar[r]& H_0(\hatcalB)\oplus H_0(\calp \cup C)\ar[r] &H_0(\hatV)\ar[r] &0,
}
\EndEquation
where $\alpha:
H_1(\redcale^+) \to H_1(\hatcalB)\oplus H_1(\calp \cup C)$
is not necessarily injective. If we denote the image of the
homomorphism $\alpha$ by $J$, we obtain an exact sequence
\Equation\label{after MV}
0\to J\to H_1(\hatcalB)\oplus H_1(\calp \cup C)\to H_1(\hatV)\to
H_0(\redcale^+)\to H_0(\hatcalB)\oplus H_0(\calp \cup C)\to H_0(\hatV)\to
0
\EndEquation
with $0$ at each end, so that the alternative sum of the dimensions of
the $\FF_\redp$-vector spaces in the sequence (\ref{after MV}) is $0$. Thus if we set
$m=\dim J$, we have
\Equation\label{becoming}
m-(\rk_\redp(\hatcalB)+\rk_\redp(\calp\cup C))
+\rk_\redp(\hatV)-\frakc(\redcale^+)+(\frakc(\hatcalB)+\frakc(\calp\cup C))-\frakc(\hatV)=0.
\EndEquation 
According to \ref{I said Aristotle}, $\hatcalB$ has the homotopy type
of a finite disjoint union of $1$-spheres, and thus
$\rk_\redp(\hatcalB)=\frakc(\hatcalB)$. We have observed that $V$ is
connected,
i.e. $\frakc(V)=1$. Hence (\ref{becoming}) simplifies to
\Equation\label{always}
\rk_\redp(\hatV)=\rk_\redp(\calp\cup C)+\frakc(\redcale^+)-\frakc(\calp\cup C)-m+1.
\EndEquation

We now proceed to estimate some of the terms on the right side of
(\ref{always}). It follows from \ref{hats and bats} that
$\redcale\cap\redcale'=\emptyset$ and hence
$\frakc(\redcale^+)=\frakc(\redcale)+\frakc(\redcale')$. 
The components of
$\redcale$ are annuli on the closed orientable surface $\partial C$. It
was pointed out in the discussion in \ref{acylindrical} 
that for  an acylindrical pair  $(C,\redcale)$, no component of $(\partial
C)-\inter\redcale$ is an annulus; hence the core curves of the components of
$\redcale$  are pairwise non-parallel. 
The
cardinality of any set of disjoint, non-parallel simple closed curves
on $\partial C$ is bounded above by $3\chibar(\partial
C)/2=3\chibar(C)$. Hence
\Equation\label{every bit helps}
\frakc(\redcale^+)\le3\chibar(C)+\frakc(\redcale').
\EndEquation

Of course, since $\calp$ and $C$ are disjoint by \ref{hats and bats},
we have
\Equation\label{diffidently}
\rk_\redp(\calp\cup C)=\rk_\redp(\calp)+\rk_\redp(C) \text{ and } \frakc(\calp\cup
C)=\frakc(\calp)+\frakc(C),
\EndEquation
and $H_1(\calp\cup C)$ is canonically identified with
$H_1(\calp)\oplus H_1(C)$.

Now recall that the homomorphism $\alpha$ in (\ref{MV}) is defined to
be $i_*+j^+_*$, where $i:\redcale^+\to\hatcalB$ and $j^+:\redcale^+\to\calp\cup
C$ are the inclusion homomorphisms. Hence $m=\dim J$ is bounded below
by the dimension of the image of $j^+$. But the image of $j^+$ is in
turn canonically identified with the direct sum of the images of $j_*$
and $j'_*$, where $j:\redcale\to C$ and $j':\redcale'\to\calp$ are the
inclusion homomorphisms. In particular we have
$m\ge\dim\image(j'_*)$. Since $\redcale'$ is the inverse image of
$\partial S$ under the $I$-bundle projection $q:\calp\to S$, the image
of $j'$ is isomorphic to the image of the inclusion homomorphism
$H_1(\partial S)\to H_1(S)$, whose dimension is $\frakc(\partial
S)-\frakc(S)=\frakc(\redcale')-\frakc(\calp)$. Hence:
\Equation\label{monkey pix}
m\ge\frakc(\redcale')-\frakc(\calp).
\EndEquation

Combining (\ref{always}) with (\ref{every bit helps}),
(\ref{diffidently}), and (\ref{monkey pix}), we obtain
\Equation\label{blatt}
\rk_\redp(\hatV)\le\rk_\redp(\calp)+\rk_\redp(C)+3\chibar(C)-\frakc(C)+1.
\EndEquation

Now note that since $S$ is a compact orientable surface, each of whose
components has strictly negative Euler characteristic, 
we have
$\rk_\redp(S)\le2\chibar(S)$.
Since $q:\calp\to S$ is a homotopy equivalence, it follows that
$\rk_\redp(\calp)\le2\chibar(\calp)$. Note also that since
$\calb\cup\calp=W$ and $\calb\cap\calp=\redcale'$, and since the
components of $\calp$ and $\redcale'$ have Euler characteristic $0$, we
have $\chibar(\calp)=\chibar(W)$. Combining these observations with
(\ref{blatt}), we obtain
\Equation\label{maybelline}
\rk_\redp(\hatV)\le\rk_\redp(C)+2\chibar(W)+3\chibar(C)-\frakc(C)+1.
\EndEquation
The first  assertion of the proposition
follows immediately from
(\ref{yes i do}) and (\ref{maybelline}).

To prove the remaining  assertions,
first consider the case in which
$C\ne\emptyset$. In this case we have $\frakc(C)\ge1$, so that the
 first assertion gives
$\rk_\redp(\psi)\le\rk_\redp(C)+2\chibar(W)+3\chibar(C)$. Since
$\chibar(W)\ge0$ 
by 
\ref{derry down down},
it follows that
$\rk_\redp(\psi)\le\rk_\redp(C)+3(\chibar(W)+\chibar(C))$. But since $V=C\cup W$, and
since each component of $\redcale=C\cap W$ has Euler characteristic $0$,
we have $\chibar(W)+\chibar(C)=\chibar(V)$, and hence $\rk_\redp(\psi)\le\rk_\redp(C)+3\chibar(V)$
in this case.

Next we consider the
case in which $C=\emptyset$ and $\chi(W)<0$. In this case the 
first assertion gives
$\rk_\redp(\psi)\le2\chibar(W)+1$; 
since $\chibar(W)\ge1$, it follows that
$\rk_\redp(\psi)\le3\chibar(W)$.

There remains the case in which $C=\emptyset$ and $\chi(W)\ge0$;
according to 
\ref{derry down down}
we must
  then
in fact have $\chi(W)=0$. 
The 
first assertion then gives
$\rk_\redp(\psi)\le1$, and the proof is complete.
\EndProof

\Proposition \label{use cylinders} 
Let $(V,C,W)$ be a generalized \abt.
Let $\psi$ be a (\redPL) map from $V$ to a simple $3$-manifold $Q$.
Suppose that 
the generalized bookish pair $(W,C\cap W)$ is \realize d by a
generalized book of $I$-bundles $\calw$ such that for every
$\jay\in\calf_\calw$ (see \ref{books and such}),
the map $\psi|\jay:\jay\to Q$
is $\pi_1$-injective.
Then
for every prime $p$ we have (in the notation of \ref{redirects here})
$$
\rk_\redp(\psi;V,Q)
\le\rk_\redp(C)+2\chibar(W)+3\chibar(C)-\frakc(C)+1.
$$
Furthermore, we have $\rk_\redp(\psi;V,Q)\le\max(1,\rk_\redp(C)+3\chibar(V))$,
and if either $C\ne\emptyset$ or $\chi(W)<0$, we have $\rk_\redp(\psi;V,Q)\le\rk_\redp(C)+3\chibar(V)$.
\EndProposition

\Proof
  We first prove the proposition under the additional hypothesis that
each
page of $\calw$ has strictly negative
Euler characteristic.

We wish to apply Lemma \ref{new where the cylinders are} with the
given $V$, $C$, $W$, $\redcale$, $\calw$, $Q$,  and $\psi$. For this purpose
we must
verify  that if $B$ is any binding of $\calw$, 
there exist a space 
$\scrS^B$ 
of the homotopy type of $S^1$, and maps
$\delta^B:B\to\scrS^B$ and $\epsilon^B:\scrS^B\to Q$, such that
$\epsilon^B\circ\delta^B=\psi|B$.

According to 
the definition of a generalized book of $I$-bundles (\ref{books and
  such}),
the
binding $B$ of $\calw$ is equipped with a Seifert fibration.
Choose a point
$x$ lying in a non-singular fiber of the given binding $B$, and set $y=f(x)$.  We denote by $J$ the fiber of
$B$ containing $x$. By definition we have $J\in\calf_\calw$. Hence by
hypothesis, 
the homomorphism
$(\psi|J)_\sharp: \pi_1(J,x)\to\pi_1(Q,y)$ is
injective. Hence if we set $f=\psi|B:B\to Q$, and  denote by $Z$ the image of  the inclusion
homomorphism $\pi_1(J,x)\to\pi_1(B,x)$, the subgroup
$f_\sharp(Z)$ of 
$\pi_1(Q,y)$ is infinite cyclic. 
We denote by $L$ the subgroup
$f_\sharp(\pi_1(B,x))$ of $\pi_1(Q,y)$. Since  $Z$ is a normal subgroup
of $\pi_1(B,x)$, the infinite cyclic subgroup $f_\sharp(Z)$ of
$L$ is also normal.
Since $Q$ is a simple $3$-manifold, it now follows from Lemma \ref{because
  g and c-r} that $L$ is infinite
cyclic. Hence
there exist a covering
space $\redPi :\tQ\to Q$ such that $\pi_1(\tQ)$ is infinite cyclic,
and a map $\tf:B\to\tQ$ such that 
$\redPi \circ\tf=f$.
Since the simple $3$-manifold $Q$
is  aspherical by \ref{keep it simple},
$\tQ$ has the homotopy type of $S^1$.
Thus the required condition holds if we set 
$\scrS^B=\tQ$,
$\delta^B=\tf$, and $\epsilon^B=\redPi $.

The conclusions of the proposition now follow immediately from those
of Lemma \ref{new where
  the cylinders are}.
This completes the proof under the assumption that
each
page of $\calw$ has strictly negative
Euler characteristic.

To prove the proposition in general, we invoke Proposition \ref{make
  'em negative} to give a generalized
book of $I$-bundles 
$\calw_0$ \realiz ing $(W,\redcale)$, such that
$\calf_{\calw_0}=\calf_{\calw}$ and each
page of $\calw_0$ has strictly negative
Euler characteristic.
Since $\calf_{\calw_0}=\calf_{\calw}$, in particular 
the map $\psi|\jay:\jay\to Q$
is $\pi_1$-injective
for every  $\jay\in\calf_{\calw_0}$. The conclusion therefore follows
upon applying the special case already proved, with $\calw_0$ in place
of $\calw$.
\EndProof

\Proposition\label{no cylinders}  
Let $(\redX,C,W)$ be \anabt. 
Then
for any prime $p$ we have
$$
\rk_\redp(X)
\le\rk_\redp(C)+2\chibar(W)+3\chibar(C)-\frakc(C)+1.
$$
If in addition we assume that $X$ is not a solid torus, then
$$\rk_\redp(\redX)
\le
\rk_\redp(C)+3\chibar(X).$$
\EndProposition

\Proof
We  apply Lemma \ref{new where the cylinders are}, taking $V=Q=X$, 
taking $\psi$ to be the identity map, 
using the $C$ and $W$   given by the hypothesis, and choosing
$\calw$  to be any strict book of $I$-bundles that realizes the strict bookish pair
$(W, C\cap W)$.
The definition of a strict book of $I$-bundles implies that each page
of $\calw$ has strictly negative Euler characteristic, as required for
Lemma  \ref{new where the cylinders are}.
The existence, for each binding $B$, of a space $\scrS^B$ and maps
$\delta^B$ and $\epsilon^B$ having the properties stated in the
hypothesis of Lemma \ref{new where the cylinders are} is trivial in
this situation: since $B$ is a solid torus and therefore has the
homotopy type of $S^1$, we need only set $\scrS^B=B$, and define
$\delta^B$ to be the identity map of $B$ and $\epsilon^B:B\to X$ to be the inclusion.

Since  in
the present context we have $\rk_\redp(\psi;V,Q)=\rk_\redp(\redX)$
for any given prime $p$,
it now
follows from  the first assertion of Lemma \ref{new where the cylinders are} that
$
\rk_\redp(X)
\le\rk_\redp(C)+2\chibar(W)+3\chibar(C)-\frakc(C)+1$,
which is the first 
assertion of the present  proposition.
If
$C\ne\emptyset$ or $\chi(W)<0$, then the second assertion of Lemma \ref{new where the cylinders are}
gives $\rk_\redp(X)\le\rk_\redp(C)+3\chibar(V)$, which is the conclusion of
the second assertion of the present proposition.
Now suppose that $C=\emptyset$ and that
$\chi(W)\ge0$. 
Since $\chi(W)\ge0$, the strict book of $I$-bundles $\calw$ has no
pages. Since $X=W$ is connected, $\calw$ consists of a single binding,
and hence $X$ is a solid torus. This contradicts the hypothesis of the
second assertion.
\EndProof

\Corollary\label{simple no cylinders}  
For any simple $3$-manifold  $\redX$
and any prime $p$, 
we have
$$
\rk_\redp(X)
\le\rk_\redp(C
_X)+2\chibar(W)+3\chibar(C_X)-\frakc(C_X)
+1,
$$
and
$$\rk_\redp(\redX)
\le
\rk_\redp(C_X)+3\chibar(X).$$
\EndCorollary

\Proof
Since $X$ is simple, it is not a solid torus. According to 
\ref{see X}, $(X,C_X,W_X)$ is \anabt.
The result therefore follows from
 Proposition \ref{no
  cylinders}.
\EndProof

Proposition \ref{no cylinders}  also yields the following corollary,
which was proved more directly (in the case $p=2$) as Lemma 2.11
of \cite{acs-singular}.
(The connectedness hypothesis was missing from the statement of Lemma 2.11 of
\cite{acs-singular}, but it is used in the proof given in \cite{acs-singular}, and
holds in the applications given there.)

\Corollary\label{no C}
If $\calw$ is a strict book of $I$-bundles such that $|\calw|$ is  connected, 
then for every prime $p$
we have
$$\rk_p(|\calw|)
\le
2\chibar(|\calw|)+1.$$
\EndCorollary

\Proof
Set $W=|\calw|$. Then $(W,\emptyset,W)$ is \anabt, and the result
follows from
the first assertion of
Proposition \ref{no cylinders}.
\EndProof

\section{Simplifying double curve intersections}\label{good reps}

\Lemma\label{better curve reduction}
Let $Q$ be a simple 3-manifold, 
let $\redPi  : \tQ\to Q$ be a two-sheeted covering, and let $\tau : \tQ\to \tQ$
denote the non-trivial deck transformation. Let $F\subset \inter\tQ$ be a
incompressible closed surface. 
Suppose that (a) $F$ meets $\tau(F)$ transversally, 
and (b) for
every surface $F_1\subset \tQ$,
  (ambiently)
isotopic to $F$, such that $F_1$ meets
$\tau(F_1)$ transversally, we have $\frakc(F \cap\tau(F ))\le
\frakc(F_1 \cap\tau(F_1))$. 
Then 
every component of $F \cap\tau(F )$ is a homotopically non-trivial
simple closed curve in $\tQ$.
\EndLemma

\Proof
Lemma 4.1 of \cite{acs-singular} asserts that if $Q$, $\tQ$ and $\tau$ satisfy the hypotheses of the
present lemma, then every
incompressible, surface in $\inter M$ is isotopic to a surface $F$ for
which the conclusion of the present lemma 
holds.

In the first paragraph of the proof of  \cite[Lemma
4.1]{acs-singular}, it is pointed out that any connected,
incompressible surface satisfying the hypothesis of that lemma is
isotopic to a surface $F$ satisfying Conditions (a) and (b) of the
present lemma, and it is explained that the rest of the proof is
devoted to showing that for such a surface $F$, every component of $F \cap\tau(F )$ is a homotopically non-trivial
simple closed curve in $\tQ$. Thus the proof of  \cite[Lemma
4.1]{acs-singular} establishes the present lemma.

(The hypothesis in \cite[Lemma
4.1]{acs-singular} that the surface has positive genus is an
inadvertent redundancy, and
should be ignored. A connected incompressible surface in a simple
manifold has strictly positive genus by definition, and for a disconnected
incompressible surface the genus is not defined, but each component of
such a surface has strictly positive genus.)

\EndProof

\Lemma\label{bestest curve reduction}
Let $Q$ be a simple, compact, orientable 3-manifold, 
let $\redPi : \tQ\to Q$ be a two-sheeted covering, and let $\tau : \tQ\to \tQ$
denote the non-trivial deck transformation. Let $R$ be a compact,
connected three-dimensional submanifold of $\inter\tQ$ such that
$\partial R$ is incompressible in $\tQ$. Then $R$ is isotopic to a
submanifold $\redRstar$ of $\tQ$ such that (1) 
$\partial\redRstar$ 
meets 
$\tau(\partial\redRstar)$ 
transversally,  (2) 
every component of $(\partial \redRstar) \cap\tau(\partial \redRstar )$ is a homotopically non-trivial
simple closed curve in $\tQ$, and (3) no component of
$\redRstar\cap\tau(\partial \redRstar )$ is a boundary-parallel annulus in $\redRstar$.
\EndLemma

\Proof
We fix a submanifold $F$ of $\tQ$, isotopic to $\partial R$, such that
$F$ meets $\tau(F)$ transversally. We may suppose
$F$ to be chosen in such a way that for every surface $\redFprime\subset \tQ$
which is
  (ambiently)
isotopic to $\partial R$, and such that $\redFprime$ meets
$\tau(\redFprime)$ transversally, we have $\frakc(F \cap\tau(F ))\le
\frakc(\redFprime \cap\tau(\redFprime))$. Then it follows from Lemma \ref{better
  curve reduction} that
every component of $F \cap\tau(F )$ is a homotopically non-trivial
simple closed curve in $\tQ$.

 Since $F$ is isotopic to $\partial R$,
there is a submanifold $R^*$ of $\tQ$, isotopic to $R$, with $\partial
R^*=F$. The properties of $F$ that we have stated imply that Conclusions
(1) and (2) of the present lemma hold with this choice of $R^*$. It
remains to verify Conclusion (3).

Assume that (3) is false. Then there is a submanifold $P$ of
$\redRstar$ such that (\reddestA) $\Fr_{R^*} P$ is an annulus and a component
of  $R^*\cap\tau(F)$, and (\reddestB) $P$ is
a parallellism (see \ref{spanish onion}) between $\Fr_{R^*}P$ and some
annulus in $F$. 
Among all
submanifolds of $R^*$ satisfying (\reddestA) and (\reddestB), we may suppose $P$ to
have been chosen to be minimal with respect to inclusion. 

We may write $\partial P=A\cup B$, where $A=P\cap\partial R^*$ and
$B=\Fr_{R^*}P$ are annuli; we have $A\cap B=\partial A=\partial
B$. We denote the components of $A\cap B$ by $\gamma_1$ and
$\gamma_2$. The annulus $A$ is contained in $F=\partial R^*$, while the
annulus $B$ 
is contained in $ \tau(F)$ and is properly embedded in $R^*$. We
claim that
\Equation\label{crable}
P\cap F=A\text{ and }P\cap \tau(F)=B.
\EndEquation

To prove (\ref{crable}), we first note that, since $F$ and $\tau(F)$
meet transversally, $A$ and $B$ are components of
$P\cap F$ and $P\cap \tau(F)$ respectively. We need to show that there are
no other components. Suppose to the contrary that $E$ is a component of
$P\cap F$ distinct from $A$ or a component of $P\cap \tau(F)$ distinct
from $B$. Then $E$ is a properly embedded $2$-manifold in $P$. If $E$
is closed, it is a component of $F$ or $\tau(F)$ and is therefore
incompressible in $\tQ$; since $E$ is contained in the solid torus
$P$, this is impossible. Now suppose that $\partial E\ne\emptyset$. If
$E$ is a component of $P\cap F$, its non-empty boundary is contained
in $B$. This is impossible because $B$ is properly embedded in $R^*$.

Now suppose that $E$ is a component of $P\cap\tau(F)$. Then $E$ is
properly embedded in $P$ and $\partial E\subset\inter A$. On the other
hand, since every
component of $F\cap\tau(F)$ is a homotopically non-trivial simple
closed curve, $E$ is $\pi_1$-injecctive in $\tau(F)$. In view of the
incompressibility of  $F$, the curve $E$ is in fact $\pi_1$-injective in $\tQ$
and therefore in $P$. Since $P$ is a solid torus, $\pi_1(E)$ is
cyclic. The non-triviality of the
components of $F\cap\tau(F)$ implies that $E$ is not a disk; hence it
is a $\pi_1$-injective annulus in $P$. Now Property (II) of the submanifold $P$ implies that
$P$ may be given the structure of an $I$-bundle over an annulus in
such a way that $A$ is a component of the horizontal boundary of
$P$. We  apply 
Proposition \ref{holiday-proposition}, 
letting the $I$-bundle $P$ play the role of $N$
in that 
proposition, 
and letting $E$ play the role of 
$\scrA$. 
We deduce that either (i) $E$ is isotopic, by an isotopy
which is constant on the vertical boundary $V$ of $P$, to a vertical
annulus, or (ii) $E$ is boundary-parallel in $P-V$. Alternative (i) is
impossible because $\partial E$ is contained in the single component
$A$ of the horizontal boundary. If (ii) holds, and if $P_1\subset P-V$ denotes a
parallelism between $E$ and an annulus in $\partial(P-V)$, then the properties (I)
and (II) persist when $P$ is replaced by $P_1$; this contradicts
the minimality of $P$. Thus (\ref{crable}) is proved.

(Note that the second equality in (\ref{crable}) implies that
$P$ is
contained either in $\tau(R^*)$ or in $\overline{\tQ-\tau(R^*)}$. 
This fact will not be used directly in the arguments.)

Next, we claim:
\Claim\label{after crable}
No component of $F\cap\tau(F)$ is invariant under $\tau$.
\EndClaim

To prove \ref{after crable}, suppose that $C$ is a $\tau$-invariant
component of $F\cap\tau(F)$.  Since $\tau$ has no fixed
point, it must preserve the orientation of $C$. Since, in addition,
the involution
$\tau$ interchanges the $2$-manifolds $F$ and $\tau(F)$, and since $C$
is a transverse component of intersection of these $2$-manifolds,
$\tau$ must reverse the orientation of
$\tQ$. This contradicts  the orientability of $Q$.

Combining (\ref{crable}) and \ref{after crable}, we will now show:

\Claim\label{el dorado}
Either $\tau(P)=P$ or $\tau(P)\cap P=\emptyset$. Furthermore, in the
former case, $P$ is a component of  $R^*\cap\tau(R^*)$.
\EndClaim

To prove \ref{el dorado}, we set $M_1=F\cap\tau(F)$,
$M_2=(F\cup\tau(F))-(F\cap\tau(F))$, and
$M_3=\tQ-(F\cup\tau(F))$. Since $F$ and $\tau(F)$ intersect
transversally, $M_i$ is an $i$-manifold for $i=1,2,3$. A submanifold
of $\tQ$ will be called a {\it stratum} if it is a component of $M_i$
for some $i$. The strata of $\tQ$ form a partition of $\tQ$ in the
sense that $\tQ$ is their set-theoretical disjoint union. Each of the
$M_i$ is invariant under $\tau$, and hence $\tau$ maps each stratum of
$\tQ$ onto a stratum.

It follows from (\ref{crable}) that $\inter P$, $\inter A$, $\inter
B$, $\gamma_1$ and $\gamma_2$ are strata of $\tQ$. The union of these
strata is $P$. 

Set $I=P\cap\tau(P)$. Then $I$ is a closed, $\tau$-invariant subset of
$\tQ$. Since $P$ is a union of strata, and $\tau$ maps strata to
strata, and the strata form a partition of $\tQ$, the set $I$ is a union
of strata. The strata contained in $I$ constitute a subset of those
contained in $P$, namely $\inter P$, $\inter A$, $\inter
B$, $\gamma_1$ and $\gamma_2$. 

If $I$ contains the stratum $\inter P$ (so that $I=P$), or if
$I=\emptyset$, then the first assertion of \ref{el dorado} is true. We
must rule out the other possibilities, namely that (i) $I$ contains at
least one of the strata $\inter A$, $\inter
B$, but does not contain $\inter P$, or (ii) $I$ contains at
least one of the strata $\gamma_1$, $\gamma_2$, but does not contain
$\inter A$ or $\inter
B$. 

Suppose that (i) holds. In the subcase where $\inter A\subset I$, the
stratum $\tau(\inter A)$ is also contained in $I$, and in particular in
$P$. But the stratum
is two-dimensional and is contained in $\tau(F)$; and of the five
strata contained in $P$, the only one that is two-dimensional and
contained in $\tau(F)$ is $\inter B$. Hence $\tau(\inter
A)=\inter B$. The same argument, with the roles of $A$ and $B$, and of
$F$ and $\tau(F)$, reversed, shows that in the subcase $\inter
B\subset I$ we  have $\tau(\inter
B)=\inter A$. Thus in either subcase, $\tau$ interchanges $A$ and
$B$, and hence leaves $\partial P$ invariant. Since the stratum
$\inter P$ is a component of $\tQ-\partial P$, the stratum
$\tau(\inter P)$ must also be a component of $\tQ-\partial P$. The
latter stratum must be distinct from $\inter P$ since
$I\not\supset\inter P$. Hence the stratum $\tau(\inter P)$ is equal to
$\tQ- P$. Since the strata of $\tQ$ form a partition, there must be
exactly six strat in this case: $\inter P$, $\tQ-P$, $\inter A$, $\inter
B$, $\gamma_1$ and $\gamma_2$. Of these, the only two-dimensional
stratum contained in $F$ is $\inter A$. But it follows from the
definition of strata that $F$ is the closure of the union of the
two-dimensiona strata which it contains. Hence $F=A$. This is a
contradiction, since $A$ is an annulus and $F$ is a closed surface.

Now suppose that (ii) holds. Then $I$ is equal to either $\gamma_1$,
$\gamma_2$ or $\gamma_1\cup\gamma_2$. 
Now each $\gamma_j$ is a component of 
$F\cap\tau(F)$; hence by \ref{after crable}, neither $\gamma_1$ nor
$\gamma_2$ can be $\tau$-invariant. Since $I$ is $\tau$-invariant,
we must have $I=\gamma_1\cup\gamma_2$, and $\tau$ must interchange 
$\gamma_1$ and
$\gamma_2$. 

In particular, $\gamma_1\cup\gamma_2=\partial A=\partial B$ is
$\tau$-invariant. Hence $\gamma_1\cup\gamma_2=(\partial
A)\cap\partial \tau(B)\subset A\cap\tau(B)\subset
P\cap\tau(P)=I=\gamma_1\cup\gamma_2$. This implies that 
$A\cap\tau(B)=\gamma_1\cap\gamma_2$. Thus the annuli $A$ and
$\tau(B)$, which are both contained in the $2$-manifold $F$, intersect
precisely in the boundary of $A$. Hence $F$ has a torus
component. This is impossible since $F$ is  incompressible in the
simple $3$-manifold $\tQ$. This completes the proof of the first
assertion of \ref{el dorado}.

To prove the second assertion of \ref{el dorado}, we note that in the
case where $P$ is $\tau$-invariant we have $P=\tau(P)\subset\tau(R^*)$ and hence $P\subset
R^*\cap\tau(R^*)$. Furthermore, $P$ is a  connected manifold, and
$R^*\cap\tau(R^*)$ is a manifold since $F$ and $\tau(F)$ meet transverally. Thus to prove that  $P$ is a component of
$R^*\cap\tau(R^*)$ in this case, it suffices to show that
$\Fr_{\tQ}P\subset\Fr_{\tQ}(R^*\cap\tau(R^*))$. We have $A\subset
F=\Fr_{\tQ}R^*$ and $A\subset P\subset 
R^*\cap\tau(R^*)$; hence $A\subset
\Fr_{\tQ}(R^*\cap\tau(R^*))$. Likewise, we have $B\subset
\tau(F)=\Fr_{\tQ}\tau(R^*)$ and $B\subset P\subset 
R^*\cap\tau(R^*)$; hence $B\subset
\Fr_{\tQ}(R^*\cap\tau(R^*))$. We therefore have $\Fr_{\tQ}P=A\cup
B\subset\Fr_{\tQ}(R^*\cap\tau(R^*))$, as required. This completes the
proof of \ref{el dorado}.

Now let $N$ be a regular neighborhood
of $P$ in $R^*$. We set $A'=N\cap F$ and
$B'=\Fr_{R^*}N$. 
We claim:

\Claim\label{regular guy}
The regular neighborhood $N$ of $P$ in $R^*$ may be chosen in such a
way that $\tau(F)\cap B'=\emptyset$ and $\tau(B')\cap B'=\emptyset$.
\EndClaim

To prove \ref{regular guy}, we show that each of the two stated
conditions holds when the regular neighborhood $N$
is small enough. To establish the first  assertion, we note that
$B$ is a component of $\tau(F)\cap R^*$ by Property (I) of the
submanifold $P$; hence when  the regular neighborhood $N$ of $P$ in
$R^*$ is small
enough,  $B'$ is disjoint from $\tau(F)\cap
R^*$ and hence from $\tau(F)$.

To establish the second  assertion, we observe that by 
\ref{el dorado} we have either (i) $\tau(P)\cap P=\emptyset$ or (ii)
$\tau(P)=P$. If (i) holds, it is obvious that when $N$ is sufficiently
small we have $\tau(N)\cap N=\emptyset$, and in particular
$\tau(B')\cap B'=\emptyset$. Now suppose that (ii) holds. Then by 
\ref{el dorado}, $P$ is a component of $R^*\cap\tau(R^*)$.
Hence when  the regular neighborhood $N$ of $P$ in
$R^*$ is small
enough,  $B'=\Fr_{R^*}N$ is disjoint from $R^*\cap\tau(R^*)$ and hence
from $\tau(R^*)$. In particular $B'$ is disjoint from $\tau(B')$.
This completes the proof of \ref{regular guy}.


We now suppose $N$ to have been chosen so that the conditions stated
in \ref{regular guy} hold. We set $G=\overline{F-A'}$. We also set $R'=\overline{R^*-N}$, so that $\partial R'=F'\doteq 
G\cup
B'$. Property (II) of the submanifold $P$ implies that $R'$ is
isotopic to $R^*$.

Consider the intersection
\Equation\label{bust 'em up}
F'\cap\tau(F')=(G\cap\tau(G))\cup
(G\cap\tau(B'))\cup(\tau(G)\cap B')\cup(B'\cap\tau(B')).
\EndEquation
The  fourth term on the right-hand side of (\ref{bust 'em
  up}) is empty by \ref{regular guy}.
The third term is contained in $\tau(F)\cap B'$ and is therefore empty by
\ref{regular guy}.
The second term is the image of the third term under
the involution $\tau$, and is therefore also empty. Hence
$F'\cap\tau(F')=G\cap\tau(G)$. 

We have
$G\cap\tau(G)\subset F\cap\tau(F)$. On the other hand, since
$A\subset\inter A'$, the simple closed curves $\gamma_1$ and $\gamma_2$, which are components of $F\cap\tau(F)$, are disjoint
from $G$ and hence from $G\cap\tau(G)$. Thus $G\cap\tau(G)\subset
(F\cap\tau(F))-(\gamma_1\cup\gamma_2)$, so that
$\frakc(F'\cap\tau(F'))=\frakc(G\cap\tau(G))\le
\frakc(F\cap\tau(F))-2$. Since $F'=\partial R'$ is isotopic to $\partial
R^*=F$, this contadicts our choice of  $F$.

\EndProof

\section{Pushing down submanifolds from a two-sheeted covering
  space}\label{pushing section}

\NotationRemarks\label{before gist}
Let $\redPsi$ be a compact PL space, and let
$\sigma$ be a fixed-point-free PL involution of a
compact PL
subset $U$   of  $\redPsi$. (A given component of $U$ may or may not be $\sigma$-invariant.)
We may define an
equivalence relation $\sim$ on the space $\redPsi$ by stipulating that $x\sim y$
if and only if either (a) $x=y$, or (b) $x\in U$ and $y=\sigma(x)$. We will denote
by $\Lambda _{\redPsi,\redUsigma}$ the quotient space $\redPsi/\sim$, and by
$\lambda_{\redPsi,\redUsigma}:\redPsi\to \Lambda _{\redPsi,\redUsigma}$ the
quotient map. 
Then $\Lambda _{\redPsi,\redUsigma}$ inherits a PL structure
from $\redPsi$, and the map $\lambda_{\redPsi,\redUsigma}$ is PL.

Note that
$\lambda_{\redPsi,\redUsigma}^{-1}(\lambda_{\redPsi,\redUsigma}(U))=U$, 
that $\lambda|U:U\to \lambda_{\redPsi,\redUsigma}(U)$ is a two-sheeted
covering map,
and that the map $\lambda_{\redPsi,\redUsigma}|(\redPsi-U)$ is a PL
homeomorphism onto its image.

Note also that:
\Claim\label{fever}
If 
$\Theta $ 
is a compact PL subset of $\redPsi$ with 
$U\subset \Theta $, then
the space $\Lambda_{\Theta ,U,\sigma}$ and the map
$\lambda_{\Theta ,U,\sigma}:\Theta \to\Lambda_{\Theta ,U,\sigma}$ are defined, and
there is a unique PL homeomorphism
 $\zeta:
\Lambda_{\Theta ,U,\sigma}\to \lambda_{\redPsi,U,\sigma}(\Theta )$
such that $\zeta \circ\lambda_{\Theta ,U,\sigma}=\lambda_{\redPsi,U,\sigma}|\Theta $.
\EndClaim
\EndNotationRemarks

\Number\label{how we use it}
The construction of \ref{before gist} will be used in the
context of covering space theory. Let $\redPi:\tQ\to Q$ be a two-sheeted covering of a compact PL
space $Q$, let $\tau:\tQ\to\tQ$ denote the non-trivial deck
transformation, and let $\redPsi $ be a compact PL subset of $\tQ$. Then
the compact PL set $\redPsi \cap\tau(\redPsi )$ is $\tau$-invariant. If
$U$ is a union of components of $\redPsi \cap\tau(\redPsi )$ which is
itself $\tau$-invariant, then $\sigma\doteq \tau|U$ is a PL involution of
$U$, and the PL space $\Lambda _{\redPsi ,\redUsigma}$ and the PL map
$\lambda_{\redPsi ,\redUsigma}$  are defined.
\EndNumber

\Number\label{second potato}
The discussion in \ref{how we use it} applies in particular to the
situation
in which $Q$ is a $3$-manifold, 
and
$\redPsi$ is a three-dimensional submanifold of
the interior of
 the two-sheeted
covering $\tQ$  such that
$\partial \redPsi $ and $\tau(\partial \redPsi )$ meet
transversally. 
Note that under these conditions, 
$\Psi\cap\tau(\Psi)$ is a
$3$-manifold, $\Fr_\Psi(\Psi\cap\tau(\Psi))$ and
$\Fr_{\tau(\Psi)}(\Psi\cap\tau(\Psi))$ are properly embedded
$2$-manifolds in $\Psi$ and $\tau(\Psi)$ respectively, and
$\Lambda
_{\redPsi ,\redUsigma}$ is a
$3$-manifold. If $\redPsi$ is compact or connected, then  $\Lambda
_{\redPsi ,\redUsigma}$
is, respectively,
compact or connected. If $Q$ is orientable, then $\Lambda
_{\redPsi ,\redUsigma}$ is
orientable.
Furthermore, if 
we are given  a $2$-manifold
$\calg\subset\partial \redPsi $ which is
disjoint from
$U$, then
 $\lambda_{\redPsi
  ,\redUsigma}|\calg$ is a homeomorphism onto its image, and 
$(\Lambda
_{\redPsi ,\redUsigma}, \lambda_{\redPsi ,\redUsigma}(\calg))$ is a
$3$-manifold pair (see \ref{real pairs}). 
\EndNumber

\Number\label{first potato}
Consider the special case of \ref{how we use it} in which we take $U$
to be the full intersection $\redPsi \cap\tau(\redPsi )$. In this case there is
a unique PL homeomorphism $\theta: \Lambda _{\redPsi ,\redUsigma}\to
\redPi(\redPsi )$ such that 
$\redPi|\redPsi $, regarded as a map of $\redPsi $ onto
$\redPi(\redPsi )
$, 
is equal to
$\theta\circ\lambda_{\redPsi ,\redUsigma}$.
\EndNumber


\Lemma\label{gisty numbers}
Let $\redPsi$ be a compact PL space, and let
$\sigma$ be a fixed-point-free PL involution of a
compact PL
subset $U$   of  $\redPsi$. Set $\Lambda=\Lambda
_{\redPsi,\redUsigma}$ and $\lambda=\lambda_{\redPsi,\redUsigma}:\redPsi\to\Lambda$. Then
$$
\chibar(\Lambda)=\chibar(\redPsi)-\frac12\chibar(U).
$$
Furthermore, if we denote the number of components of $U$ that are not
$\sigma$-invariant (which is clearly even) by $2\mni$, and if we denote the number of components of $U$ that are 
$\sigma$-invariant by $\minv$, then for any prime $p$, the codimension
in $H_1(\Lambda;\FF_p)$ of $\lambda_*(H_1(\redPsi;\FF_p))\subset
H_1(\Lambda;\FF_p)$ is at most $\minv+\mni$.
\EndLemma

\Proof
We set $V=\lambda(U)\subset\Lambda$ and $q=\lambda|U:U\to V$. Then $q$ is a
two-sheeted covering map.
We let $M$ denote
the mapping cylinder of 
$q$. By definition, $M$ is obtained from the disjoint
union $(U\times [0,1])\discup V$ by identifying $(x,1)$ with
$q(x)\in V$ for every $x\in U$. We regard $V$ as a subspace of $M$. The map of $(U\times [0,1])\discup
V$ to $V$ defined by $(x,t)
\mapsto 
q(x)$ for $x\in U,t\in[0,1]$ and $y
\mapsto
y$ for $y\in V$ induces  a canonical deformation
retraction $r:M\to V$.

Now 
define a space $\redLambdaplus$ 
to be the quotient of the disjoint union $\redPsi\discup M$
formed by gluing $U\times\{0\}\subset M$ to $U\subset\redPsi$
via the homeomorphism $(x,0)\mapsto x$. We identify $\redPsi$ and $M$
with subspaces of $\redLambdaplus$, so that $\redLambdaplus=\redPsi\cup M$, and define a map
$\phi:\redLambdaplus\to\Lambda$ by 
setting $\phi|\redPsi=\lambda$ and defining
$\phi|M$ to be the composition of $r$ with the inclusion
$V\hookrightarrow\Lambda$; note that $\phi$ is well-defined on $U=\redPsi\cap M$.
We now apply Proposition
\ref{from-brown-book}, taking
$X=\Lambda^+$,
$A=M$, $Q=\Lambda$, $B=V$, and $h=\phi$; the pair $(\Lambda^+,M)$ has
the homotopy extension property because it is a CW pair. This shows 
that $\phi$ is a
homotopy equivalence between $\redLambdaplus$ and $\Lambda$. 

Since $U$ is a two-sheeted covering of $V$, we have
$\chibar(M)=\chibar(V)=\chibar(U)/2$. 
The subspaces $\redPsi$ and $M$ are underlying sets of subcomplexes of
a CW complex with underlying space $\redLambdaplus$, and
$\redPsi\cap M$ is homeomorphic to $U$. Hence
$$\chibar(\Lambda)=\chibar(\redLambdaplus)=\chibar(\redPsi)+\chibar(M)-\chibar(U)=\chibar(\redPsi)+\frac12\chibar(U)-\chibar(U)=\chibar(\redPsi)-\frac12\chibar(U),$$
which gives the first assertion.

To prove the second assertion, we first note that
$\lambda:\redPsi\to\Lambda$ is the composition of the inclusion
$\redPsi\to\redLambdaplus$ with the homotopy equivalence $\phi$. Hence the codimension
in $H_1(\Lambda;\FF_p)$ of $\lambda_*(H_1(\redPsi;\FF_p))\subset
H_1(\Lambda;\FF_p)$ is equal to the codimension in
$H_1(\redLambdaplus;\FF_p)$ of the image of the inclusion homomorphism
$H_1(\redPsi;\FF_p)\to H_1(\redLambdaplus;\FF_p)$. It follows from the exact
homology sequence of the pair $(\redLambdaplus,\redPsi)$ that the latter
codimension is bounded above by $\dim H_1(\redLambdaplus,\redPsi;\FF_p)$, which by
the excision principle is equal to $\dim H_1(M,U)$. Hence it suffices
to show that $\dim H_1(M,U)\le \minv+\mni$.

We index the $\sigma$-invariant components of $U$ as
$U_1,\ldots,U_{\minv}$, and we index the remaining components of $U$ as
$U'_1,\ldots,U'_{\mni},U''_1,\ldots,U''_{\mni}$, where $\sigma$
interchanges $U_i'$ and $U_i''$ for each $i\in\{1,\ldots,{\mni}\}$. We
may then index the components of $M$ as
$M_1,\ldots,M_{\minv},M'_1,\ldots,M'_{\mni}$, in such a way that
$M_i\supset U_i$ for $i=1,\ldots,{\minv}$ and $M_i\supset U_i'\cup
U_i''$ for $i=1,\ldots,{\mni}$. We have
\Equation\label{flesh garden}
H_1(M,U;\FF_p)\cong
\bigg(\bigoplus_{i=1}^{\minv}H_1(M_i,U_i;\FF_p)\bigg)\oplus
\bigg(\bigoplus_{i=1}^{\mni}H_1(M'_i,U_i'\cup U_i'';\FF_p)\bigg).
\EndEquation
For $i=1,\ldots,\minv$, the space $M_i$ is the mapping cylinder of the
two-sheeted covering map $q|U_i:U_i\to V_i$, where
$V_i=U_i/\langle\sigma\rangle$. As $U_i$ is connected, we have $\dim
H_1(M_i,U_i;\FF_p)\le1$, with equality if and only if $p=2$. For
$i=1,\ldots,\mni$, the pair $(M'_i,U_i'\cup U_i'')$ is homeomorphic to
$(U_i'\times[0,1], U_i'\times\{0,1\})$, 
As $U_i'$ is connected, we have $\dim
H_1(M_i',U_i\cup U_i'';\FF_p)=\dim H_1 (U_i'\times[0,1],
U_i'\times\{0,1\})=1$. The required inequality $\dim H_1(M,U)\le
\minv+\mni$ now follows from (\ref{flesh garden}).

\EndProof

\Number\label{composition phenomenon}
Now let $\redPsi$ be a compact PL space, and let
$\sigma$ be a fixed-point-free PL involution of a
compact PL
subset $U$   of  $\redPsi$. 
Suppose that $U$ is written as a
disjoint union $U_\redone\discup U_\redtwo$, where each $U_i$ is a union of
components of $U$, and each $U_i$ is $\sigma$-invariant. 
For $i=\redone,\redtwo$
let $\sigma_i$ denote the PL involution $\sigma|U_i$ of $U_i$. 
Let us set $\Lambda=
\Lambda _{\redPsi,\redUsigma}$ and 
$\lambda=\lambda_{\redPsi,\redUsigma}:\redPsi\to\Lambda$. 
Let us also set 
$\redOmega_i=
\Lambda _{\redPsi,U_i,\sigma_i}$ and 
$\lambda_i=\lambda_{\redPsi,U_i,\sigma_i}:\redPsi\to\redOmega_i$. 
According
to \ref{before gist}, 
$\lambda_\redone|(\redPsi-U_1)$ 
is a homeomorphism onto
its image; in particular,
$\lambda_\redone$ restricts to a PL homeomorphism of $U_\redtwo$ onto a PL subset
$U_\redtwo\redstar$ of $\redOmega_1$. 
Conjugating the involution $\sigma_\redtwo$ of $U_\redtwo$
by the homeomorphism $\lambda_\redone|U_\redtwo:U_\redtwo\to U_\redtwo\redstar$, we obtain a PL
involution $\sigma_\redtwo\redstar$ of $U_\redtwo\redstar$. Set 
 $\Lambda^\dagger=
\Lambda _{\redOmega_1, U_\redtwo\redstar,\sigma_\redtwo\redstar}$ and 
$
\redlambdatwostar
=\lambda_{\redOmega_1,
  U_\redtwo\redstar,\sigma_\redtwo\redstar}:\redOmega_1
\to\Lambda^\dagger$.
Then there is a unique PL homeomorphism $\mu:\Lambda^\dagger\to\Lambda$ such
that $\mu\circ
\redlambdatwostar
\circ\lambda_\redone=\lambda$.
\EndNumber

\Lemma\label{gist}
Let $\tQ$ be a two-sheeted covering of a simple $3$-manifold $Q$, and
let  $\tau:\tQ\to\tQ$ denote the non-trivial deck transformation. Let
$W$ be a three-dimensional submanifold of $\inter \tQ$, and $\redcale$ a
two-dimensional submanifold of $\partial W$
such
 that $(W,\redcale)$
is a strict bookish pair. Suppose that $\partial W$ and $\tau(\partial W)$
meet transversally in $\tQ$, 
so that $W\cap\tau(W)$ is a $\tau$-invariant
three-dimensional
submanifold of $\tQ$. Let $U$ be a union of
components of 
$W\cap\tau(W)$, and suppose  that  $U$ is $(W,\redcale)$-vertical
and is itself $\tau$-invariant.
Suppose that  $U$ is disjoint from $\redcale$,
that each
component of $U$ has Euler characteristic $0$, and that $\Fr_WU$ is 
$\pi_1$-injective
in $\tQ$.
 Set $\sigma=\tau|U$. (Thus
by \ref{how we use it}, the PL space $\Lambda_{W,\redUsigma} $ and the
PL map
$\lambda_{W,\redUsigma}:W\to \Lambda_{W,\redUsigma} $ are
well-defined; and by \ref{second potato},
$\lambda_{W,\redUsigma}|\redcale$ is a homeomorphism onto its image, and 
$(\Lambda_{W,\redUsigma} ,\lambda_{W,\redUsigma}(\redcale))$
is a $3$-manifold pair.) 
Then $(\Lambda_{W,\redUsigma} ,\lambda_{W,\redUsigma}(\redcale))$ is a generalized bookish pair.
%
Furthermore, 
if $\calw$ is any strict book of $I$-bundles \realiz ing
$(W,\redcale)$, 
such that $U$ is 
strongly $\calw$-vertical
(such a realization exists by \ref{verticacklepuss}),
then $(\Lambda_{W,\redUsigma} ,\lambda_{W,\redUsigma}(\redcale))$ is \realiz ed by some
generalized book of $I$-bundles  $\calL$ with the following
properties:
\begin{itemize}
\item there is a 
$\calw$-vertical submanifold  of $W$ which is $\pi_1$-injective in
$W$,
is contained in $\calp_\calw$,
and 
is mapped  homeomorphically
onto $\calp_\calL$ by $\lambda_{W,\redUsigma}$; and
\item for every element
$F$ of $\calf_\calL$, there is a
simple closed curve $\redUpsilon\subset W$,
which is $\pi_1$-injective in $W$ and
is contained in either $\calb_\calw$ or in $\calp_\calw$, such
that 
$\lambda_{W,\redUsigma}(\redUpsilon)$ 
is a simple closed curve isotopic to $F$, and 
$\lambda_{W,\redUsigma}|\redUpsilon:\redUpsilon\to \lambda_{W,\redUsigma}(\redUpsilon)$ 
is a covering map
of degree at most $2$.
\end{itemize}

\EndLemma




\Proof
Let $\calw$ be any strict book of $I$-bundles \realiz ing
$(W,\redcale)$, such that $U$ is 
strongly
$\calw$-vertical.  
Set $\calb=\calb_\calw$, $\calp=\calp_\calw$, and
$\redcalA=\redcalA_\calw$. 
Set $\Lambda =\Lambda _{W,\redUsigma}$, and
$\lambda=\lambda_{W,\redUsigma}:W\to \Lambda $. 
We must show that $(\Lambda, \lambda(\redcale))$ is a generalized bookish
pair, that $\chi(\lambda^{-1}(L))=\chi(L)$, and that the two bulleted conditions of the statement hold.

We claim:
\Claim\label{new shnide}
Each component of $U$ is a solid torus, and each component of
$\Fr_W U$ is an annulus which is $\pi_1$-injective in 
$W$ and hence in $U$.
\EndClaim

To prove 
\ref{new shnide}, we first
note
that since 
$\Fr_W U$ is $\pi_1$-injective in $\tQ$, it is
in particular $\pi_1$-injective 
in $W$ and hence
in $U$. Now consider an arbitrary 
component $K$ of $U$. Since $K$ is $\calw$-vertical, and is disjoint
from $\redcalA$,
either (a) 
$K$ is contained in a page $P$ of $\calw$, which has the
structure of an $I$-bundle over a surface $S_P$, and
$K$ has the form $q_P^{-1}(T)$, where $q_P:P\to S_P$ denotes the
bundle projection, and $T$ is a subsurface of $S_P$ whose frontier
is a properly embedded one-dimensional submanifold of $S_P$; 
or
(b) $K$ is contained in a binding $B$ of $\calw$, which has the
structure of a Seifert fibered solid torus over a disk $\Delta_B$, and
$K$ has the form $q_B^{-1}(F)$, where $q_B:B\to \Delta_B$ denotes the Seifert
fibration, and $F$ is a subsurface of $\Delta_B$ whose frontier
is a properly embedded one-dimensional submanifold of $\Delta_B$.  If (a) holds, then
since $\chibar(K)=0$ by hypothesis, the surface $S_P$ must be
an annulus or a M\"obius band; this implies that $K$ is a solid torus,
and that the components of $\Fr_WK$ are annuli. Now suppose that (b)
holds. If some component of $\Fr_{\Delta_B} F$ is a simple closed
curve, then some component of $\Fr_W K$ is a  torus; this is
impossible, because $\Fr_WU$ is $\pi_1$-injective in $\tQ$ by
hypothesis, and $\pi_1(\tQ)$ has no rank-two free abelian subgroups by
the definition of simplicity. Hence every component of $\Fr_{\Delta_B} F$
is an arc, and hence $F$ is a disk. It now follows that, in this case
as well, $K$ is a solid torus,
and the components of $\Fr_WK$ are annuli. This completes the proof of
\ref{new shnide}.

Let $N$ be a regular neighborhood of $U$
relative to $W$. 
By hypothesis
 $U$ 
 is disjoint  from $\redcale$,
and  is strongly $\calw$-vertical  and hence disjoint  from
 $\redcalA$.
By \ref{new shnide}, the components of $\Fr_WU$ are annuli
and are $\pi_1$-injective in $W$. Hence
 we may suppose $N$ to have been chosen in
such a way that
\Claim\label{the customers say that}
$N$ is
$\calw$-vertical and is disjoint from $\redcalA\cup\redcale$,
and the
components of $\Phi\doteq \Fr_WN$ are annuli and are $\pi_1$-injective in
$W$ and hence in $N$.
\EndClaim

Now set $U_! =\lambda(U)\subset\Lambda$.
According to
\ref{before gist}, we have
$\lambda^{-1}(U_!)=U$, and  the map $\lambda|(W-U)$ is a 
homeomorphism onto its image. 
Since
$N$ is a regular
neighborhood of $U$ relative to $W$, it follows that:
\Claim\label{borneo}
The map $\lambda|\overline{W-N}$ is a
homeomorphism onto its image, and $N_!\doteq \lambda(N)$ is a regular
neighborhood of $U_!$ in $\Lambda$. Furthermore, we have
$\lambda^{-1}(N_!)=N$; in particular $\lambda$ maps
$\Phi$ onto $\Phi_!\doteq \Fr_\Lambda N_!$.
\EndClaim

We claim:

\Claim\label{cornuquackier}
Each component of $N_! $ is a solid torus, and each component of $\Phi_!$
is an annulus which is $\pi_1$-injective in $N_! $.
\EndClaim

To prove \ref{cornuquackier}, first note that $\lambda|U:U\to U_!$ is a
covering of the (possibly disconnected) orientable $3$-manifold
$U_!$ according to \ref{before gist}. 
Since each component of $U$ is  a solid torus by \ref{new
  shnide}, it follows from
Proposition \ref{it goes down good}
that each component of $U_!$ is a solid torus.
Since $N_!$ is a regular
neighborhood of $U_!$ in $\Lambda$ by \ref{borneo}, it follows
that each component of $N_! $ is a solid torus.

Now consider an arbitrary component $F$ of 
$\Phi_!$. It follows from
\ref{borneo} that $\lambda$ maps some component $\reddenedG$ of 
$\Phi$ homeomorphically onto $F$. 
Since $G$ is an annulus by \ref{the customers say that},
$F$ is also an
annulus. 
On the other hand,  \ref{the customers say that} also
implies that 
$\reddenedG$ is $\pi_1$-injective in
$N$. 
Since the covering map $\lambda|U:U\to U_!$ is $\pi_1$-injective,
the map $\lambda|N:N\to N_!$ between regular neighborhoods of $U$ and
$U_!$ is $\pi_1$-injective. It now follows that
$\lambda|\reddenedG:\reddenedG\to N_!$ is $\pi_1$-injective. Since
$\lambda$ maps  $\reddenedG$ homeomorphically onto $F$, we deduce that $F$ is
$\pi_1$-injective in $N_!$. This completes the proof of \ref{cornuquackier}.

According to the definition of a book of $I$-bundles, we have
$\calb\cup\calp=W$ and $\calb\cap\calp=\redcalA$. Since
$N\cap\redcalA=\emptyset$ by \ref{the customers say that}, each component
of $N$ is contained either in $\calp$ or in $\calb$.
We may therefore write $N=N^\calp\discup N^\calb$, where
$N^\calp\subset \calp$ and $N^\calb\subset\calb$. We set
$\Phi^\calb=\Phi\cap\calb=\Fr_\calb N^\calb$ and $\Phi^\calp=\Phi\cap\calp=\Fr_\calp N^\calp$, so that $\Phi=\Phi^\calb\discup\Phi^\calp$.
We also set
$\calb^-=\overline{\calb-N^\calb}$ and
  $\calp^-=\overline{\calp-N^\calp}$, so that
    $\calb^-\cup\calp^-=\overline{W-N}$ and
    $\calb^-\cap\calp^-=\redcalA$.

We have $W=\calb^-\cup \calp^-\cup N$, $\calb^-\cap\calp^-=\redcalA$,
$\calb^-\cap N=\Phi^\calb$, and $\calp^-\cap N=\Phi^\calp$. The
submanifolds $\redcalA$, $\redcale$, $\Phi^\calb$ and
$\Phi^\calp$ are pairwise disjoint.
    It now follows from \ref{borneo}
that:
\Claim\label{beyond} 
The submanifolds $\calb^-$, $\calp^-$, $\redcalA$, $\redcale$, $\Phi^\calb$, and
$\Phi^\calp$ of $W$ are mapped homeomorphically by $\lambda$ onto respective
submanifolds $\calb^-_!$, $\calp^-_!$, $\redcalA _!$, $\redcale_!$, $\Phi^\calb _!$, and
$\Phi^\calp _!$ of $\Lambda$. Each of
the submanifolds 
$\calb^-$, $\calp^-$, $\redcalA$, $\redcale$, $\Phi^\calb$ and
$\Phi^\calp$
is the full
preimage of its image under $\lambda$. We
have $\Lambda=\calb^-_!\cup \calp^-_!\cup N_!$, 
$\calb^-_!\cap\calp^-_!=\redcalA _!$,
$\calb^-_!\cap N_!=\Phi^\calb_!$, and 
$\calp^-_!\cap N_!=\Phi^\calp_!$. The
submanifolds $\redcalA_!$, $\redcale_!$, $\Phi^\calb_!$, and
$\Phi^\calp_!$ are pairwise disjoint.
\EndClaim

Since $(W,\cale)$ is in particular a generalized bookish pair,
it follows from 
the definitions of a generalized book of $I$-bundles (\ref{books and such}) and of
a generalized bookish pair (\ref{new-gbp})
that the components of $\redcalA$ and $\redcale$ are annuli. Hence by
\ref{beyond}, the components of $\redcalA_!$ and $\redcale_!$ are also
annuli. As the components of $\Phi_!=\Phi^\calb_!\discup
\Phi^\calp_!$ are annuli by \ref{cornuquackier}, we may state:
\Claim\label{new thing}
The components of $\redcalA_!$, $\redcale_!$, $\Phi^\calb_!$, and
$\Phi^\calp_!$ are annuli.
\EndClaim

According to the definition of a  book of $I$-bundles, $\calp$ has the
structure of an $I$-bundle over a $2$-manifold whose components are of
non-positive Euler characteristic, and  $\redcalA$ is the vertical
boundary of $\calp$.
It follows from  \ref{the
  customers say that} that $N^\calp$ is vertical in the $I$-bundle
$\calp$. 
Furthermore, since
$\Phi^\calp$ is a union of components of $\Phi$, it 
also follows from  \ref{the
  customers say that} that
$\Phi^\calp$
is $\pi_1$-injective 
in $W$ and hence
in $\calp$. The verticality of
$N^\calp$ in $\calp$, together with the $\pi_1$-injectivity of its
frontier, implies that
$\calp^-$ may be given the structure of an $I$-bundle over a surface
whose components are
of non-positive Euler characteristic, in such a way
that the 
vertical
boundary of $\calp^-$ is the disjoint union of $\redcalA$ and
$\Phi^\calp$. Combining this observation with \ref{beyond}, we deduce
that:
\Claim\label{poodles}
$\calp^-_!$ may be given the structure of an $I$-bundle over a surface whose components are
of non-positive Euler characteristic, in such a way
that the vertical boundary of $\calp^-_!$ is the disjoint union of $\redcalA_!$ and
$\Phi^\calp_!$.
\EndClaim

The definitions of a strict book of $I$-bundles,  and of a strict bookish
pair, also give that each component
of $\calb$ is a solid torus, and that $\calb$ is equipped with a
Seifert fibration
in which $\redcalA\cup\redcale$ is  vertical.
It follows from  \ref{the
  customers say that} that $N^\calb$ is vertical in the Seifert
fibration of
$\calb$. 
Furthermore, since
$\Phi^\calb$ is a union of components of $\Phi$, it 
also follows from  \ref{the
  customers say that} that
the components of 
$\Phi^\calb=\Fr_\calb \calb^-$ are annuli.
These observations imply
that each component of
$\calb^-$ is a solid torus, 
and that the components of $\redcalA$, $\redcale$ and
$\Phi^\calb$ are  contained in $\partial \calb$ and are
$\pi_1$-injective in $\calb$. Combining this observation with \ref{beyond}, we deduce:
\Claim\label{toodle-oo}
Each component of
$\calb^-_!$ is a solid torus, and the components of $\redcalA_!$, $\redcale_!$ and
$\Phi^\calb_!$ are 
contained in $\partial \calb_!$ and are
$\pi_1$-injective in $\calb_!$.
\EndClaim

According to \ref{cornuquackier} and
\ref{toodle-oo}, the components of both $N_! $ and $\calb^-_!$ are
solid tori.
By \ref{beyond} we have $\Phi^\calb_!=N_!\cap \calb^-_!$, and by
\ref{new thing} and \ref{toodle-oo}, the components of 
$\Phi^\calb_!$ are annuli, contained in $\partial \calb_!$ and in
$\partial N_!$, and
$\pi_1$-injective in $\calb_!$ and in
$N_!$. 
It follows that $N_!\cup\calb^-_!$
admits a Seifert fibration in which the components of
$N_! $, $\calb^-_!$  and $\Phi^\calb_!$ are vertical.

According to \ref{beyond}, \ref{new thing} and \ref{toodle-oo},
$\redcale_!$ and $\redcalA_!$ are disjoint from each other and from
$\Phi^\calb_!$; and the components of $\redcale_!$ and $\redcalA_!$ are annuli in 
$\partial\calb_!^-$
which are
$\pi_1$-injective  in $\calb_!^-$.
It follows that $\redcale_!$ and $\redcalA_!$ are contained in the boundary
of $N_!\cup\calb^-_!$, and that we
may choose the Seifert fibration of 
$N_!\cup\calb^-_!$ in such a way that $\redcale_!$ and $\redcalA_!$ are vertical.

We have $\Phi^\calp_!\subset\Phi=\Fr_\Lambda N_!\subset\partial N_!$,
and according to \ref{cornuquackier} and \ref{beyond}, 
the components
of $\Phi^\calp_!$ are annuli, disjoint from $\Phi^\calb_!$, and
$\pi_1$-injective in $N_!$. It follows that 
$\Phi^\calp_!$ is contained in the boundary
of $N_!\cup\calb^-_!$, and that we
may choose the Seifert fibration of 
$N_!\cup\calb^-_!$ in such a way that 
$\Phi^\calp_!$ is vertical.

By \ref{beyond}, the manifold  $\Lambda$ is the union  of the Seifert
fibered space $N_!\cup\calb^-_!$ and the $I$-bundle
$\calp^-_!$.
It also
follows from \ref{beyond} that  the intersection of 
$N_!\cup\calb^-_!$ with  $\calp^-_!$
is $\Phi^\calp_!\cup\redcalA_!$.
Since  $\Phi^\calp_!$ and $\redcalA_!$ are
vertical in the Seifert fibration of $N_!\cup\calb^-_!$,
and their
union is the vertical boundary of
$\calp^-_!$ by \ref {poodles}, we may define  a generalized book of
$I$-bundles $\calL$ with $|\calL|=\Lambda$ by setting
$\calb_\calL=N_!\cup\calb^-_!$ and $\calp_\calL=\calp^-_!$.
Since $\redcale_!$ is vertical in the Seifert fibration of 
$N_!\cup\calb^-_!$, it is vertical in the book of $I$-bundles
$\calL$. Furthermore, the components of $\redcale_!$ are annuli by
\ref{new thing}; and since $\redcale_!$ is contained in the boundary of
$N_!\cup\calb^-_!$ and is disjoint from $\redcalA_!$, it is contained in
the boundary of $\Lambda$. This shows that $(\Lambda,\redcale_!)$ is a
generalized bookish pair, and is realized by $\calL$.

It remains to show that  the realization
$\calL$ of
$(\Lambda,\redcale_!)$ satsifies the two bulleted conditions of the statement.

To verify the first bulleted condition, we observe that, according to
\ref{beyond}, $\lambda$ maps $\calp^-$ homeomorphically onto
$\calp^-_!=\calp_\calL$. The manifold
$\calp^-=\overline{\calp-N^\calp}\subset\calp\subset W$ is
$\calw$-vertical because $N^\calp$ is vertical in $\calp$ by \ref{the
  customers say that}. Furthermore,
$\Fr_{\calp}\calp^-=\Phi^\calp$ is a union of components of
$\Phi$, and is therefore $\pi_1$-injective in $W$---and in
particular in $\calp$---by \ref{the
  customers say that}.
Hence  $\calp^-$ is $\pi_1$-injective in
$\calp$. As
$\calp$ is in turn $\pi_1$-injective in $W$ 
by
\ref{E or A is injective},
we may now deduce that $\calp^-$ is $\pi_1$-injective in
$W$. Thus the first bulleted condition is established.

To establish the second bulleted condition, let $F\in\calf_\calL$ be
given. By the definition  of $\calf_\calL$, there is a simple closed
curve $\Upsilon_!$, isotopic to $F$ in $\Lambda$, such that
$\Upsilon_!$ is a non-singular fiber
of a binding of $\calL$. In particular we have $\Upsilon_!\subset\calb_\calL$. Since
$\calb_\calL$ is the union of the vertical submanifolds $N_!$ and
$\calb^-_!$, and since $\Upsilon_!$ is a fiber, 
$\Upsilon_!$ must be contained in either
$N_!$ or $\inter\calb^-_!$; furthermore, as the vertical submanifolds
$N_!$ and 
$\inter\calb^-_!$ are Seifert-fibered solid tori,
the fiber $\Upsilon_!$ is
$\pi_1$-injective in $N_!$ or $\inter\calb^-_!$. In the case $\Upsilon_!\subset N_!$, since $N_!$ is
a regular neigbhorhood of $U_!$ by \ref{borneo}, we may choose $\Upsilon_!$
within its isotopy class so that $\Upsilon_!\subset \inter U_!$. Thus
in any case there is a component $H_!$ of
either $U_!$ or $\calb^-_!$ such that $\Upsilon_!\subset \inter H_!$,
and $\Upsilon_!$ is $\pi_1$-injective in $H_!$.

If $H_!$ is a component of $U_!$, it follows from \ref{before gist}
that
 $\lambda|\lambda^{-1}(H_!): \lambda^{-1}(H_!) \to H_!$ is a
two-sheeted covering map (with possibly disconnected domain).  If $H_!$ is a
component of $\calb^-_!$, it follows from \ref{beyond} that
$\lambda^{-1}(H_!)$ is a component of $\calb^-$, and that
$\lambda|\lambda^{-1}(H_!): \lambda^{-1}(H_!) \to H_!$ is a
homeomorphism. 
Hence, in either case, if we select a component  $H$  of $\lambda^{-1}(H_!)$, the map
$\lambda|H:H\to H_!$ is a covering map of degree at most $2$. It
follows that if we choose a  component $\Upsilon$ of
$\lambda^{-1}(\Upsilon_!)$ contained in $H$, the map $\lambda|\Upsilon:\Upsilon\to\Upsilon_!$ is a
covering map of degree at most $2$. In particular we have $\lambda(\Upsilon)=\Upsilon_!$.
Since $\Upsilon_!$ is $\pi_1$-injective in $H_!$, the 
curve $\Upsilon$ is $\pi_1$-injective in $H$.

We claim that $H$ is $\pi_1$-injective in $W$. If  $H_!$ is a
component of $U_!$, then $H$ is a component of $U$, since
$\lambda^{-1}(U_!)=U$
by \ref{before gist};
the asserted $\pi_1$-injectivity in this case then follows from \ref{new shnide}. If
$H_!$ is a component of $\calb^-_!$, then 
$H$ is a component of $\calb^-$, and
to prove $\pi_1$-injectivity
of $H$ in $W$ it suffices to prove that $\Fr_W\calb^-$ is $\pi_1$-injective in $W$. But
$\Fr_W\calb^-$ is the disjoint union of $\redcalA_\calw$ with $\Phi^\calb$;
and  $\redcalA_\calw$ is $\pi_1$-injective 
by
\ref{E or A is injective},
while  $\Phi^\calb$ is a disjoint union of components of $\Phi$, and
is therefore $\pi_1$-injective in $W$ by \ref{the customers say that}.
This completes the proof that $H$ is
$\pi_1$-injective in $W$. 

Since we have seen that $\Upsilon$ is
$\pi_1$-injective in $H$, it now follows  that $\Upsilon$ is
$\pi_1$-injective in $W$.

Finally, we must show that $\Upsilon$ is contained in either $\calb$
or $\calp$. Since $\Upsilon\subset\lambda^{-1}(\Upsilon_1)$, and $\Upsilon_!$ is contained in either
$N_!$ or $\inter\calb^-_!$, we have either
$\Upsilon\subset\lambda^{-1}(N_!)$ or
$\Upsilon\subset\lambda^{-1}(\inter\calb^-_!)$. If
$\Upsilon\subset\lambda^{-1}(N_!)$, then since $\lambda^{-1}(N_!)=N$
by \ref{borneo}, and $N$ is disjoint from $\cale$ by \ref{the
  customers say that}, we have $\Upsilon\cap\cale=\emptyset$, and
hence $\Upsilon$ is contained in either $\calb$
or $\calp$. If $\Upsilon\subset\lambda^{-1}(\inter\calb^-_!)$, then
it follows from \ref{beyond} that
$\Upsilon\subset\inter\calb^-\subset\calb$. Thus the second bulleted condition of the
statement is established.
\EndProof

\Lemma\label{key lemma} 
Let $Q$ be a simple $3$-manifold, and let  $\redPi:\tQ\to Q$ be a
connected
two-sheeted covering of $Q$ (so that $\tQ$ is simple by 
\ref{keep it simple}).
Let $R\subset\inter\tQ$ be a
connected
three-dimensional 
submanifold such that
$\partial R$ is  non-empty and is incompressible in $\tQ$ (so that $R$
is simple by 
Remark
\ref{more fruit}, and hence $C_R$,  $W_R$ and $\cale_R$ are defined by \ref{see X}).
Let $K$ be a compact, connected PL space
such that $\pi_1(K)$ is freely indecomposable, and let $f:K\to \tQ$ be
a $\pi_1$-injective \redPL map with  $f(K)\subset R$. Then $\redPi\circ f:K\to Q$
is homotopic to a \redPL map whose image is contained in a
compact, connected three-dimensional
 submanifold $T$ of
$\inter Q$, 
where $\partial T$ is 
incompressible (so that $T$ is simple by 
Remark
\ref{more fruit}, 
and hence $C_T$, $W_T$ and $\cale_T$ are defined by \ref{see X}),
and one of the following alternatives holds:
\Alternatives
\item
 $\chibar(T)=\chibar(R)$,
and the pairs $(C_T,\redcale_T)$ and $(C_R,\redcale_R)$ are homeomorphic, 
or 
\item
$\chibar(T)<\chibar(R)$, and 
for every prime $p$ we have
  $$\rk_\redp(T)\le \rk_\redp(C_R)+4\chibar(R)+\chibar(T).$$
\EndAlternatives
\EndLemma

\begin{unnumberedremark}\label{a key goes with it}
\textnormal{
The conclusion of Lemma \ref{key lemma} does not assert that
$\partial T\ne\emptyset$, or equivalently that $T\ne Q$.
However, when Alternative (i) of the conclusion holds,
we must have $ T\ne Q$; this follows, for example, from
the equality $\chibar(T)=\chibar(R)$, since the simple manifold $R$ with
non-empty boundary satisfies $\chibar(R)>0$ by
\ref{keep it simple},
whereas a closed $3$-manifold has Euler characteristic $0$.
When Alternative (ii) holds, there appears to be no general reason why
$\partial T$ should be non-empty.
}
\end{unnumberedremark}

\Proof[Proof of Lemma \ref{key lemma}]
We set $W=W_R$, $C=C_R$, and $\redcale=\redcale_R$. We recall from 
Proposition
\ref{new characterization} and Subsection \ref{see X} that $(C,W)$ is \anabs\ of $R$ which, a
priori, is defined only up to isotopy. Recall also that, according to
the definition of \anabs, $(W,\redcale)$ is a strict bookish pair,
$(C,\redcale)$ is an acylindrical pair, and $\redcale=C\cap W$.

It follows from Lemma \ref{bestest curve reduction} that, after
replacing $R$ and $f$ by $h(R)$ and $h\circ f$, where $h$ is some
self-homeomorphism of $\tQ$ that is isotopic to the identity,
we may assume that the following conditions hold:
(1) 
$\partial R$
meets 
$\tau(\partial R)$ 
transversally,  (2) 
every component of $(\partial R) \cap\tau(\partial R )$ is a homotopically non-trivial
simple closed curve in $\tQ$, and (3) no component of
$R\cap\tau(\partial R )$ is a boundary-parallel annulus in $R$.

We set $U=R\cap\tau(R)$.
It follows from Condition (1) that  
$U$
is a compact $3$-manifold, and that if $F$ is a component of
$\partial U$, 
then either (i) $F$ contains a component
of $(\partial R)\cap\tau(\partial R)$, or (ii) $F$ is a component of
$\partial R$ or $\tau(\partial R)$. If (i) holds, then Condition (2)
implies that $F$ contains a simple closed curve which is homotopically
non-trivial in $\tQ$, and hence $F$ is not a $2$-sphere. If (ii) holds,
the incompressibility of $\partial R$ implies that $F$ is not  a
$2$-sphere. This shows:

\Claim\label{no spheres}
No boundary component of the compact $3$-manifold 
$U$ 
is a
$2$-sphere.
\EndClaim

Since the Euler characteristic of any compact orientable $3$-manifold
is half the Euler characteristic of its boundary, it follows from
\ref{no spheres} that every component of 
$U$
has non-positive Euler characteristic. Hence we may write
\Equation\label{either or}
U
=U\redsubone\discup\Uneg
\EndEquation
where $U\redsubone$ is the union of all components of 
$U$ 
having
Euler characteristic $0$, and $\Uneg$ is the union of all
components of strictly negative Euler characteristic. 
Either or both of the
$U_i$
 may
be empty. 

For $i=\redone ,\redtwo$, set $\Phi_i=\Fr_RU_i$. 
It follows from Condition (1) that $\Phi\redsubone\discup\Phi\redsubtwo=\Fr_R(R\cap\tau(R))=
R\cap\tau(\partial R)$. 
Furthermore, it follows from Conditions (1) and (3) that 
$\Phi_i$ is a two-dimensional submanifold of
$\tau(\partial R)$, 
that the boundary curves of $\Phi\redsubone\discup\Phi\redsubtwo$ are the
components of $(\partial R)\cap\tau(\partial R)$, and that:
\Claim\label{beirut} 
$\Phi_i\subset R$ is a properly embedded $2$-manifold, no component of
which is a boundary-parallel annulus in $R$.
\EndClaim

Condition
(2) asserts that every component of $\partial \Phi_i$
is a homotopically non-trivial curve in $\tQ$, and in particular in
$\tau(\partial R)$. Hence the subsurface $\Phi_i$
is $\pi_1$-injective in $\tau(\partial R)$. But the hypothesis of
the lemma
implies that $\tau(\partial R)$ is incompressible in $\tQ$, and
therefore:
\Claim\label{anatolia}
$\Phi_i$
is $\pi_1$-injective in $\tQ$, and in particular in $R$.
\EndClaim
The homotopic non-triviality of the components of $\partial \Phi_i$ in
$\tQ$ also implies:
\Claim\label{iz mir}
No component
of $\Phi_i$ is a disk. 
\EndClaim

The next step will be to apply 
the second assertion of Proposition \ref{new
  characterization} (together with Subsection \ref{see X});
we let $R$ and $U\redsubone$ play the respective roles of $X$ and $U$
in 
the second assertion of Proposition \ref{new
  characterization}.
The simple manifold
$R$ has non-empty boundary by the hypothesis of the present lemma, and
by  definition each component
of $U\redsubone$ has Euler characteristic $0$. The other hypotheses of
the second assertion of Proposition \ref{new
  characterization}
follow from 
\ref{no spheres}
and from
the case $i=\redone $ of \ref{beirut}, \ref
{anatolia}, and \ref {iz mir}.
Hence 
the  \abs\ $(C,W)=(C_R,W_R)$  of $R$
may be chosen within
its isotopy class in such a way that:
\Claim\label{now it's vertical}
$U\redsubone$ is contained in 
$W$,
is 
$(W,\redcale)$-vertical, and is disjoint from $\cale$.
\EndClaim


We now apply the observations in
Subsections 
\ref{how we use it}---\ref{first
  potato},
letting $R$ play the role of $\redPsi$.


As we observed in those subsections,
the manifold  $U=\reallyredR \cap\tau(\reallyredR )$ is
$\tau$-invariant; as in those subsections, we let $\sigma$ denote
the involution $\tau|U$. 
Now it
follows from the definition of $U_1$ and $U_2$ that each of them is
$\sigma$-invariant. Thus we are in the situation described in Subsection \ref{composition
  phenomenon}, with $R$ playing the role of $\redPsi$ in that
subsection, and with $U$, $U_1$, $U_2$, and
$\sigma$  defined as above.
As in that subsection,  we 
let $\sigma_i$ denote the PL involution $\sigma|U_i$ of $U_i$ for
$i=\redone,\redtwo$; we set $\Lambda=
\Lambda _{\reallyredR,\redUsigma}$ and 
$\lambda=\lambda_{\reallyredR,\redUsigma}:\reallyredR\to\Lambda$, and
we set 
$\redOmega_i=
\Lambda _{\reallyredR,U_i,\sigma_i}$ and 
$\lambda_i=\lambda_{\reallyredR,U_i,\sigma_i}:\reallyredR\to\redOmega_i$. 
As was observed in Subsection  \ref{composition
  phenomenon},
$\lambda_\redone$ restricts to a PL homeomorphism of $U_\redtwo$ onto a PL subset
$U_\redtwo\redstar$ of $\redOmega_1$, and 
conjugating the involution $\sigma_\redtwo$ of $U_\redtwo$
by the homeomorphism $\lambda_\redone|U_\redtwo:U_\redtwo\to U_\redtwo\redstar$, we obtain a PL
involution $\sigma_\redtwo\redstar$ of
$U_\redtwo\redstar$. Furthermore, according to the main observation in Subsection  \ref{composition
  phenomenon}, if we set 
 $\Lambda^\dagger=
\Lambda _{\redOmega_1, U_\redtwo\redstar,\sigma_\redtwo\redstar}$ and 
$
\redlambdatwostar
=\lambda_{\redOmega_1,
  U_\redtwo\redstar,\sigma_\redtwo\redstar}:\redOmega_1
\to\Lambda^\dagger$, then there is a unique PL homeomorphism $\mu:\Lambda^\dagger\to\Lambda$ such
that $\mu\circ
\redlambdatwostar
\circ\lambda_\redone=\lambda$.

It follows from the discussion in
\ref{second potato} 
that $\Lambda$, $\Omega_1$, $\Omega_2$ and $\Lambda^\dagger$ are
compact, connected $3$-manifolds.

Since $U$ is the full intersection $R\cap \tau(R)$, it
follows from
\ref{first potato}, 
applied with $R$ in the
role of $\redPsi$,
that there is
a unique homeomorphism $\theta: \Lambda \to
\redPi(R)$ such that $\redPi|R$, regarded as a map of $R$ onto
$\redPi(R)
\subset\inter Q
$, is equal to
$\theta\circ\lambda$. Let us set $\kappa=\theta\circ\mu$, so that
\Claim\label{compose 'em}
$\kappa:\Lambda^\dagger\to \redPi(R)$ is a  homeomorphism, and
$\kappa\circ\redlambdadagger\circ\lambda\redsubone=\redPi|R$.
\EndClaim

Let $j:\redPi(R)\to Q$ denote the inclusion map, and let
$\psi$  denote the map $j\circ\kappa\circ\redlambdadagger:\redRdagger\to Q$. 
The following commutative diagram of spaces summarizes some
information that will be used below:
\Equation\label{all them maps}
\xymatrix{
R
\ar[r]^{\lambda_1}
\ar@/^4pc/[rrrrr]^{\redPi|R}
&\redRdagger\ar[r]^{\redlambdatwostar}
\ar@/_5pc/[rrrr]^{\psi}
&\Lambda^\dagger\ar[r]^{\mu}_{\cong}
\ar@/_2pc/[rr]_{\cong}^\kappa
&\Lambda\ar[r]^-{\theta}_{\cong}&\redPi(R) \ar@{^{(}->}[r]^-{j}&Q
}
\EndEquation


Since the definition of $U_1$ implies that $\chi(U_1)=0$, we may apply
Lemma  \ref{gisty numbers}, 
with $R$, $U_1$ and $\sigma_1$ playing the
respective roles of $\redPsi$, $U$ and $\sigma$ (so that $\redRdagger$ plays the role of
$\Lambda$ in that lemma), 
to deduce that
\Equation\label{grand scheme}
\chibar(\redRdagger)=\chibar(R).
\EndEquation

Set $\redWdagger=\lambda\redsubone(W)$,
$\redCdagger=\lambda\redsubone(C)$, and $\redcaladagger=\lambda\redsubone(\redcale)$.
Since 
$(C,W)$ is in particular a splitting of $R$, we have
$R=C\cup W$, and 
hence
\Equation\label{the onion}
\redRdagger=\redCdagger\cup \redWdagger.
\EndEquation

On the other hand, 
the definition of $\redcale=\redcale_R$ (see \ref{see X}) gives $C\cap
W=\redcale$. By \ref{now it's vertical} we have
$U\redsubone\subset W-\redcale$; 
and
according to \ref{before gist}, $\lambda\redsubone|(R-U\redsubone)$ is injective and
$\lambda\redsubone^{-1}(\lambda\redsubone(U\redsubone))=U\redsubone$. It
follows that 
\Equation\label{ruggles of red cap}
\redCdagger\cap \redWdagger= \redcaladagger.
\EndEquation
The injectivity of $\lambda\redsubone|(R-U\redsubone)$, and the
fact that $U\redsubone$ is contained in  $W-\redcale=R-C$, also imply: 
\Claim\label{new redolence}
The map
$\lambda\redsubone|C$ is a  homeomorphism of $C$ onto $\redCdagger$. In particular, 
the pairs $(C,\redcale)$ and $(\redCdagger,\redcaladagger)$ are
homeomorphic, and hence 
$(\redCdagger,\redcaladagger)$ is an acylindrical pair.
\EndClaim


Since $U\redsubone\subset W$ by \ref{now it's vertical},
it follows from \ref{fever} that
$\redWdaggerstar\doteq \Lambda_{W,U\redsubone,\sigma\redsubone}$ 
and
$\lambda_{W,U\redsubone,\sigma\redsubone}:W\to \redWdaggerstar$
are defined, 
and that there is
a unique  homeomorphism $\zeta:\redWdaggerstar\to \redWdagger$ such that
$\zeta \circ
\lambda_{W,U\redsubone,\sigma\redsubone}
=\lambda\redsubone|W$.
In particular, if we set $\redcaladaggerstar=
\lambda_{W,U\redsubone,\sigma\redsubone}
(\redcale)$, 
the
pairs $(\redWdaggerstar,\redcaladaggerstar)$ and $(\redWdagger,\redcaladagger)$ are 
homeomorphic. But since \ref{now it's vertical} also gives that $U\redsubone$
is 
$(W,\redcale)$-vertical, 
and since
$\Phi_1=\Fr_RU_1$ is $\pi_1$-injective in $\tQ$ by \ref{anatolia}, it follows from
Lemma \ref{gist}, applied with
$\tQ$, $Q$, $\tau$, $\redcale$ and $W$ defined as above, and
with $U\redsubone$ playing the role of $U$ (so that $\sigma\redsubone$ plays the  role
of $\sigma$),
that $(\redWdaggerstar,\redcaladaggerstar)$ is a 
generalized bookish pair.
Hence:
\Claim\label{new oy}
$(\redWdagger,\redcaladagger)$ is a generalized bookish pair.
\EndClaim

It follows from \ref{the onion}, \ref{ruggles of red cap}, \ref{new
  redolence} and \ref{new
  oy} that
\Claim\label{batfish}
 $(\redCdagger,\redWdagger)$ is a generalized \abs\ of the compact, connected, orientable
$3$-manifold $\redRdagger$; that is,
$(\redRdagger,\redCdagger,\redWdagger)$ is a generalized \abt.
\EndClaim

According to \ref{verticacklepuss}, we may fix a strict
book of $I$-bundles
 $\calw$ 
\realiz ing the strict bookish pair
$(W,\redcale)$, such that 
$U$ is strongly
$\calw$-vertical.
Then 
Lemma \ref{gist} also gives a generalized book of $I$-bundles
$\redcalwdaggerstar$ realizing $(\redWdaggerstar,\redcaladaggerstar)$, 
such that 
$\calp_{\redcalwdaggerstar}$ is the homeomorphic image
under
$
\lambda_{W,U\redsubone,\sigma\redsubone}
$ of
a 
$\calw$-vertical submanifold of $W$ which is $\pi_1$-injective in
$W$; 
and such that,
for every element 
$F$ of $\calf_{\redcalwdaggerstar}
$, there is a
simple closed curve $\redUpsilon\subset W$, 
which is $\pi_1$-injective in $W$ and
is contained in either $\calb_\calw$ or in $\calp_\calw$, such
that 
$\redUpsilon_1^\diamondsuit\doteq \lambda_{W,U\redsubone,\sigma\redsubone}(\redUpsilon)$ 
is a simple closed curve isotopic to $F$, and 
$\lambda_{W,U\redsubone,\sigma\redsubone}|\redUpsilon:\redUpsilon\to
\redUpsilon_1^\diamondsuit$ 
is a covering map
of degree at most $2$.
Transporting the generalized book of $I$-bundles
$\redcalwdaggerstar$ via the homeomorphism $\zeta$ gives a generalized book of $I$-bundles
$\redcalwdagger$ realizing $(\redWdagger,\redcaladagger)$, and it now follows that:

\Claim\label{va avoy}
$\calp_{\redcalwdagger}$ is the homeomorphic image under
$\lambda_1$ of
a 
$\calw$-vertical 
submanifold of $W$ which is $\pi_1$-injective in
$W$;
 and 
for every element
$F$ of $\calf_{\redcalwdagger}$, there is a
simple closed curve $\redUpsilon\subset W$, 
which is $\pi_1$-injective in $W$ and
is contained in either $\calb_\calw$ or in $\calp_\calw$, such
that 
$\redUpsilon_1\doteq \lambda_1(\redUpsilon)$ 
is a simple closed
curve isotopic to $F$ in $\redcalwdagger$, and 
$\lambda_1|\redUpsilon:\redUpsilon\to \redUpsilon_1$ 
is a covering map
of degree at most $2$.
\EndClaim


The following fact will be useful.
\Claim\label{thirsty}
Let $G$ be a submanifold of $\redRdagger$ such that either
$G\in\calf_{\redcalwdagger}$, or $G=\redCdagger$, or
$G=\calp_{\redcalwdagger}$. 
Then there is a $\pi_1$-injective
submanifold $\redUpsilon$ of $R$ such
that 
$\redUpsilon_1\doteq \lambda_1(\redUpsilon)$ 
is a submanifold of
$\redRdagger$ isotopic to 
$G$, 
and either ($\alpha$)
$\lambda_1|\redUpsilon:\redUpsilon\to \redUpsilon_1$ 
is a homeomorphism, or ($\beta$)
$G\in\calf_{\redcalwdagger}$ and 
$\lambda_1|\redUpsilon:\redUpsilon\to \redUpsilon_1$ 
is a covering map
of degree $2$. Furthermore, in both cases,
$\psi|G:G\to Q$ is $\pi_1$-injective.
\EndClaim

To prove 
the first assertion of
\ref{thirsty}, first note that if 
$G=\redCdagger$ or
$G=\calp_{\redcalwdagger}$, it follows from 
\ref{new redolence}
or \ref{va avoy}   that $G$ is itself the homeomorphic image
under $\lambda_1$ of either $C$ or a
$\calw$-vertical submanifold of $W$ which is $\pi_1$-injective in
$W$.
Now $C$ is $\pi_1$-injective in $R$ by 
\ref{triad defs}.
Since $W$ is also $\pi_1$-injective in $R$ by
\ref{triad defs},
any 
$\calw$-vertical submanifold of $W$ which is $\pi_1$-injective in
$W$ is likewise $\pi_1$-injective in $R$. 
Thus the first assertion of \ref{thirsty} is established for the cases $G=\redCdagger$ and
$G=\calp_{\redcalwdagger}$. In the case $G\in\calf_{\redcalwdagger}$,
it follows from \ref{va avoy} that there is a $\pi_1$-injective
simple closed curve $\redUpsilon\subset W$ such
that 
$\redUpsilon_1\doteq \lambda_1(\redUpsilon)$ 
is a submanifold of
$\redRdagger$ isotopic to $G$, and either Condition ($\alpha$) or Condition
($\beta$) of \ref{thirsty} holds. Again using that $W$ is also
$\pi_1$-injective in $R$, we deduce that $\redUpsilon$ is
$\pi_1$-injective in $R$, and the first assertion of \ref{thirsty} is
established in this case as well.

To prove the second assertion of \ref{thirsty}, 
we may assume, in
light of the first assertion, that 
there is a $\pi_1$-injective
submanifold $\redUpsilon$ of $R$ such
that $\lambda_1(\redUpsilon)=G$, 
and either ($\alpha'$)
$\lambda_1|\redUpsilon:\redUpsilon\to \redUpsilon_1$ 
is a homeomorphism, or ($\beta'$)
$G\in\calf_{\redcalwdagger}$ and $\lambda_1|\redUpsilon:\redUpsilon\to G$ is a covering map
of degree $2$.
Now
$R$ is $\pi_1$-injective in $\tQ$ by the incompressibility of
$\partial R$, and  
$\redPi:\tQ\to Q$ 
is $\pi_1$-injective 
since $\redPi$ is a covering map. Hence 
$\redPi|\redUpsilon :\redUpsilon \to Q$ 
is
$\pi_1$-injective. But 
it follows from  
the commutativity of the diagram (\ref{all them maps})
that $\psi\circ\lambda_1|\redUpsilon =\redPi|\redUpsilon $. 
Thus $\psi\circ\lambda_1|\redUpsilon: \redUpsilon \to Q$ is
$\pi_1$-injective. 
If ($\alpha'$) holds, it follows that $\psi|G:G\to Q$ is
$\pi_1$-injective, which is the second assertion of \ref{thirsty}. If
($\beta'$) holds, it now follows that $(\psi|G)_\sharp:\pi_1(G)\to\pi_1( Q)$
restricts to an injection on an index-two subgroup of its domain
$\pi_1(G)$. But since $G\in\calf_{\redcalwdagger}$ in this case, $G$
is a simple closed curve and  $\pi_1(G)$ is therefore infinite cyclic. Hence
$(\psi|G)_\sharp:\pi_1(G)\to\pi_1( Q)$ is itself injective, and the
proof of \ref{thirsty} is complete.

For the rest of the proof, we distinguish two cases.

{\bf Case I: $U\redsubtwo=\emptyset$.}

In this case we have $U\redsubtwo\redstar=\emptyset$, and hence $\redlambdadagger$ is a 
homeomorphism.
By \ref{compose 'em}, it follows that
$\kappa\circ\redlambdadagger:\redRdagger\to \redPi(R)$
is a  homeomorphism.
Using the commutativity of the diagram (\ref{all them maps}), we
deduce that, if we
set 
$\redOmega\reddaggerprime=\redPi(R)$, $C\reddaggerprime=\redPi(C)$, $W\reddaggerprime=\redPi(W)$ and
$\redcale\reddaggerprime=\redPi(\redcale)$, then
$\kappa\circ\redlambdadagger$
maps $\redRdagger$, $\redCdagger$, $\redWdagger$ and $\redcaladagger$
homeomorphically onto $\redOmega\reddaggerprime$, $C\reddaggerprime$, $W\reddaggerprime$ and
$\redcale\reddaggerprime$ respectively.  
Using
\ref{ruggles of red cap},
\ref{new redolence}, and
\ref{batfish},
we now deduce:
\Claim\label{new wrap it up}
$(\redOmega\reddaggerprime, C\reddaggerprime, W\reddaggerprime)$ is a generalized \abt, homeomorphic to $(\redRdagger, \redCdagger, \redWdagger)$, and
$C\reddaggerprime\cap
W\reddaggerprime=\redcale\reddaggerprime$. Furthermore, 
$\kappa\circ\redlambdatwostar\circ\lambda_1$ (which is $\Pi|R$ regarded
as a map from $R$ to $\redRdagger'$)
maps $C$ and $\redcale$
homeomorphically onto 
$C_1'$ and $\redcale_1'$
respectively.
\EndClaim



We may transport $\redcalwdagger$ via the homeomorphism $
(\kappa\circ\redlambdadagger)|\redWdagger:\redWdagger\to W\reddaggerprime$ to obtain
a generalized book of $I$-bundles $\calw\reddaggerprime$ \realiz ing
the generalized bookish pair
 $(W\reddaggerprime,\redcale\reddaggerprime)$.
According to \ref{thirsty}, the restrictions of $\psi$ to
both $\redCdagger$ and $\calp_{\redcalwdagger}$, 
and to each element of 
$\calf_{\calw_1}$,  
are
$\pi_1$-injective. Since $\psi$ maps $\redCdagger$ and 
$\calp_{\redcalwdagger}$
homeomorphically onto $
C\reddaggerprime
$ and
$\calp_{\calw\reddaggerprime}$ respectively, and maps each element of
$\calf_{\calw_1}$ 
homeomorphically onto an element  of $\calf_{\calw\reddaggerprime}$, we deduce that 
\Claim\label{sarah}
$C\reddaggerprime$ and
$\calp_{\calw\reddaggerprime}$ are $\pi_1$-injective in $Q$, and each
element of $\calf_{\calw\reddaggerprime}$ is also $\pi_1$-injective in $Q$.
\EndClaim

According to \ref{new wrap it up}, the pair 
$(C_1',\redcale_1')$ 
is acylindrical. According to \ref{acylindrical}, it follows that no
component of 
$\partial C_1'$ 
is a sphere or torus. Since
$C\reddaggerprime$ is $\pi_1$-injective in $Q$ by \ref {sarah}, we deduce:
\Claim\label{by the way}
No component of 
$C_1'$ 
is contained in a solid torus in $Q$.
\EndClaim

We now claim:
\Claim\label{convergence vague}
There is a submanifold 
$W_1^+$ 
of $\inter Q$ with the following properties:
(\trulyredA) $\redWdaggerplus\supset W\reddaggerprime$, (B) $\partial \redWdaggerplus$ is a union of
components of $\partial W\reddaggerprime$, (C) each component of
$\overline{\redWdaggerplus-W_1'}$ 
is a solid torus,
(D) $\partial W_1^+\supset\redcale_1'$, 
and (E)
the $3$-manifold pair 
$(W_1^+,\redcale_1')$ 
is a strict
bookish pair.
\EndClaim

To prove \ref{convergence vague}, we first note that $(W\reddaggerprime,\redcale\reddaggerprime)$ is a generalized bookish pair by \ref{new wrap it
  up}, and is \realiz ed by $\calw\reddaggerprime$; that $W_1'\subset
\redRdagger'\subset\inter Q$; and that
$\calp_{\calw\reddaggerprime}$ and the
elements of $\calf_{\calw\reddaggerprime}$ are $\pi_1$-injective in
$Q\supset W_1'$ by \ref{sarah}. Thus all the 
hypotheses of Lemma \ref{like 2.5} hold with
$(W\reddaggerprime,\redcale\reddaggerprime)$ playing the role of
$(W,\redcale)$. Hence 
one of the alternatives ($\redalpha$) or ($\redbeta$) of Lemma \ref{like 2.5} holds when
$(W,\redcale)$ is replaced by
$(W\reddaggerprime,\redcale\reddaggerprime)$.

Suppose that Alternative ($\redalpha$) holds, i.e. that there are distinct components
$V_1$ and $V_2$ of 
$W_\reddenedone'$ 
such that the
component $\Xi$ of $\overline{Q-V_1}$ containing $V_2$ is a solid torus 
whose boundary is contained in $V_1$. 
Then in particular $W_1'$ is disconnected; since $\redRdagger'$ is connected,
we have $W_1'\ne \redRdagger'$. Since $(\redRdagger',C_1',W_1')$ is a $3$-manifold
triad by \ref{new wrap it up}, we have $\redRdagger'=C_1'\cup W_1'$, and the
connectedness of $\redRdagger'$ now implies that each component of $W_1'$
meets some component of $C_1'$. In particular,
$V_2$ meets some component of 
$C_1'$, 
and this
component of  
$C_\reddenedone'$ 
must be contained in the solid torus $\Xi$;
but this contradicts \ref{by the way}. 
Hence
Alternative ($\redbeta$) of Lemma
\ref{like 2.5} must hold.

Thus there is a $3$-manifold pair
$(W_1^+,\redcale_1^-)$ 
such that Conditions
(1)--(5) of the conclusion of Lemma \ref{like 2.5} hold with
$(W_1^+,\redcale_1^-)$ 
and
$(W\reddaggerprime,\redcale\reddaggerprime)$ in place of $(W^+,\redcale^-)$
and $(W,\redcale)$. 

According to Condition (5), 
$\redcale_1^-$ 
is the union of
all components of 
$\redcale_1'$
contained in 
$\partial
W_1^+$. 
Suppose that $A$ is a component of $\redcale_\reddenedone'$
which is not contained in 
$\partial
W_1^+$. 
Then it follows from Conditions (1) and (3) that $A$ is
contained in the boundary of some component $L$ of 
$\overline{W_\reddenedone^+-W_\reddenedone'}$, and that $L$ is a
component of 
$\overline{Q-W_\reddenedone'}$. Since, by  \ref{new wrap it up},
$(\redOmega\reddaggerprime, C\reddaggerprime, W\reddaggerprime)$ is a
$3$-manifold triad and
$C\reddaggerprime\cap W\reddaggerprime=\redcale\reddaggerprime$, we must 
have $A\subset J\subset L$ for some component $J$ of
$C\reddaggerprime$. But $L$ is a solid torus by Condition (4), and the
inclusion $J\subset L$ therefore contradicts \ref{by the way}. It
follows that $\redcale_\reddenedone^-=\redcale_\reddenedone'$.

Thus the  $3$-manifold pair $(W_\reddenedone^+,\redcale_\reddenedone^-)$ may be
renamed $(W_\reddenedone^+,\redcale_\reddenedone')$. Note that Conditions (1), (3)
and (4) of the conclusion of Lemma \ref{like 2.5} respectively imply
that $W_\reddenedone^+$ has the properties (\trulyredA), 
(B) and (C) stated in \ref{convergence vague}. Since
$(W_\reddenedone^+,\redcale_\reddenedone')$ is a $3$-manifold pair, $W_\reddenedone^+$
has Property (D). It remains to verify Property
(E). If $W_\reddenedone^+$ fails to have this property, it follows from Condition (2) of
Lemma \ref{like 2.5} that $(W_\reddenedone^+,\redcale_\reddenedone')$ has a
degenerate component.

Suppose that $(\Theta_\reddenedone^+,E_\reddenedone')$  is a degenerate component of
$(W_\reddenedone^+,\redcale_\reddenedone')$. Then $\Theta_\reddenedone^+$ is a solid torus,
$E_\reddenedone'\subset\partial \Theta_\reddenedone^+$ is a single annulus, and the inclusion
homomorphism $H_1(E_\reddenedone',\ZZ)\to H_1(\Theta_\reddenedone^+,\ZZ)$ is an isomorphism.

According to \ref{new wrap it up}, 
$\Pi$
maps $\redcale$ homeomorphically onto $\redcale_\reddenedone'$. Hence there is a
unique component $E$ of $\redcale$ such that $\Pi(E)=E_\reddenedone'$, and
$\Pi|E:E\to E_\reddenedone'$ is a homeomorphism. Let $\Theta$ denote the
component of $W$ containing $E$. 
Then
$\Pi(\Theta)\subset\Theta_1^+$. Since $E_\reddenedone'$ is the
only component of $\redcale_\reddenedone'$ contained in $\Theta_\reddenedone^+$, and
since $\Pi|\redcale:\redcale\to\redcale_\reddenedone'$ is a homeomorphism, $E$ is
the only component of $\redcale$ contained in $\Theta$; that is,
$(\Theta,E)$ is a component of $(W,\redcale)$. 

Now it follows from Proposition \ref{yes they are} that the component $\Theta$ of $W$ is  $\pi_1$-injective in $R$,
and
$R$ is in turn $\pi_1$-injective in $\tQ$ by the incompressibility of
$\partial R$. Hence the covering map $\Pi$ restricts to a $\pi_1$-injective
map from $\Theta$ to $Q$. In particular, $\Pi|\Theta$ may be regarded
as a $\pi_1$-injective
map  $\eta:\Theta\to \Theta_1^+$. Since
$\Theta_1^+$ is a solid torus, $\pi_1(\Theta)$ is cyclic. Since $R$ is
simple, and the components of $\Fr_R\Theta$ are essential annuli by
Proposition \ref{yes they are}, $\Theta$ is irreducible and is not a
ball; hence $\Theta$ is a solid torus.


Let us fix a generator $\gamma$ of
$H_1(\Theta,\ZZ)$. The definition of a strict bookish
pair also implies that $(\Theta,E)$ is non-degenerate; and since
$\Theta$ is a solid torus, and $E$ is a single annulus, it now follows
that some oriented core curve of $E$ represents the element $n\gamma$
of $H_1(\Theta,\ZZ)$, where $n$ is an integer strictly greater than $1$.

Since $\eta$
maps $E$ homeomorphically onto $E_\reddenedone'$, 
some oriented core curve of $E$ represents the element $n\eta_*(\gamma)$
of $H_1(\Theta_\reddenedone^+,\ZZ)$. In particular, the element of 
$H_1(\Theta_\reddenedone^+,\ZZ)$ represented by this core curve is
divisible by $n>1$. This contradicts the fact that the inclusion
homomorphism $H_1(E_\reddenedone';\ZZ)\to H_1(\Theta_\reddenedone^+;\ZZ)$ is an
isomorphism. The proof of \ref{convergence vague} is now complete.


Since $(\redOmega\reddaggerprime, C\reddaggerprime, W\reddaggerprime)$ is a
$3$-manifold triad by \ref{new wrap it up}, we have
$C_\reddenedone'\subset\overline{Q-W_\reddenedone'}$; and the properties (\trulyredA),
(B) and (C) of $W_\reddenedone^+$ stated in \ref{convergence vague} imply
that each component of $\overline{Q-W_\reddenedone'}$ is either a component
of $\overline{Q-W_\reddenedone^+}$ or a solid torus. Since, according to
\ref{by the way}, no component of $C_\reddenedone'$ can be contained in a
solid torus, we must have $C_\reddenedone'\subset\overline{Q-W_\reddenedone'}$,
so that $C_\reddenedone'\cap W_\reddenedone^+=C_\reddenedone'\cap
W_\reddenedone'=\redcale_\reddenedone'$. If we set $\redRdagger^+=C_\reddenedone'\cup
W_\reddenedone^+
\subset\inter Q
$, it 
follows that:
\Claim\label{old college try}
 $(\redRdagger^+,C_\reddenedone',
W_\reddenedone^+)$ is a compact,  orientable $3$-manifold triad, $\redRdagger'\subset
\redRdagger^+$, 
and  every component of $\overline{\redRdagger^+-\redRdagger'}$ is a solid torus
whose boundary is contained in $\partial \redRdagger'$. 
\EndClaim
Since the $3$-manifold triad $(\redOmega\reddaggerprime,
C\reddaggerprime, W\reddaggerprime)$ is compact, connected and orientable
by \ref{new wrap it up}, it follows from \ref{old college try} that $(\redRdagger^+,C_\reddenedone',
W_\reddenedone^+)$ is
also a compact, connected,  orientable $3$-manifold triad. 

It also follows from \ref{old college
  try}  that
$\chibar(\redRdagger^+)=\chibar(\redOmega\reddaggerprime)=\chibar(\redRdagger)$,
which with \ref{grand scheme} gives
$\chibar(\redRdagger^+)=\chibar(R)$.

But
$(C_\reddenedone',\redcale_\reddenedone')$ is an acylindrical pair by 
\ref{new wrap it up}, and $(W_\reddenedone^+,\redcale_\reddenedone')$ is a strict
bookish pair according to Property (E) of $W_\reddenedone^+$; hence the
$3$-manifold triad $(\redRdagger^+,C\reddaggerprime,\redWdaggerplus)$ is
\anabt. To summarize:
\Claim\label{square onion}
$(\redRdagger^+,C\reddaggerprime,\redWdaggerplus)$ is \anabt, and $C_1'\cap
W_1^+=\redcale\reddaggerprime$. Furthermore, we have $\chibar(\redRdagger^+)=\chibar(R)$.
\EndClaim


Next, we claim:

\Claim\label{easy way out}
$\redRdagger^+$ is not a solid torus.
\EndClaim

To prove \ref{easy way out}, note that  $\redRdagger^+\supset
\redRdagger'=\Pi(R)\supset\Pi(f(K))$. Note also that $f$ is $\pi_1$-injective
by hypothesis, and that the covering map $\Pi$ is
also $\pi_1$-injective. Hence $\Pi\circ f$ is $\pi_1$-injective as a
map from $K$ to $Q$, and in particular is $\pi_1$-injective when
regarded as a map from $K$ to $\redRdagger^+
$. If $\redRdagger^+$ were a solid torus, it
would follow that $\pi_1(K)$ is isomorphic to a subgroup of an
infinite cyclic group, a contradiction to the hypothesis that $K$ is
freely indecomposable. Thus \ref{easy way out} is proved.



  Since $f(K)\subset R$, we have $(\redPi\circ f)(K)\subset
\redPi(R)=\redOmega\reddaggerprime\subset \redRdagger^+$. 
Thus:
\Equation\label{new rap around the clap}
(\redPi\circ f)(K)\subset \redRdagger^+.
\EndEquation

Consider
the subcase in which $\partial \redRdagger^+$ is incompressible in $Q$.
In this subcase, it follows from Remark \ref{more fruit} (and the
connectedness of $\redRdagger^+$, which has been established and is
included in the statement of \ref{square onion}) that $\redRdagger^+$ is simple. According to  \ref{square onion},  $(C\reddaggerprime,\redWdaggerplus)$ is
\anabs\ of the simple manifold $\redRdagger^+$, and $C\reddaggerprime\cap \redWdaggerplus=\redcale\reddaggerprime$. It
therefore follows from  Proposition
\ref{new characterization} and Subsection \ref{see X} that the splitting
$(C\reddaggerprime,\redWdaggerplus)$ of $\redRdagger^+$ is isotopic to 
$(C_{\redRdagger^+},W_{\redRdagger^+})$. In particular the pair $(C_{\redRdagger^+},\redcale_{\redRdagger^+})$ is
homeomorphic to $(C\reddaggerprime,\redcale\reddaggerprime)$, and therefore to $(C,\redcale)$ by
\ref{new wrap it up} and 
\ref{new redolence}. Since \ref{square onion} gives
$\chibar(\redRdagger^+)=\chibar(R)$, and since
$\partial
\redRdagger^+$ is
incompressible, we deduce 
that Alternative (i) of the
conclusion of the present lemma holds in this subcase if we set $T=\redRdagger^+$.

To complete
the proof in Case I, it remains to consider the subcase in which
$\partial \redRdagger^+$ is compressible in $Q$. In this 
subcase, we apply
Lemma \ref{surgery lemma}, 
with $Q$ and $K$ as above, taking
$T_0$ to be the connected manifold $\redRdagger^+$, and setting
$g_0=\redPi\circ f$. 
We have $(\redPi\circ f)(K)\subset \redRdagger^+$ 
by 
(\ref{new rap around
  the clap}),
and since $f$ is $\pi_1$-injective by hypothesis and the
covering map $\redPi$ is $\pi_1$-injective, $\redPi\circ f$ is also
$\pi_1$-injective.
Since $\redRdagger^+$ admits the \abs\ 
$(C\reddaggerprime,\redWdaggerplus)$ according to \ref{square onion}, and is not a solid
torus 
according to \ref {easy way out},
it follows from
Proposition \ref {no nothing} that
no component of $\partial \redRdagger^+$ is a sphere or a torus. 
Since in
addition $\partial \redRdagger^+$ is compressible, 
Lemma \ref{surgery lemma} gives
a 
compact, connected three-dimensional 
submanifold  $T$  of 
$\inter\redQ $ 
such that 
(a) $\partial T$ is incompressible in $\redQ $ (and hence $T$ is
simple
by Remark \ref{more fruit}),
(b$'$) 
$\chibar(T)<\chibar(\redRdagger^+)$, and
(c) 
  $\rk_\redp(T)\le \rk_\redp(\redRdagger^+)$
for every prime $p$,
and   a \redPL map 
$g :K\to\inter \redQ $ 
such that (e) $g$ is homotopic to $\redPi\circ
f$ and $g (K)\subset T$.

Since 
$\chibar(\redRdagger^+)=\chibar(R)$ by \ref{square onion}, Property (b$'$)
of $T$ 
gives 
\Equation\label{harnack}
\chibar(T)<\chibar(R).
\EndEquation

Now since 
$(\redRdagger^+,C\reddaggerprime,\redWdaggerplus)$ is \anabt\ by \ref {square onion}, and $\redRdagger^+$ is not
a solid torus by \ref{easy way out},  the
second assertion of Proposition \ref{no cylinders} gives $\rk_p(\redRdagger^+)
\le
\rk_p(C\reddaggerprime)+3\chibar(\redRdagger^+)$
for every prime $p$.
But Property (c) of $T$ gives   $\rk_\redp(T)\le \rk_\redp(\redRdagger^+)$
for every prime $p$;
by \ref{square onion}
we have $\chibar(\redRdagger^+)=\chibar(R)$; and by \ref{new wrap it up} and \ref{new
  redolence}, 
$C\reddaggerprime$ is homeomorphic to $C$.
Hence
\Equation\label{binz}
\rk_\redp(T)\le \rk_\redp(C)+3\chibar(R)
\EndEquation
for every prime $p$.

Recalling from 
\ref{keep it simple}
that the simplicity of $R$ and $T$ implies that $\chibar(R)$
and $\chibar(T)$ are non-negative, we see that (\ref{binz}) is
stronger than the inequality $\rk_\redp(T)\le
\rk_\redp(C_R)+4\chibar(R)+\chibar(T)$ appearing in Alternative (ii) of
the conclusion of the present lemma. 
With the inequality (\ref{harnack}), this shows that Alternative (ii)
holds in the subcase where $\partial \redRdagger^+$ is compressible, and the
proof in Case I is complete.


{\bf Case II: $U\redsubtwo\ne\emptyset$.}

In this case we will denote the  number of components of $U_2\redstar$ that are not
$\sigma_2\redstar$-invariant (which is clearly even) by $2\mni$, and we will denote the number of components of $U_2\redstar$ that are 
$\sigma_2\redstar$-invariant by $\minv$. We set $m=\mni+\minv$. Since $U_2\redstar$ is by definition
homeomorphic to $U_2$, and is therefore non-empty in this case, at
least one of the quantities $\mni$ and $\minv$ is strictly
positive. Hence
\Equation\label{om}
m>0.
\EndEquation

According to the definition of $U_2$, each component of $U_2$ (and
hence each component of $U_2\redstar$) has strictly negative Euler
characteristic. Furthermore, each $\sigma_2\redstar$-invariant component of
$U_2\redstar$ admits a fixed-point free involution, and therefore has even
Euler characteristic. It follows that
\Equation\label{goner}
\chibar(U_2\redstar)\ge 2m.
\EndEquation

We now apply Lemma \ref{gisty numbers}, letting $\redRdagger$, $U_2\redstar$ and
$\sigma_2\redstar$ play the respective roles of $\redPsi$, $U$ and $\sigma$
(so that $\Lambda^\dagger$ and $\redlambdadagger$ play the roles of $\Lambda$ and
$\lambda$). 
The first assertion of Lemma \ref{gisty numbers} gives
$\chibar(\Lambda^\dagger)=\chibar(\redRdagger)-\chibar(U_2\redstar)/2$, which with
(\ref{goner}) implies that $\chibar(\Lambda^\dagger)\le\chibar(\redRdagger)-m$. Since
$\chibar(\redRdagger)=\chibar(R)$ by \ref{grand scheme}, 
and since $\Lambda^\dagger$ and $\redPi(R)$ are homeomorphic by \ref{compose 'em},
this may be rewritten
as
\Equation\label{I'll give you what happened}
\chibar(\redPi(R))\le\chibar(R)-m.
\EndEquation
From (\ref{om}) and (\ref{I'll give you what happened}) we deduce
\Equation\label{progress illuminating asia}
\chibar(\redPi(R))<\chibar(R).
\EndEquation

For the remainder of the proof in Case II, we let $p$ denote an
arbitrary prime.
The second assertion of Lemma \ref{gisty numbers} 
gives:
\Claim\label{oh that beautiful number}
The 
codimension of $\lambda_*(H_1(\redRdagger;\FF_\redp))$
in $H_1(\Lambda^\dagger;\FF_\redp)$ is at most $m$.
\EndClaim

The next step is an application of Proposition \ref{use
  cylinders}. According to \ref{batfish},
$(\redRdagger,\redCdagger,\redWdagger)$ is \anabt, and by
\ref{ruggles of red cap}
we have $\redCdagger\cap \redWdagger= \redcaladagger$.

The generalized book of $I$-bundles
$\redcalwdagger$ \realiz es $(\redWdagger,\redcaladagger)$.
According to \ref{thirsty}, $\psi|\jay:\jay\to Q$ is
$\pi_1$-injective for every $\jay\in\calf_{\redcalwdagger}$. Thus all the
hypotheses of Proposition \ref{use cylinders} hold if we let $(\redRdagger,\redCdagger,\redWdagger)$
and $\redcalwdagger$ play the respective roles of $(V,C,W)$
and $\calw$, 
and define $Q$ and $\psi$ as above.
It therefore
follows 
from Proposition \ref{use cylinders} that
\Equation\label{tunnel of lerf}
\rk_\redp(\psi;\redRdagger,Q)
\le \max(1,\rk_\redp(\redCdagger)+3\chibar(\redRdagger)).
\EndEquation

Since $\redCdagger$ is homeomorphic to $C$ by
\ref{new redolence},
 and $\chibar(\redRdagger)=\chibar(R)$ by (\ref{grand scheme}),
we may rewrite (\ref{tunnel of lerf}) in the form
$\rk_\redp(\psi;\redRdagger,Q)
\le\max(1, \rk_\redp(C)+3\chibar(R))$. But 
since the 
simple manifold $R$ has a non-empty boundary by
hypothesis, we have $\chibar(R)>0$ by
\ref{keep it simple}.
Hence
\Equation\label{funnel of smurf}
\rk_\redp(\psi;\redRdagger,Q)
\le \rk_\redp(C)+3\chibar(R).
\EndEquation

Consider the  diagram
$$\xymatrix{
H_1(\redRdagger;\FF_\redp)
\ar[rd]^{\psi_*}\ar[dd]_{\lambda_*}&\\
&
H_1(Q;\FF_\redp)
\\
H_1(\Lambda^\dagger;\FF_\redp)
\ar[ru]_{(j\circ\kappa)_*}&
}
$$
which is commutative by the definition of $\psi$ (cf. (\ref{all them maps})). 
The images of $\psi_*$ and $(j\circ\kappa)_*$ are subspaces of the
vector space $H_1(Q;\FF_\redp)$. By commutativity we have
$\image\psi_*\subset\image(j\circ\kappa)_*$. Since, according to \ref
{oh that beautiful number}, $\image\lambda_*$ has codimension at most
$m$ in 
 $H_1(\Lambda^\dagger;\FF_\redp)$, the codimension of $\image\psi_*$ in
 $\image(j\circ\kappa)_*$  is also at most $m$. But according to
 (\ref{funnel of smurf}), the dimension of $\image\psi_*$ is at most
$\rk_\redp(C)+3\chibar(R)$. 
Hence the dimension of
$\image(j\circ\kappa)_*$ is at most 
$\rk_\redp(C)+3\chibar(R)+m$. 
Now
since $\kappa$ is a homeomorphism by \ref{compose 'em}, and $j:\redPi(R)\to
Q$ is the inclusion map, the dimension of
$\image(j\circ\kappa)_*$ is equal to $\rk_\redp(\redPi(R),Q)$ by the
definitions given in \ref{redirects here}. Since we have 
$m\le\chibar(R)-\chibar(\redPi(R))$ by (\ref
{I'll give you what happened}),
we now deduce that
\Equation\label{van der beck}
\rk_\redp(\redPi(R),Q)\le
\rk_\redp(C)+4\chibar(R)-\chibar(\redPi(R)).
\EndEquation

We now apply Lemma \ref{surgery lemma},
taking $K$ and $Q$ as above,
and setting $g_0=\redPi\circ f$ and $T_0=\redPi(R)\subset\inter
Q$. Since $R$ is connected, $T_0$ is also connected. The inclusion
$g_0(R)\subset T_0$ is an equality in this context. The
$\pi_1$-injectivity of $g_0$ follows from the hypotheses that $f$ is
$\pi_1$-injective and that $\redPi$ is a covering map. The manifold
$T_0=\redPi(R)$ is covered by $R\cup\tau(R)$, which by \ref{no spheres} has
no sphere boundary component; hence no boundary component of $T_0$ is
a sphere. Lemma \ref{surgery lemma} therefore gives 
a compact submanifold  $T$  of $\inter\redQ $ and a \redPL map $\redg :K\to\inter \redQ $  such that 
(a)
$\partial T$ is incompressible in $\redQ $;
(b) $\chibar(T)\le\chibar(\redPi(R))$;
(d) $\rk_\redp(T,Q)\le\rk_\redp(\redPi(R),Q)$; 
and (e) 
$g$ is homotopic to  $\redPi\circ f$, and $g (K)\subset T$. 
Combining Properties (b) and (d) of $T$ with (\ref{progress
  illuminating asia}) and (\ref{van der beck}), we obtain:
\Equation\label{pour chanter en choeur}
\chibar(T)<\chibar(R)\text{ and }\rk_\redp(T,Q)\le\rk_\redp(C)+4\chibar(R)- \chibar(T).
\EndEquation

Next, combining the second inequality in (\ref{pour chanter en
  choeur}) with Lemma \ref{insurance}, we obtain
\Equation\label{ah mais c'est toi}
\rk_\redp(T)\le\rk_\redp(T,Q)+\tg(\partial
T)\le\rk_\redp(C)+4\chibar(R)- \chibar(T) +\tg(\partial T),
\EndEquation
where $\tg(\partial T)$ is defined by \ref{tgdef}. Since $\partial T$ is incompressible in the simple manifold $\redQ $
according to property (a) of $T$, each component of $\partial T$ is a
closed orientable surface of genus at least $2$; this implies that
$\tg(\partial T)\le\chibar(\partial T)=2\chibar(T)$. From (\ref{ah
  mais c'est toi}) we therefore obtain
\Equation\label{oui c'est moi}
\rk_\redp(T)\le\rk_\redp(C)+4\chibar(R)+
\chibar(T).
\EndEquation

From properties (a) and (e) of $T$, together with (\ref{oui c'est
  moi}) and the first inequality in (\ref{pour chanter en choeur}), we
see that Alternative (ii) of the conclusion of the lemma always holds
in Case II.
\EndProof

\section{Towers and the proof of the main theorem}\label{tower section}

\Number\label{ads def}
Let $K$ and $L$ be finite simplicial complexes, and let $\phi : K\to
L$ be a simplicial map. As in \cite[Definition 5.1]{acs-singular}, we
define the {\it degree of singularity} of $\phi$, denoted  $\DS(\phi)$, to be the number of ordered pairs $(v, w)$ of vertices of $K$ such that $v\ne w$ but $\phi(v) = \phi(w)$.
If $f$ is any PL map from a compact PL space $X$ to a PL space $Y$,
then as in  \cite[Definition 5.1]{acs-singular}, we define the {\it absolute degree of singularity} of $f$, denoted $\ADS(f)$, by
$\ADS(f ) = \min_{K,L,\phi} \DS(\phi)$, 
where $(K,L,\phi )$ ranges over all triples such that $K$ and $L$
are simplicial complexes whose respective underlying PL spaces are $X$
and $f(X)$, and $\phi:K\to L$ is a simplicial map whose underlying PL
map is
$f : X\to f(X)$.

\EndNumber

\Number\label{tower defs}
As in \cite[Definition 8.1]{acs-singular}, we define a {\it  tower} of
$3$-manifolds, to be a $(3n + 2)$-tuple
$\calt = (M_0,N_0,\Pi_1,M_1,N_1,\Pi_2,\ldots,\Pi_n,M_n,N_n)$, where $n$
(called the {\it height} of the tower)  is a non-negative integer,
$M_0,\ldots, M_n$ are compact, connected, orientable  $3$-manifolds,
$N_j$ is a compact, connected 3-dimensional PL submanifold of $M_j$
for $j = 0,\ldots, n$, and $\Pi_j :M_j \to N_{j-1}$ is a covering map for
$j=1,\ldots,n$. (We have used $\Pi_j$ here in place of the notation
$p_j$ which is used in  \cite{acs-singular}, to
avoid confusion with the use of $p$ for a prime elsewhere in this
section.) As in \cite{acs-singular}, we shall refer to $M_0$ as  the
{\it base} of the tower $\calt$ and to $N_n$ as its {\it top}. 

If $\iota_j : N_j \to M_j$ denotes the inclusion map for $j =
0,\ldots,n$, the map $h_\calt=\iota_0 \circ \Pi_1 \circ\iota_1 \circ\Pi_2 \circ\cdots\circ\Pi_n
\circ\iota_n :N_n\to M_0$ is called the {\it tower map} associated to
the tower $\calt$.
If $\phi$ is a PL map from a compact PL space $K$ to the base $M_0$ of
the tower $\calt$, then, as in \cite[Definition 8.7]{acs-singular}, we
define a {\it homotopy-lift} of $\phi$ through the tower $\calt$ to be
a PL map $\tphi:K\to N_n$ such that $h_\calt\circ\tphi$ is homotopic to
$\phi$.
The homotopy-lift $\tphi$ of $\phi$ is termed {\it good} if the inclusion homomorphism $\pi_1(\tphi(K)) \to\pi_1(N_n)$ is surjective.

As in \cite{acs-singular}, if $\calt =
(M_0,N_0,\Pi_1,M_1,N_1,\Pi_2,\ldots,\Pi_n,M_n,N_n)$ is a tower, we
shall refer to any tower of the form
$(M_0,N_0,\Pi_1,M_1,N_1,\Pi_2,\ldots,\Pi_m,M_m,N_m)$, with $m\le n$,
as a {\it truncation} of $\calt$; and we shall say that a tower
$\calt^+$ is an {\it extension} of a tower $\calt$ if $\calt$ is a
truncation of $\calt^+$.
\EndNumber

\Number\label{from acs}
According to \cite[Subsection 8.2]{acs-singular}, if 
$$
\calt=(M_0,N_0,\Pi_1,M_1,N_1,\Pi_2,...,\Pi_n,M_n,N_n)
$$
is a tower, then
there is an index $j_0$  with $ 0\le j_0 \le n$ such that $M_j$ is
closed when $0\le j\le j_0$ and $\partial M_j\ne\emptyset$ when $j_0
<j\le n$. Furthermore, in the latter case, for each $j<j_0$ we have
$N_j=M_j$; in particular, for each $j<j_0$, the map $\Pi_{j+1}$ is a
covering map between the closed manifolds $M_{j+1}$ and $M_j$.
\EndNumber

\Number\label{C'est si bon}
As in \cite[Definition 8.4]{acs-singular}, we shall say that a tower 
$\calt = (M_0,N_0,\Pi_1,M_1,N_1,\Pi_2,...,\Pi_n,M_n,N_n)$  is {\it good} if it has the following properties:
\begin{itemize}
\item $M_j$ and $N_j$ are irreducible for $j=0,\ldots,n$;
\item $\partial N_j$ is a (possibly empty) incompressible surface in
  $\inter M_j$ for $j=0,\ldots, n$;
\item the covering map $\Pi_j : M_j\to N_{j-1}$ has degree 2 for $j = 1,\ldots,n$; and
\item for each $j$ with $2 \le j \le n$ such that $M_j$ is closed, the
  map 
$\Pi_j\circ\Pi_{j-1} : M_j\to M_{j-2}$, which is a four-fold covering
map by \ref{from acs}, defines a regular covering space whose deck
group is isomorphic to $\ZZ/2\ZZ\times \ZZ/2\ZZ$.
\end{itemize}
\EndNumber

\Proposition\label{perfect is good}
Suppose that
$\calt = (M_0,N_0,\Pi_1,M_1,N_1,\Pi_2,...,\Pi_n,M_n,N_n)$  is a good
tower, and that its base $M_0$ is simple. Then $M_j$ and $N_j$ are
simple for $j=0,\ldots,n$.
\EndProposition

\Proof
According to the definition of a good tower, $\partial N_j$ is
incompressible in $\inter M_j$ for each $j$. The definition also
implies that $N_j$ is irreducible, and in particular connected. Hence by Remark \ref{more
  fruit}, if $M_j$ is simple for a given $j$ then $N_j$ is also
simple. On the other hand, if $N_j$ is simple for a given $j<n$ then
the compact covering space $M_{j+1}$ of $N_j$ is simple by \ref{keep
  it simple}. Since $M_0$ is
simple by hypothesis, the result now follows by induction.
\EndProof

\Lemma\label{new eight point nine}
Let $K$ be a compact, connected PL space such that $\pi_1(K)$ is
freely indecomposable. Let 
$$\calT=(M_0,N_0,\Pi_1,...,N_n)$$
be a good tower of $3$-manifolds.
Assume that either 
(\wasB) $n\ge1$ and $\rk_2(N_{n-1})
\ge\rk_2(K)+2$, or 
(\wasA) $\partial N_n\ne\emptyset$, or 
(C) $n=0$.
Let $\phi : K\to Y$ be
a $\pi_1$-injective PL map,
 and let $\tphi:K\to N$ be a good
  homotopy-lift of $\phi$ through the tower
 $\calT$. Assume that $\tphi_*:H_1(K,\FF_2)\to H_1 (N_n;\FF_2)$ is not
 surjective. Then there exist a height-$(n + 1)$ extension $\calT^+$
 of $\calT$ which is a good tower, and a
good
homotopy-lift $\tphi^+$ of $\phi$ through the tower $\calT^ +$, such that $\ADS(\tphi^+) < \ADS(\tphi)$.
\EndLemma

\Proof
This is an application of
Lemma 8.8 
of \cite{acs-singular}. In that
lemma, one considers a compact, connected PL space $K$ such that
$\pi_1(K)$ is freely
indecomposable, a good tower of $3$-manifolds
$\calt=(M_0,N_0,\redPi_1,\ldots,N_n)$, a $\pi_1$-injective \redPL map $\phi:K\to
M_0$, and a good homotopy-lift $\tphi$ of $\phi$ through $\calt$. The lemma
gives  conditions under which 
 there exist a height-$(n + 1)$ extension $\calT^+$
 of $\calT$ which is a good tower, and a
good
homotopy-lift $\tphi^+$ of $\phi$ through the tower $\calT^ +$, such that $\ADS(\tphi^+) < \ADS(\tphi)$.

Specifically, the lemma asserts that such a $\calT^+$ and $\tphi^+$
exist if
there is a
two-sheeted covering space  $\redPi_{n+1}:M_{n+1}\to N_n$ 
such that $\tphi:K\to
N_n$ admits a lift to  $M_{n+1}$,
and if one of the following (mutually exclusive) conditions is satisified:
\begin{enumerate}
\item [($\alpha$)] $n\ge1$, the manifold $N_n$ is closed 
(so that in the notation of Subsection \ref{from acs} we have $j_0=n$, and hence
$\Pi_{i+1}:M_{i+1}\to M_i$ is a covering map of closed manifolds for $0\le i\le n$), 
and the covering map
$\redPi_n\circ \redPi_{n+1} : M_{n+1}\to M_{n-1}$ 
is regular and has covering group isomorphic to $(\ZZ/2\ZZ)\times
(\ZZ/2\ZZ)$; or
\item [($\beta$)] $\partial M_{n+1}\ne\emptyset$; or
\item [($\gamma$)] $n=0$.
\end{enumerate}

(The stray phrase
``of rank $k\ge2$'' which appears in the hypothesis of \cite[Lemma
8.8]{acs-singular} seems to be the result of an editing error, and
should be ignored. It is not mentioned in the
proof or applications of Lemma 8.8 given in \cite{acs-singular}, and
in any case it could be deduced from the other hypotheses.)

Under the hypotheses of the present lemma, since $\tphi_*:H_1(K,\FF_2)\to H_1 (N_n;\FF_2)$ is not
 surjective,  there exists 
a
two-sheeted covering space  $\redPi^0_{n+1}:M_{n+1}^0\to N_n$ such that $\tphi:K\to
N_n$ admits a lift to  $M_{n+1}^0$. 
If Alternative (\wasA) of the hypothesis holds, then the two-sheeted
covering $M_{n+1}^0$ of $N_n$ has non-empty boundary, and hence
($\beta$) holds when we set $M_{n+1}=M_{n+1}^0$. Likewise, Alternative
(C) of the hypothesis implies ($\gamma$). Thus if (\wasA) or (C) holds, we
may apply Lemma 8.8 of \cite{acs-singular}, taking 
$M_{n+1}=M_{n+1}^0$ and $\redPi_{n+1}=\redPi_{n+1}^0$, to obtain the required
conclusion.

Now suppose that Alternative (\wasB) of the hypothesis holds but that (\wasA)
does not. In this case
we will show that there exists 
a
two-sheeted covering space  $\redPi^1_{n+1}:M_{n+1}^1\to N_n$ such that $\tphi:K\to
N_n$ admits a lift to  $M_{n+1}^1$, and such that ($\alpha$) holds
when we set $M_{n+1}=M_{n+1}^1$ and $\redPi_{n+1}=\redPi_{n+1}^1$. Applying Lemma 8.8 of \cite{acs-singular}, taking 
$M_{n+1}=M_{n+1}^1$ and $\redPi_{n+1}=\redPi_{n+1}^1$, will then give the required
conclusion in this case.

Since (\wasA) does not hold, 
$N_n$
is closed in this case. It then follows that in the notation of
Subsection \ref{from acs} above, we have $j_0=n$; hence
$\redPi_n:M_n\to M_{n-1}$ is a covering map of closed manifolds
.

We denote by $\tphi^-$ the \redPL map $\redPi_n\circ\tphi$ from $K$ to
$M_{n-1}=N_{n-1}$. We set $V=H_1(M_{n-1};\FF_2)$ and
$X=\tphi^-_*(H_1(K;\FF_2))\subset V$, we set $d=\dim
V=\rk_2(N_{n-1})$, and we let $q$
denote the natural homomorphism from $\pi_1(M_{n-1})$ to $V$. By
(\wasB) we have $d\ge\rk_2(K)+2\ge(\dim X)+2$; thus $\dim X\le d-2$.

The
two-sheeted cover $M_n$ of $M_{n-1}$ corresponds to a normal subgroup
of $\pi_1(M_{n-1})$ having the form $q^{-1}(Z)$, where $Z$ is some
codimension-$1$ vector subspace of $V$. Since $\tphi^-$ admits
the lift $\tphi$ to $M_n$, we have $X\subset Z\subset V$. Since in
addition we have $\dim X\le d-2<d-1=\dim Z$, there exists a
$(d-2)$-dimensional vector subspace $Y$ of $V$ with $X\subset Y\subset
Z$. The subgroup $q^{-1}(Y)$ determines a regular covering 
space $P:M^1_{n+1}\to M_{n-1}$ with covering group $(\ZZ/2\ZZ)\times
(\ZZ/2\ZZ)$. Since $q^{-1}(Y )\subset q^{-1}(Z)$, the degree-four
covering map $P$ factors as the composition of a degree-two covering
map $\redPi^1_{n+1}:M^1_{n+1}\to M_n$ with $\redPi_n$. Thus Condition
($\alpha$) holds with $M^1_{n+1}$ and
$\redPi_{n+1}^1$ playing the roles of 
$M_{n+1}$ and
$\redPi_{n+1}$. It remains to show that $\tphi$ admits a lift to $M^1_{n+1}$.

Since $\tphi^-_\sharp (\pi_1(K))\subset q^{-1}(X)\subset  q^{-1}(Y)$,
the map $\tphi^-$ admits a lift to the four-fold cover $M^1_{n+1}$ of
$M_{n-1}$. Since $M^1_{n+1}$ is a regular covering space of $M_{n-1}$,
there exist four different lifts of $\tphi^-$ to $M^1_{n+1}$. But
$\tphi^-$ can have at most two lifts to $M_n$, and each of these can
have at most two lifts to $M^1_{n+1}$. Hence each lift of $\tphi^-$
to $M_n$ admits a lift to $M^1_{n+1}$. In particular, $\tphi$  admits a lift to $M^1_{n+1}$.

(The last two paragraphs were adapted from the final step of the proof of \cite[Lemma 8.9]{acs-singular}.)
\EndProof

\Lemma\label{that's a good one}
Let $K$ be a compact, connected PL space such that $\pi_1(K)$ is
freely indecomposable, and let $\phi$ be a $\pi_1$-injective \redPL map from $K$
to a simple manifold  $Y$.  Assume that either 
(\wasb) $\rk_2(Y)\ge\rk_2(K)+3$
or 
(\wasa) $\partial
Y\ne\emptyset$. Then there exist a good
tower $\calt$ with base $Y$, and a 
good
homotopy-lift $\tphi:K\to N$ of $\phi$
(where $N$ denotes the top of the tower $\calt$), such that
$\tphi_*:H_1(K;\FF_2)\to H_1(N;\FF_2)$ is surjective.
\EndLemma

\Proof
Let $\scrS$ denote the set of all ordered pairs $(\redcalU, \redtpsi)$ such
that $\redcalU$ is a good tower with base $Y$ and $\redtpsi$ is a good
homotopy-lift of $\phi$ through $\redcalU$.
According to \cite[Proposition 5.8]{acs-singular}, there exist a PL
map $\tphi_0:K\to Y$, homotopic to $\phi$, and a compact, connected
three-dimensional submanifold $N_0$ of $\inter Y$, such that (i)
$\inter N_0\supset\tphi_0(K)$, (ii) the inclusion homomorphism
$\pi_1(\tphi_0(K))\to\pi_1(N_0)$ is surjective, (iii) $\partial N_0$ is incompressible in $Y$, and (iv) $N_0$ is irreducible. According to the definitions, this means that $\calT_0\doteq (M,N_0)$ is a good tower of height 0 and that $\tphi_0$ is a good homotopy-lift of $\phi$ through $\calT_0$; thus
$(\calT_0, \tphi_0)\in\scrS$, and in particular $\scrS\ne\emptyset$.
Hence there is an element $(\newredcalT,\redtphi)$ of $\scrS$ such that
$\ADS(\redtphi)\le\ADS(\redtpsi)$ for every element $(\redcalU,\redtpsi)$ of
$\scrS$. Let us write
$\newredcalT=(M_0,N_0,\redPi_1,...,N_n)$ (where $M_0=Y$),
so that $\redtphi$ is a PL map from $K$ to $N_n$. Our goal is to show that
$(\redtphi)_*:H_1(K;\FF_2)\to H_1(N_n;\FF_2)$ is surjective; this will
immediately imply the conclusion of the lemma.

As a preliminary to proving surjectivity, we claim:

\Claim\label{gold bug}
For each $j\in\{0,\ldots,n\}$, we have either 
(\wasB) $\rk_2(N_j)
\ge\rk_2(K)+2$
or 
(\wasA) $\partial
N_j\ne\emptyset$.
\EndClaim

To prove \ref{gold bug}, we first
consider the case in which Alternative (\wasb) of the hypothesis of
the present 
lemma holds, i.e. $\rk_2(M_0)=\rk_2(Y)\ge\rk_2(K)+3$. In particular we then have
$\rk_2(M_0)\ge3$. We  invoke Lemma 8.5 of
\cite{acs-singular}, which implies that if 
$(M_0,N_0,\redPi_1,...,N_n)$ is any good tower of $3$-manifolds, and if
$r\doteq \rk_2(M_0)\ge3$, then for any index $j\in\{0,\ldots,n\}$ such that
$M_j$ is closed, we have $\rk_2(M_j)\ge r-1$. If Alternative (\wasA) of
\ref{gold bug} does not hold for a given $j$, i.e. if $N_j$ is closed,
then since $N_j$ is a submanifold of the connected manifold $M_j$, we
must have $N_j=M_j$, and in particular $M_j$ is closed. Hence
$\rk_2(N_j)=\rk_2(M_j)\ge r-1$. Since $r=\rk_2(Y)\ge\rk_2(K)+3$ in this
case, it follows that Alternative (\wasB) of \ref{gold bug} holds.

Now
consider the case in which
Alternative (\wasa) of the hypothesis of the lemma holds, i.e. $\partial
M_0=\partial Y\ne\emptyset$. In this case, in the notation of
Subsection \ref{from acs}, we have $j_0=0$; hence for $j=0,\ldots,n$
we have
$\partial
N_j\ne\emptyset$, i.e. Alternative (\wasA) of \ref{gold bug} holds. 
This
completes the proof of \ref{gold bug}.

We now turn to the proof that 
$(\redtphi)_*:H_1(K;\FF_2)\to H_1(N_n;\FF_2)$ is surjective. We apply Lemma
\ref{new eight point nine}, 
taking $\newredcalT$ and $\redtphi$ 
to be defined as above.
Since
$(\newredcalT,\redtphi)\in\cals$, the map 
$\redtphi$ 
is a good
homotopy-lift of $\phi$ through the good tower $\newredcalT$, and in
particular is a \redPL map from $K$ to $N\doteq N_n$. If $n=0$ then Alternative
(C) of the hypothesis of Lemma \ref{new eight point nine}
holds; and if $n\ge1$ then the case $j=n$ of \ref{gold
  bug} 
asserts that one of the alternatives (\wasB), (\wasA) of Lemma \ref{new eight point nine} holds. 
If we assume that 
$(\redtphi)_*$ is not surjective, it now follows from  Lemma \ref{new
  eight point nine} that there exist a height-$(n + 1)$ extension $\calT^+$
 of $\newredcalT$ which is a good tower, and a good
homotopy-lift $\tphi^+$ of $\phi$ through the tower $\calT^ +$, such
that $\ADS(\tphi^+) < \ADS(\redtphi)$. According to the definition we
have $(\calt^+,\tphi^+)\in\scrS$, and the inequality $\ADS(\tphi^+) <
\ADS(\redtphi)$ contradicts our choice of $(\newredcalT,\redtphi)$. The
 surjectivity of $(\redtphi)_*$ is thereby established.
\EndProof

\Definition\label{carried def}
Let $Z$ be a path-connected subset of a path-connected topological
space $Y$. We shall say that a subgroup $A$ of $\pi_1(Y)$ is {\it
  carried by} $Z$ if $A$ is conjugate to a subgroup of the image $P$ of
the inclusion homomorphism $\pi_1(Z)\to\pi_1(Y)$. (Here we regard the
subgroups $A$ and $P$ as being themselves defined only up to conjugacy.)
\EndDefinition

The following theorem is the main result of this paper. It was stated in a weaker form, in the introduction,
as \redTheoremC.

\Theorem\label{main desingular}
Let $Y$ be a simple 
$3$-manifold, 
and let $A$ be a finitely generated, freely indecomposable
subgroup of $\pi_1(Y)$. Set $\eta=\rk_2(A)$.  
Suppose that {\it either} (a) $\partial
Y\ne\emptyset$ {\it or} (b) 
$\rk_2(Y)\ge3\eta^2-4\eta+4$.
Then 
$A$ is carried
by some
compact,
connected,
three-dimensional
submanifold
$Z$ of $\inter Y$ 
such that 
\begin{itemize}
\item
$\partial Z$ is  non-empty and incompressible in $Y$;
\item
$0\le\chibar(Z)\le\eta-1$; and
\item
$\rk_2(Z)\le 3\eta^2-4\eta-3\chibar(Z)^2 +3\chibar(Z)+1$; in particular,
$\rk_2(Z)\le 
3\eta^2-4\eta+1$.
\end{itemize}
\EndTheorem

\Proof
Since $A$ is a finitely generated subgroup of the $3$-manifold group
$\pi_1(Y)$, it follows from the main theorem of \cite{scott} that $A$
is finitely presented. Hence $A$ is isomorphic to $\pi_1(K)$ for some
compact, connected PL
 space $K$ of dimension at most $2$. We fix a
$\pi_1$-injective \redPL map   $\phi:K\to Y$ such that the image of
$\phi_\sharp:\pi_1(K)\to\pi_1(Y)$ is conjugate to $A$. 
In
particular  $\pi_1(K)$ is isomorphic to $A$ and is therefore freely indecomposable.

We shall show that there exist a \redPL map $g:K\to Y$, homotopic to $\phi$, and a compact,
connected 
submanifold
$Z$ of $\inter Y$ with $Z\supset g(K)$, such that $Z$ satisfies the bulleted
conditions in the statement of the theorem; this will
immediately imply the conclusion.

The hypothesis implies  that either $\partial
Y\ne\emptyset$ or  $\rk_2(Y)\ge\rk_2(K)+3$.
Hence Lemma \ref{that's a
  good one} gives a good
tower 
\Equation\label{whitehall}
\calt=(M_0,N_0,\redPi_1,M_1,N_1,\redPi_2,\ldots,\redPi_n,M_n,N_n) 
\EndEquation
with base $M_0=Y$, and a 
good
homotopy-lift $\tphi:K\to N_n$ of
$\phi$, such that
$\tphi_*:H_1(K;\FF_2)\to H_1(N_n;\FF_2)$ is surjective.

Proposition \ref{perfect is good} gives:
\Claim\label{ergo simplisticus}
For $j=0,\ldots,n$, the manifolds $M_j$ and $N_j$ are simple.
\EndClaim

For $j =
0,\ldots,n$ we let
$\iota_j : N_j \to M_j$ denote the inclusion map.

For $\reallyredj=0,\ldots,n-1$ we define a \redPL map $\phi_\reallyredj:K\to N_\reallyredj$ by
$\phi_\reallyredj=\redPi_{\reallyredj+1}\circ\iota_{\reallyredj+1}\circ \redPi_{\reallyredj+2}\circ\cdots\circ\iota_{n-1}\circ \redPi_n\circ\iota_n\circ\tphi$. We set $\phi_n=\tphi$. Note that
$\phi_0$, regarded as a map from $K$ to $Y=M_0\supset N_0$, is
homotopic to $\phi$.

We claim:

\Claim\label{new higgledy}
For any index $\redj\in\{0,\ldots,n\}$, and for any connected submanifold $X$ of
$M_\redj$ with $\rk_2(X)\le3\eta^2-4\eta+2$, we have $\partial X\ne\emptyset$.
\EndClaim

To prove \ref{new higgledy}, first note that the assertion is
immediate if $\partial M_\redj\ne\emptyset$. Next note that if $\partial
M_0\ne\emptyset$, 
then in the notation of \ref{from acs} we have $j_0=0$, and hence
$\partial M_\redj\ne\emptyset$ for every $\redj$.


Now suppose that $M_0=Y$ is closed. 
In this case we apply   \cite[Lemma
8.5]{acs-singular}, which implies that
if $\calt=(M_0,N_0,\redPi_1,M_1,N_1,\redPi_2,\ldots,\redPi_n,M_n,N_n)$
is any good tower of $3$-manifolds such that $\rk_2(M_0)\ge3$, then
for any index $j\in\{0,\ldots,n\}$ we have either 
$\partial M_\redj\ne\emptyset$ or
$\rk_2(M_\redj)\ge\rk_2(M_0)-1$.
In the present situation, since $Y$ is closed, the hypothesis of the
theorem gives
$\rk_2(M_0)=\rk_2(Y)
\ge
3\eta^2-4\eta+4$. Since
$\eta$ 
is an integer,
we indeed have $\rk_2(M_0)\ge3$, 
as required for the application of \cite[Lemma
8.5]{acs-singular}.
Hence  for a given index $\redj\in\{0,\ldots,n\}$,
we have either (I) $\partial M_\redj\ne\emptyset$, so that the conclusion of
\ref{new higgledy} is true, or (II)
$\rk_2(M_\redj)\ge\rk_2(M_0)-1=\rk_2(Y)-1 \ge
3\eta^2-4\eta+3$.
Since the hypothesis of 
\ref{new higgledy} gives $\rk_2(X)\le3\eta^2-4\eta+2$, we have
$X\ne M_\redj$ in this case; as $M_\redj$ is 
simple by \ref{ergo simplisticus} and therefore connected,
it follows that
$\partial X\ne\emptyset$. This completes the proof of \ref{new higgledy}.

Next we will prove, by induction on $\redi=0,\ldots,n$:
\Claim\label{here's the thing}
There exist a connected submanifold $R_{n-\redi}$ of
  $\inter N_{n-\redi}$
and a \redPL map
$\redgamma_{n-\redi}:
K\to N_{n-\redi}$ with the following
properties:
\begin{enumerate}[(a)]
\item $\partial R_{n-\redi}$ is non-empty, and is incompressible in $N_{n-\redi}$ (so that $R_{n-\redi}$ is
  simple by \ref{ergo simplisticus} and Remark
\ref{more fruit});
\item $\chibar(R_{n-\redi})\le\eta-1$;
\item $\rk_2(C_{R_{n-\redi}})\le 3\eta^2-4\eta-3\chibar(R_{n-\redi})^2+1$;
\item $\redgamma_{n-\redi}(K)\subset R_{n-\redi}$; and
\item $\redgamma_{n-\redi}:K\to N_{n-\redi}$ 
  is homotopic to the map 
$\phi_{n-\redi}$.
\end{enumerate}
\EndClaim

To prove \ref{here's the thing}, we first consider the base case
$\redi=0$. In this case, we take $R_n$ to be a submanifold obtained from $M_n$
by removing a half-open boundary collar. Then
$\phi_n=\tphi$
is homotopic in $N_{n}$ to a \redPL map whose image is contained in
$R_{n}$. Conditions (d) and (e) of \ref{here's the thing} are then
immediate. Since $\tphi_*:H_1(K;\FF_2)\to H_1(N_n;\FF_2)$ is surjective,
we have 
$\rk_2(R_n)=\rk_2(N_n)\le\rk_2(K)=\eta$. 
Since 
the inequality $\eta\le3\eta^2-4\eta+2$ holds for every integer $\eta$, we
have $\rk_2(N_n)\le3\eta^2-4\eta+2$.
Now since the manifold $N_n$ is simple by \ref{ergo
  simplisticus}, it is in particular connected, and it therefore
follows from
\ref{new higgledy}  that $\partial N_n\ne\emptyset$.
On the other hand, according to
\ref{keep it simple}, the simplicity of $N_n$ also implies that $N_n$ is
boundary-irreducible and has no $2$-sphere boundary component. 
The
definition of $R_n$ then implies that $\partial R_n$ is incompressible in
$N_n$. Thus Condition (a) of  \ref{here's the thing} holds. 

Since $\partial R_n\ne\emptyset$ we have $H_3(R_n;\FF_2)=0$ and
connectedness gives $\dim H_0(R_n;\FF_2)=1$; hence
$\chibar(R_n)\le\rk_2(R_n)-1\le\eta-1$. This is Condition
(b). To verify Condition (c) we note that, since the 
manifold $R_n$ has non-empty boundary and is simple by Condition (a),  Lemma \ref{here again} gives
$\rk_2(C_{R_n})\le \rk_2(R_n)+\chibar(R_n)-1\le\eta+(\eta-1)-1=2\eta-2$.
Since  $3\eta^2-4\eta-3\chibar(R_{n-\redi})^2+1\ge
3\eta^2-4\eta-3(\eta-1)^2+1=2\eta-2$, Condition (c) now follows, and the
base case is established.

For the induction step, we suppose that an index
$\redi\in\{0,\ldots,n-1\}$ is given and that there exist a connected
submanifold $R_{n-\redi}$ of $\inter N_{n-\redi}\subset\inter M_{n-\redi}$ and a \redPL map 
$\redgamma_{n-\redi}
:K\to N_{n-\redi}$
such that Conditions (a)--(e) of
  \ref{here's the thing} hold.

We recall that $\redPi_{n-\redi}:M_{n-\redi}\to N_{n-(\redi+1)}$ is a two-sheeted
covering map
by the definition of a good tower, and that $N_{n-(\redi+1)}$ is
simple by \ref{ergo simplisticus}.
We also observe that since $\partial N_{n-\redi}$ is
incompressible in $M_{n-\redi}$ by 
the definition of a good tower,
the submanifold $N_{n-\redi}$ is $\pi_1$-injective in $M_{n-\redi}$
  by \ref{don't squeeze-a da fruit}.
Since $\partial R_{n-\redi}$  is
incompressible in $N_{n-\redi}$  by
Condition (a) of 
the induction hypothesis,
it follows 
that $\partial R_{n-\redi}$ is incompressible in
$M_{n-\redi}$. Condition (a) of
the induction hypothesis
also gives $\partial R_{n-\redi}\ne\emptyset$.
Note in addition that since, by Condition (e) of 
the induction hypothesis
\ref{here's the thing},
$\gamma_{n-\redi}$ is in particular homotopic to $\phi_{n-\redi}$ when regarded as a
  map of $K$ into $N_{n-\redi}$, and since $\redPi_1\circ\cdots\circ
  \redPi_{n-\redi}\circ\phi_{n-\redi}$ is homotopic to the $\pi_1$-injective map $\phi$
  by the definition of a homotopy-lift, the map $\gamma_{n-\redi}:K\to M_{n-\redi}$ is
  itself $\pi_1$-injective.
Thus
all the hypotheses of Lemma \ref{key lemma}  hold
 with $N_{n-(\redi+1)}$,
$M_{n-\redi}$, $\redPi_{n-\redi}$, $R_{n-\redi}$ and
$\redgamma_{n-\redi}$
playing
the respective roles of $Q$, $\tQ$, $\redPi$, $R$ and $f$, and with $K$
defined as above.
Lemma \ref{key lemma} therefore gives a 
simple
submanifold of $\inter N_{n-(\redi+1)}$,
which is denoted by $T$ in that lemma and will be denoted by
$R_{n-(\redi+1)}$ in the present context;
the boundary of
$R_{n-(\redi+1)}$ is incompressible in  $N_{n-(\redi+1)}$, 
and one of the following alternatives holds:
\Alternatives
\item
  $\chibar(R_{n-(\redi+1)})=\chibar(R_{n-\redi})$,
  and the pairs $(C_{R_{n-(\redi+1)}},\redcale_{R_{n-(\redi+1)}})$
and $(C_{R_{n-\redi}},\redcale_{R_{n-\redi}})$ are homeomorphic, 
or 
\item
  $\chibar(R_{n-(\redi+1)})<\chibar(R_{n-\redi})$, and 
  $$\rk_2(R_{n-(\redi+1)})\le \rk_2(C_{R_{n-\redi}})
  +4\chibar(R_{n-\redi})+\chibar(R_{n-(\redi+1)}).$$
\EndAlternatives

The conclusion of Lemma \ref{key lemma} also gives a \redPL map from $K$ to
$N_{n-(\redi+1)}$, which we shall denote here by $\redgamma_{n-(\redi+1)}$; it
has image contained in $R_{n-(\redi+1)}$ and (as a
map from $K$ to
$N_{n-(\redi+1)}$) is
homotopic to $\redPi_{n-(\redi+1)}\circ\redgamma_{n-\redi}$. Since
  $\redgamma_{n-\redi}:K\to N_{n-\redi}$ is homotopic to
$\phi_{n-\redi}$,
the map $\redgamma_{n-(\redi+1)} :K \to
N_{n-(\redi+1)}$ is homotopic to
$\phi_{n-(\redi+1)}=\redPi_{n-(\redi+1)}\circ\phi_{n-\redi}$.
Thus Conditions (d) and (e) of \ref{here's the thing} hold when $\redi$ is
replaced by $\redi+1$. To complete the induction step, we must show that
Conditions (a)--(c) of \ref{here's the thing} also persist when $\redi$ is
replaced by $\redi+1$. 
This will be done in two cases,
depending on which of the alternatives (i), (ii) holds.

First suppose that (i) holds. According to 
the remark following the
statement of Lemma \ref{key lemma}, 
the incompressible surface $\partial R_{n-(\redi+1)}$ is
non-empty, so that Condition (a) persists in this case. Furthermore,
in this case we have $\chibar(R_{n-(\redi+1)})=\chibar(R_{n-\redi})$ and
$\rk_2(C_{R_{n-(\redi+1)}})=\rk_2(C_{R_{n-\redi}})$, so that the inequalities
$\chibar(R_{n-\redi})\le\eta-1$ and  $\rk_2(C_{R_{n-\redi}})\le 3\eta^2-4\eta-3\chibar(R_{n-\redi})^2+1$
  translate directly into
  $\chibar(R_{n-(\redi+1)})\le\eta-1$ and $\rk_2(C_{R_{n-(\redi+1)}})\le
  3\eta^2-4\eta-3\chibar(R_{n-(\redi+1)})^2+1$. Thus Conditions (b) and
  (c) also persist in this case.

Now consider the case in which (ii) holds. In this case we set $\alpha =\chibar(R_{n-\redi})$ and
$\beta =\chibar(R_{n-(\redi+1)})$.  By \ref{keep it simple} we have
$\beta\ge0$; combining this with the statement of Alternative (ii), we
find
\Equation\label{what it says one}
0\le \beta \le \alpha -1
\EndEquation
and
\Equation\label{what it says two}
\rk_2(R_{n-(\redi+1)})\le
\rk_2(C_{R_{n-\redi}})+4\alpha+\beta.
\EndEquation

Since Condition (b) of the induction hypothesis gives
$\alpha\le\eta-1$, it follows from (\ref{what it says one}) that
$\beta<\eta-1$; thus 
(b)
persists. Next note that 
Condition (c) of the induction hypothesis gives $\rk_2(C_{R_{n-\redi}})\le
3\eta^2-4\eta-3\alpha^2+1$; combining this with (\ref{what it says two}), we obtain
\Equation\label{efharisto}
\rk_2(R_{n-(\redi+1)})\le 3\eta^2-4\eta-3\alpha^2 +4\alpha+\beta+1.
\EndEquation

From (\ref{efharisto}) and (\ref{what it says one}) we deduce that
$\rk_2(R_{n-(\redi+1)})\le 3\eta^2-4\eta-3\alpha^2 +5\alpha$. Since
(\ref{what it says one}) also gives $\alpha\ge1$, we have $3\alpha^2
-5\alpha\ge-2$, and hence $\rk_2(R_{n-(\redi+1)})\le 3\eta^2-4\eta+2$.
Since the simple manifold $R_{n-(\redi+1)}$ is in particular
connected, it
now follows from  \ref{new higgledy}
that the incompressible surface $\partial
R_{n-(\redi+1)}$ is non-empty; this shows that Condition (a) persists.

Since  the simple $3$-manifold
$R_{n-(\redi+1)}$ has non-empty boundary, it follows from Lemma \ref{here again} that
$\rk_2(C_{R_{n-(\redi+1)}})\le \rk_2(R_{n-(\redi+1)})+\beta-1$. Combining this
with (\ref{efharisto}), we obtain
\Equation\label{mcgintloch}
\rk_2(C_{R_{n-(\redi+1)}})\le 3\eta^2-4\eta-3\alpha^2 +4\alpha+2\beta.
\EndEquation


The function $x\mapsto3x^2-4x$ increases monotonically on
$[1,\infty]$. Since by (\ref{what it says one}) we have
$1\le\beta+1\le\alpha$, we have $3\alpha^2
-4\alpha\ge3(\beta+1)^2-4(\beta+1)=3\beta^2+2\beta-1$. It therefore
follows from (\ref{mcgintloch}) that
$$
\rk_2(C_{R_{n-(\redi+1)}})\le 3\eta^2-4\eta-3\beta^2 +1.
$$
This shows that Condition (c) of \ref{here's the thing} holds when $\redi$
is replaced by $\redi+1$. Thus the induction is complete, and \ref{here's
  the thing} is proved.


To complete the proof of the theorem, we apply the case $\redi=n$ of
\ref{here's the thing}. We set $Z=R_0$, and we define $g$ to be the
map $\gamma_0$ regarded as a \redPL map from $K$ to $Y$.

The manifold $N_0$ is $\pi_1$-injective in $M_0$ because it has an
incompressible boundary
by the definition of a good tower.
Since $\partial Z$ is
incompressible in $N_0$ by Condition (a) of \ref{here's the thing}, it
follows that $\partial Z$ is
incompressible in $M_0$. Condition (a) also assures that $\partial Z$
is non-empty, and that $Z$ is simple. Condition (b) asserts that
$\chibar(Z)\le\eta-1$. By \ref{keep it simple} we also have $\chibar(Z)\ge0$.

Condition (c) of
\ref{here's the thing} asserts that $\rk_2(C_{Z})\le
3\eta^2-4\eta-3\chibar(Z)^2+1$. But by 
the second assertion of
Corollary \ref{simple no cylinders} we
have $\rk_2(Z)
\le
\rk_2(C_Z)+3\chibar(Z)$.
Hence 
\Equation\label{nice}
\rk_2(Z)\le 3\eta^2-4\eta-3\chibar(Z)^2 +3\chibar(Z)+1.
\EndEquation
Since $\chibar(Z)\ge0$, it follows from
(\ref{nice}) that $\rk_2(Z)\le 3\eta^2-4\eta+1$. 

Condition (d) guarantees that $g(K)\subset Z$. Finally,
since Condition (e) of \ref{here's the thing} gives that  $\gamma_0$, regarded as a map from $K$ to $N_0$, is homotopic to
$\phi_0$, in
particular $g$ and $\phi$ are homotopic as maps from $K$ to $M_0=Y$. 
Thus $Z$ and $g$ have all the properties stated in the conclusion of
the theorem.

\EndProof

\bibliographystyle{plain}

\end{document}